\documentclass[preprint,12pt,review,numbers,sort&compress]{elsarticle}

\usepackage{amssymb}
\usepackage{subfigure}
\usepackage{caption}
\usepackage{multirow}
\usepackage{tgtermes}
\usepackage{amsmath}
\usepackage{placeins}
\usepackage{natbib}
\usepackage{array}
\usepackage{graphicx}
\usepackage{gensymb}
\usepackage{textcomp}
\usepackage{bm}
\usepackage{enumerate}
\usepackage[margin=2.5cm]{geometry}
\usepackage{color, soul}
\usepackage[justification=centering]{caption}
\usepackage{booktabs}
\usepackage{mathtools}

\allowdisplaybreaks[4]

\journal{ }

\begin{document}

\begin{frontmatter}
\sethlcolor{yellow}

\title{A new meshless Fragile Points Method (FPM) with minimum unknowns at each point, for flexoelectric analysis under two theories with crack propagation. Part II: Validation and discussion}

\author[add1]{Yue Guan\corref{cor1}}
\ead{yuguan@ttu.edu}
\cortext[cor1]{Corresponding author.}

\author[add2]{Leiting Dong}
\author[add1]{Satya N. Atluri}

\address[add1]{Department of Mechanical Engineering, Texas Tech University, Lubbock, TX 79415, United States}
\address[add2]{School of Aeronautic Science and Engineering, Beihang University, Beijing 100191, China}

\begin{abstract}

In the first part of this two-paper series, a new Fragile Points Method (FPM), in both primal and mixed formulations, is presented for analyzing flexoelectric effects in 2D dielectric materials. In the present paper, a number of numerical results are provided as validations, including linear and quadratic patch tests, flexoelectric effects in continuous domains, and analyses of stationary cracks in dielectric materials. A discussion of the influence of the electroelastic stress is also given, showing that Maxwell stress could be significant and thus the full flexoelectric theory is recommended to be employed for nano-scale structures. The present primal as well as mixed FPMs also show their suitability and effectiveness in simulating crack initiation and propagation with flexoelectric effect. Flexoelectricity, coupled with piezoelectric effect, can help, hinder, or deflect the crack propagation paths and should not be neglected in nano-scale crack analysis. In FPM, no remeshing or trial function enhancement are required in modeling crack propagation. A new Bonding-Energy-Rate(BER)-based crack criterion as well as classic stress-based criterion are used for crack development simulations. Other complex problems such as dynamic crack developments, fracture, fragmentation and 3D flexoelectric analyses will be given in our future studies.

\end{abstract}

\end{frontmatter}

\section{Introduction}

As a result of the recent trend of miniaturization of electromechanical systems, there has been an increasing demand for reliable and accurate theories and numerical methods for flexoelectric analysis \cite{Wang2019, Zhuang2020, Mao2014}. Part I of this study introduced a new meshless Fragile Points Method (FPM) \cite{Dong2019, Guan2020} based on a Galerkin weak form and local, simple, polynomial and discontinuous trial and test functions, and presented the theoretical foundation and numerical implementation for the FPM in analyzing flexoelectric effects in dielectric solids. Both primal and mixed FPM formulations were developed, based on two flexoelectric theories with or without the electric field gradient effect and electroelastic stress. In general, the FPM presents clear advantages in its algorithmic formulation as compared to previous numerical methods including primal and mixed FEM \cite{Beheshti2017, Amanatidou2002}, Element-Free Galerkin (EFG) method \cite{He2019} and primal and mixed Meshless Local Petrov-Galerkin (MLPG) method \cite{Sladek2020, Atluri2004} in analyzing flexoelectric behavior, as it does not require high quality meshes, has minimum number of DoFs at each Point, needs only very simple weak form integration schemes, and has great advantages in simulating crack and rupture initiation and propagation.

The present Part II of the paper presents extensive numerical studies of the primal and mixed FPM based on both full as well as  reduced theories for analyzing flexoelectric effects at multiple length scales. Linear and Quadratic (displacement) patch tests are presented in section~\ref{sec:patch}. Section~\ref{sec:continuous} shows a number of numerical examples in 2D continuous domains to illustrate the implementation and accuracy of both the primal and mixed FPMs. A discussion on the influence of the electroelastic stress is also given. Section~\ref{sec:crack_1} concentrates on the analysis of flexoelectric effects on stationary cracks. Following that, the simulation of crack initiation and propagation is exhibited in section~\ref{sec:crack_2}. A short discussion of the influences of computational parameters is given at last in section~\ref{sec:comp_para}.

\section{Patch tests} \label{sec:patch}

\subsection{Linear Displacement patch tests}

First of all, we examine the consistency of the present FPM by conducting two patch tests. A unit square domain is considered. First, a known linear displacement and its corresponding stress field are given as:
\begin{align} \label{eqn:exact_linear}
\begin{split}
& u_1 =  x - y , \\
& u_2 =  x + y.
\end{split} \quad
\begin{split}
&\sigma_{11} = \sigma_{22}  =  \frac{E}{1- \nu} , \\
& \sigma_{12} = 0,
\end{split}
\end{align}
where $E$ and $\nu$ are the Young's modulus and Poisson's ratio of the material. Two kinds of boundary conditions are considered: 1). Dirichlet boundary conditions are applied on all edges of the square; 2). the left and bottom edges of the square are subjected to the known normal displacements, while constant tensile loadings are applied on the right and top edges with $\widetilde{Q} = E/(1-\nu)$ (as shown in Fig.~\ref{fig:Ex00_Schem_01}). No electrical loading is applied.

Here we define the relative error $e \left( \mathbf{x} \right) $ for any variable $\mathbf{x}$ as:
\begin{align} \label{eqn:error}
\begin{split}
e \left( \mathbf{x} \right) = \frac{ \left\| \mathbf{x}^h - \mathbf{x} \right\|_{L^2}}{\left\| \mathbf{x} \right\|_{L^2}}, \; \text{where} \; \left\| \mathbf{x} \right\|_{L^2} = \left( \int_\Omega \mathbf{x}^\mathrm{T} \mathbf{x} \mathrm{d} \Omega \right)^{1/2},
\end{split}
\end{align}
where $\mathbf{x}$ is the exact solution, and $\mathbf{x}^h$ is the computed result. 

Three different Fragile Point distributions are used in the FPM (shown in Fig.~\ref{fig:Ex00_Part}). For all the three kinds of point distributions, the linear displacement and stress field are reproduced successfully, with a relative error smaller than $2 \times 10^{-7}$.

\subsection{Quadratic displacement patch tests}

\begin{figure}[htbp] 
  \centering 
  \subfigure[]{ 
    \label{fig:Ex00_Schem_01} 
    \includegraphics[width=0.48\textwidth]{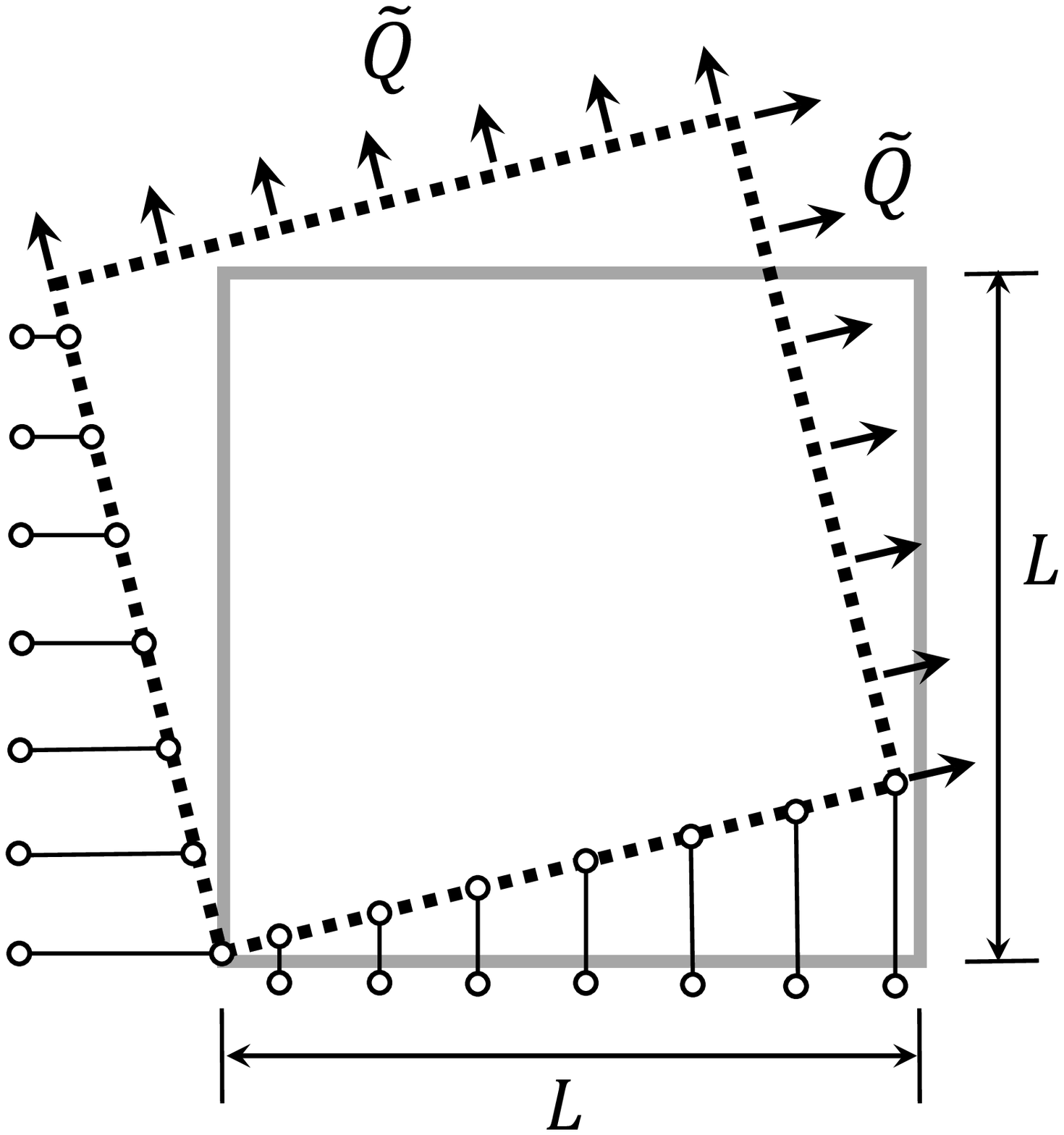}}  
    \subfigure[]{ 
    \label{fig:Ex00_Schem_02} 
    \includegraphics[width=0.48\textwidth]{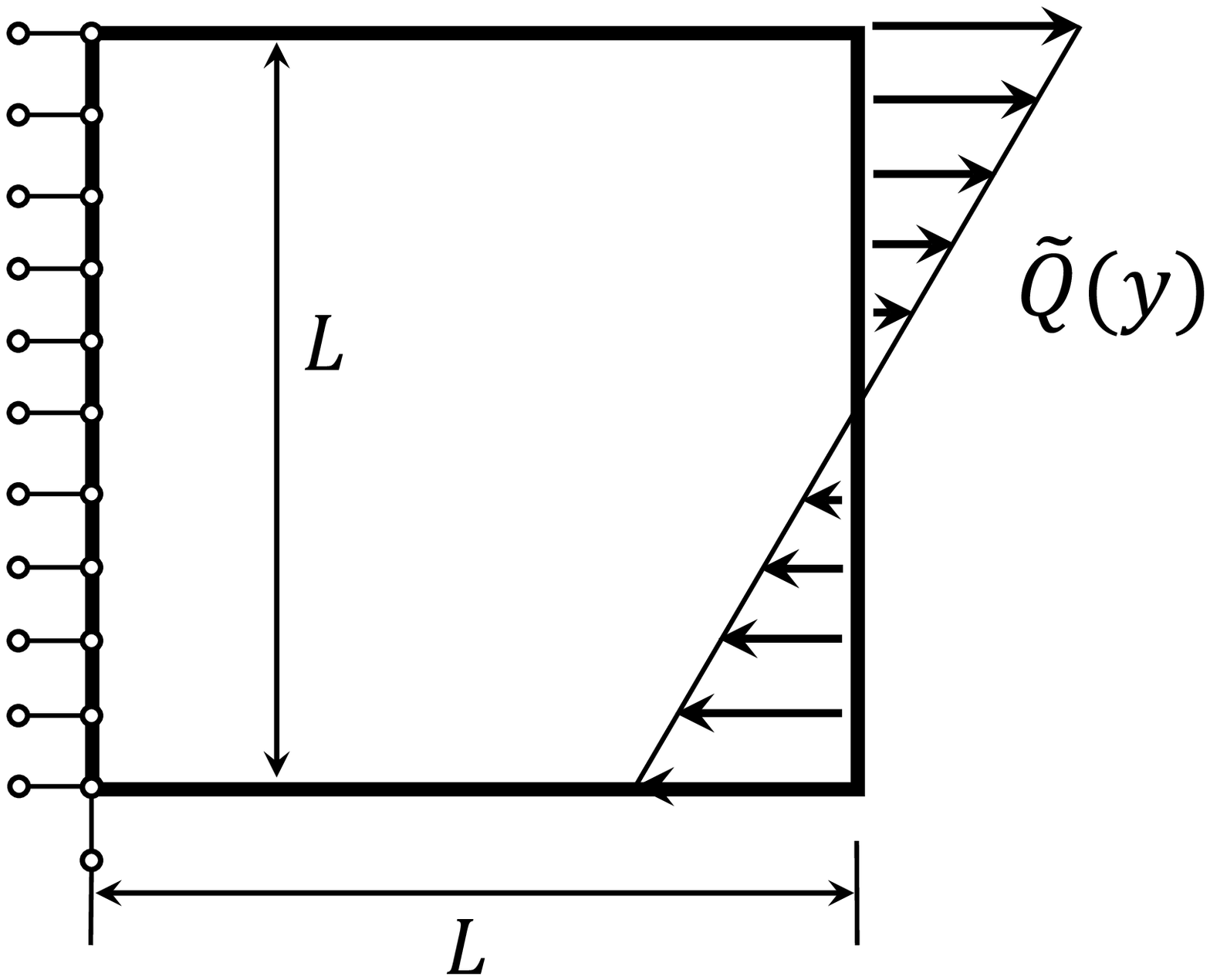}}  
  \caption{Patch test problems: a unit square domain subjected to (a) pure tensile loading. (b) a bending moment.} 
  \label{fig:Ex00_Schem} 
\end{figure}

Next, we consider a quadratic patch test proposed by \citet{Zienkiewicz1997} and \citet{Lee1993}. As shown in Fig.~\ref{fig:Ex00_Schem_02}, the plate is subjected to a constant bending moment. The left edge is constrained in $x$-direction. And the lower left vertex is fixed additionally in $y$-direction to avoid rigid body motions. The right edge has Neumann boundary conditions with a linear surface traction $\widetilde{Q} = E \left( 2 y - 1 \right)$. The exact analytical solution for the displacement and stress field can be given as:
\begin{align} 
\begin{split}
& u_1 =  2 xy - x , \\
& u_2 = - x^2 - \nu y^2 +  \nu y,
\end{split} \quad
\begin{split}
& \sigma_{11} = E \left( 2 y - 1 \right), \\
& \sigma_{22} = \sigma_{12} = 0.
\end{split}
\end{align}
The patch test describes a quadratic behavior in the displacement field.

\begin{figure}[htbp] 
  \centering 
    \subfigure[]{ 
    \label{fig:Ex00_Part_01} 
    \includegraphics[width=0.32\textwidth]{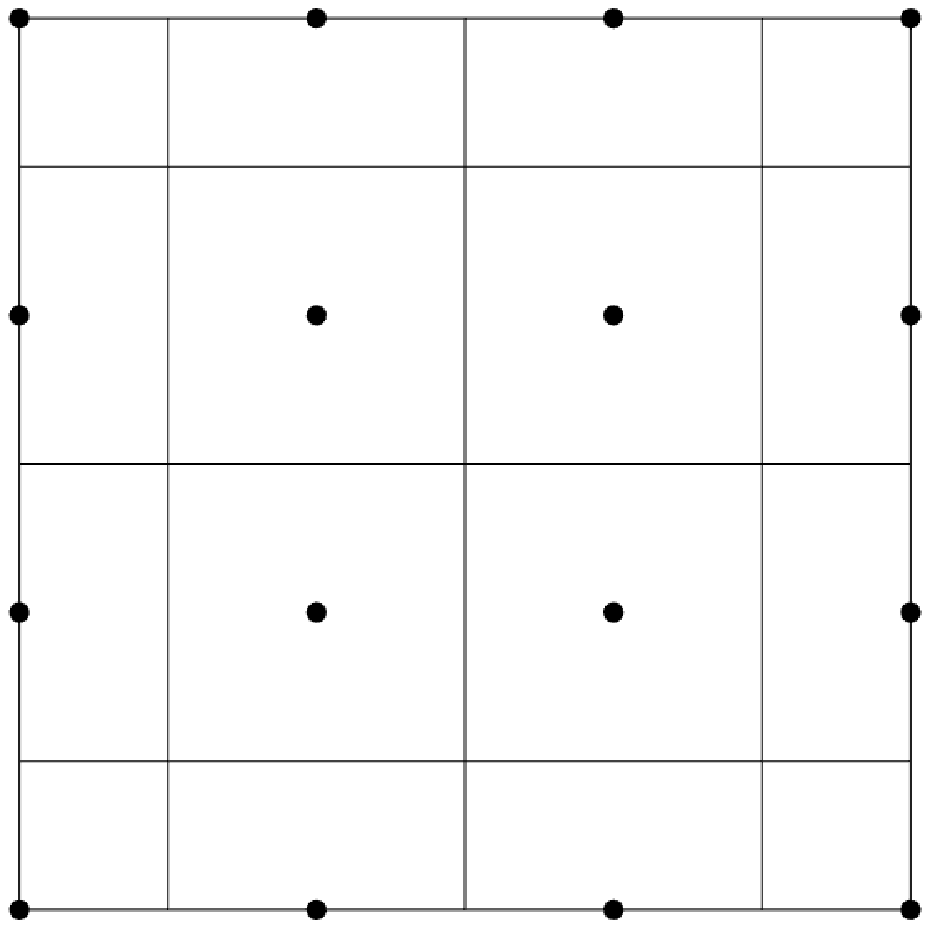}}  
    \subfigure[]{ 
    \label{fig:Ex00_Part_02} 
    \includegraphics[width=0.32\textwidth]{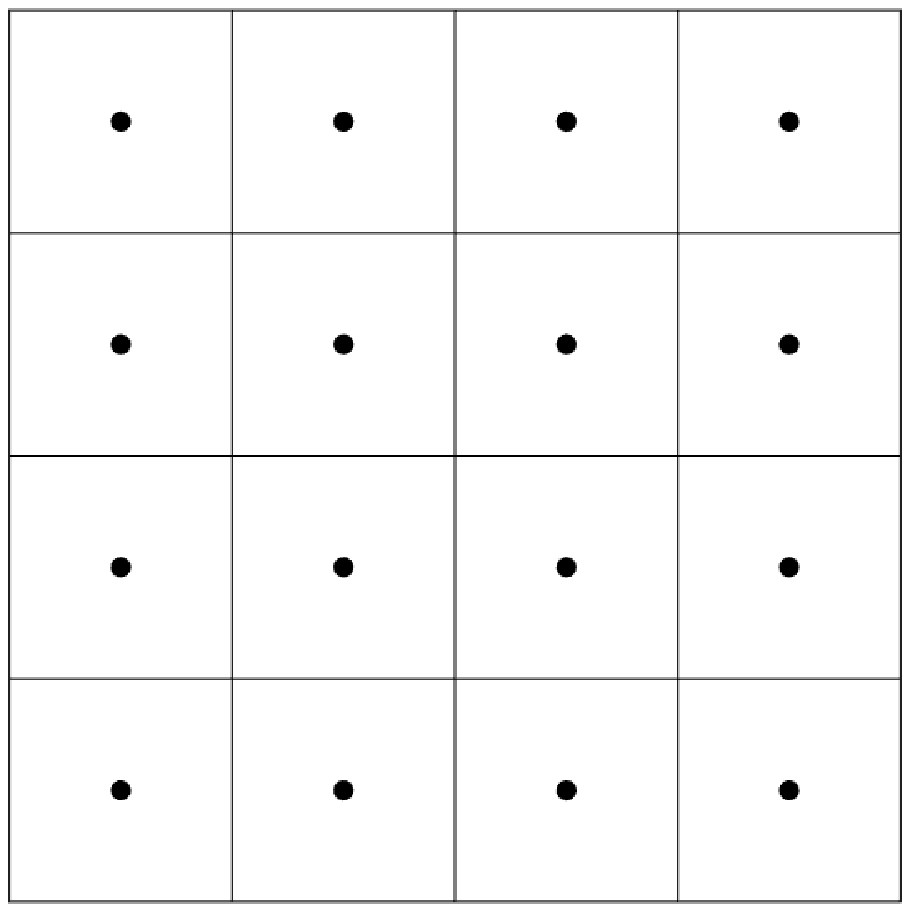}}  
     \subfigure[]{ 
    \label{fig:Ex00_Part_03} 
    \includegraphics[width=0.32\textwidth]{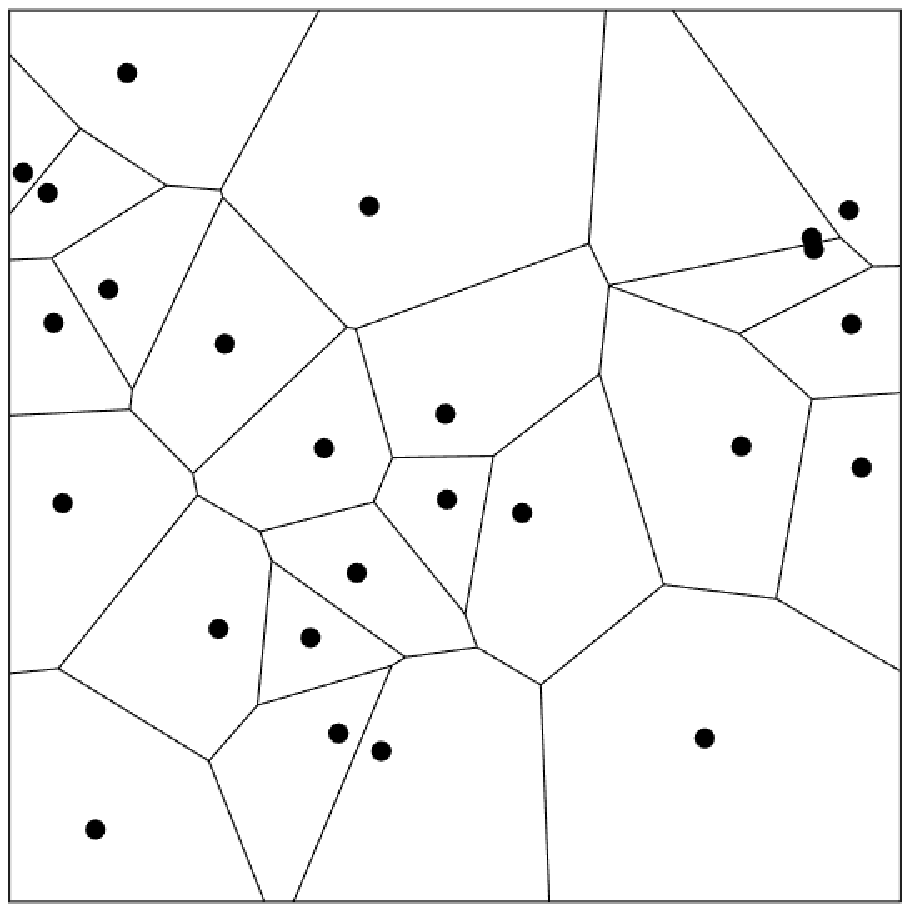}} 
  \caption{The point distribution and domain partition in the FPM. (a) 16 uniform points (4 points in the domain and 12 points on the boundaries). (b) 16 uniform points in the domain. (c) 25 random points.} 
  \label{fig:Ex00_Part} 
\end{figure}

\begin{figure}[htbp] 
  \centering 
    \subfigure[]{ 
    \label{fig:Ex00_01_01} 
    \includegraphics[width=0.32\textwidth]{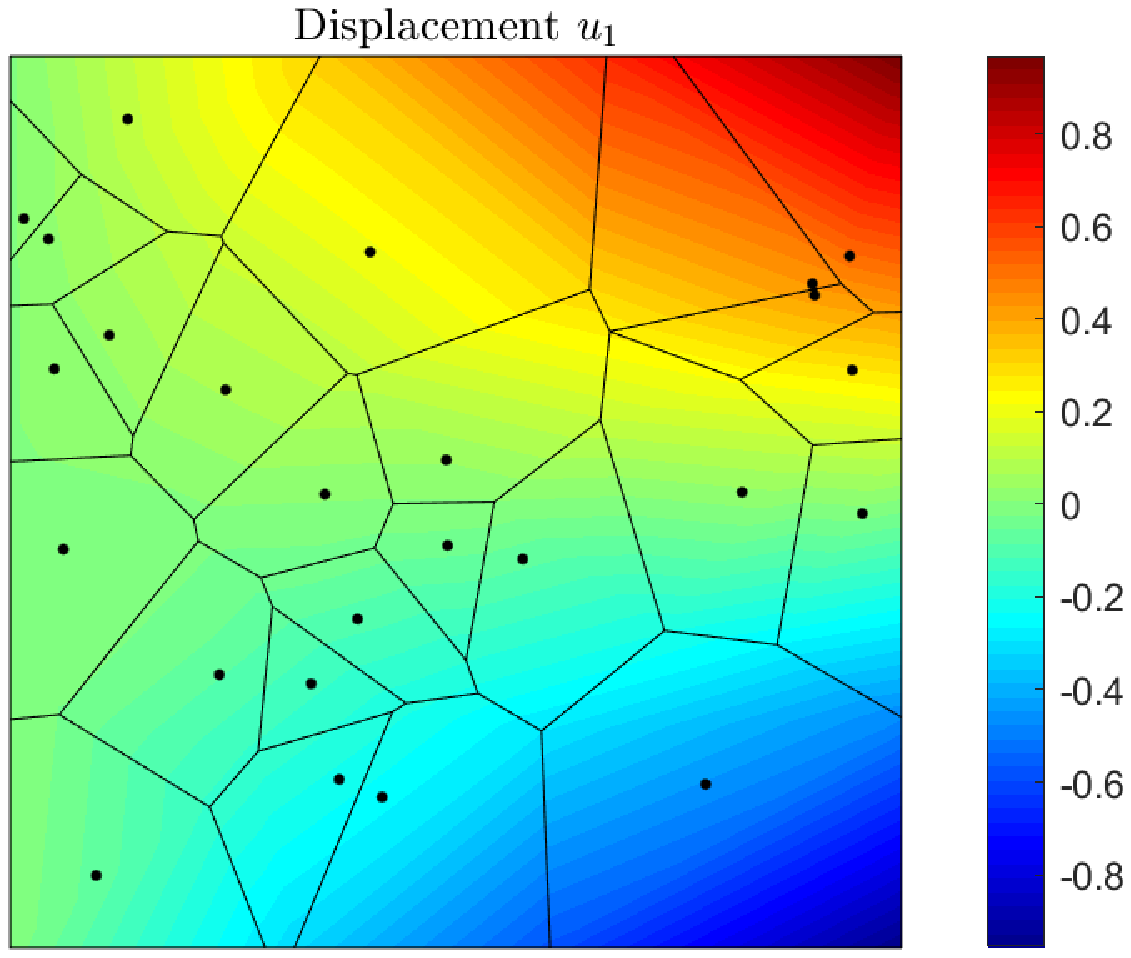}}  
    \subfigure[]{ 
    \label{fig:Ex00_01_02} 
    \includegraphics[width=0.32\textwidth]{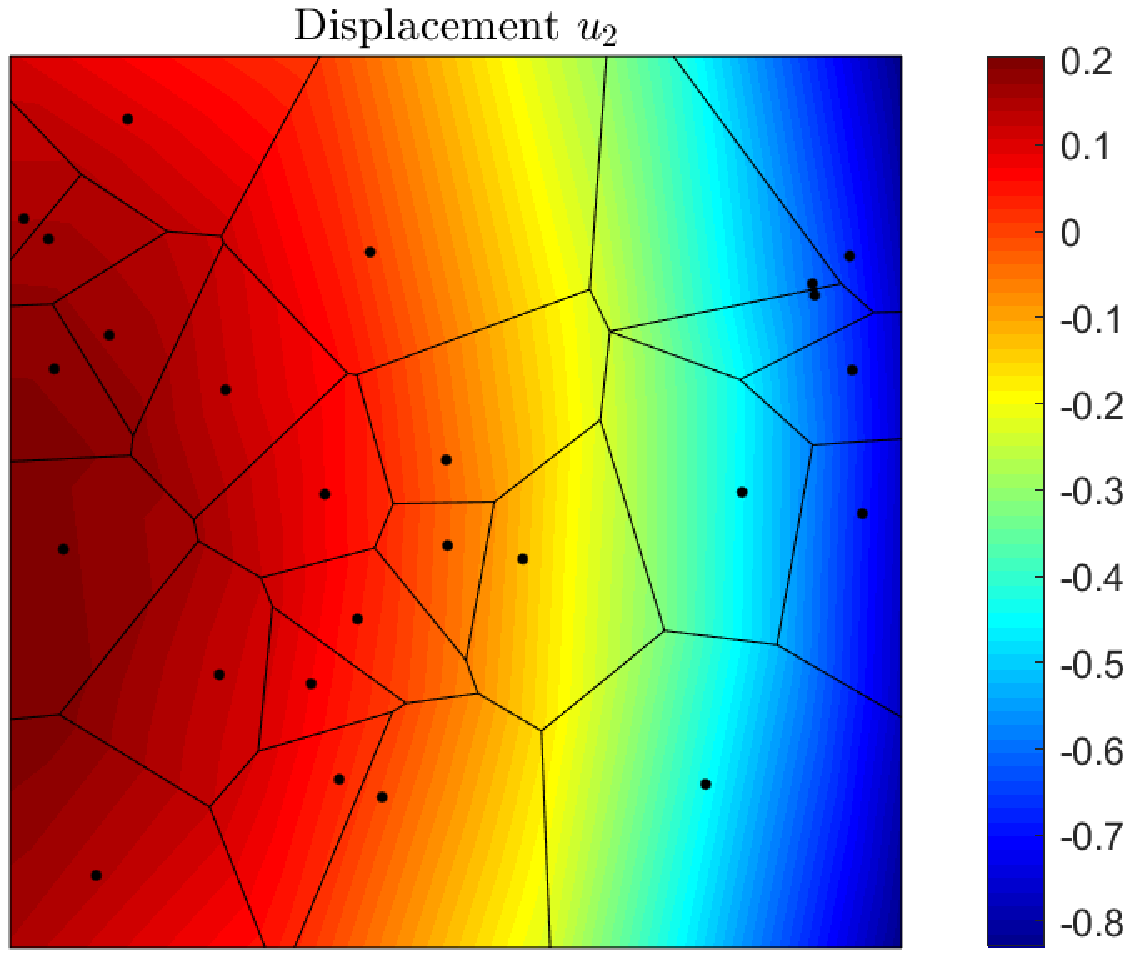}}  
     \subfigure[]{ 
    \label{fig:Ex00_01_03} 
    \includegraphics[width=0.32\textwidth]{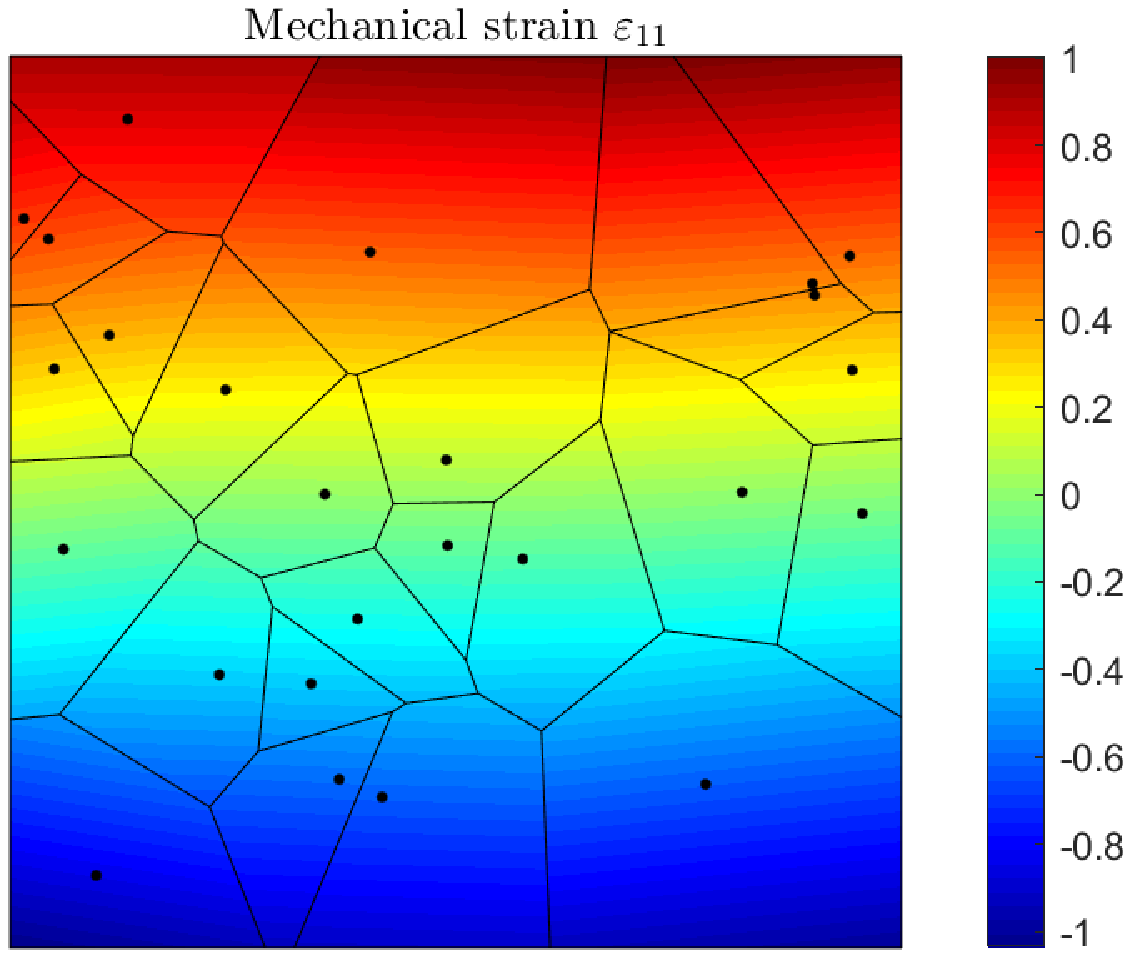}} 
  \caption{The computed solution by FPM with 25 random points. (a) Distribution of displacement $u_1$. (b) Distribution of displacement $u_2$. (c) Distribution of mechanical strain $\varepsilon_{11}$.} 
  \label{fig:Ex00_01} 
\end{figure}

The quadratic patch test is also conducted with the three different point distributions shown in Fig.~\ref{fig:Ex00_Part}. The solution of the FPM with 25 randomly distributed points is presented in Fig.~\ref{fig:Ex00_01}. For all the three kinds of point distributions, the relative errors for the displacement $\mathbf{u}$ and mechanical strain $\boldsymbol{\varepsilon}$ are less than $2 \times 10^{-3}$. 

\section{Numerical examples in continuous domains} \label{sec:continuous}

\subsection{A hollow cylinder}

The first example is an infinite length flexoelectric tube. As shown in Fig.~\ref{fig:Ex1_Schem}, this is a plane strain problem with axisymmetric boundary conditions.

\begin{figure}[htbp] 
  \centering 
    \subfigure[]{ 
    \label{fig:Ex1_Schem} 
    \includegraphics[width=0.48\textwidth]{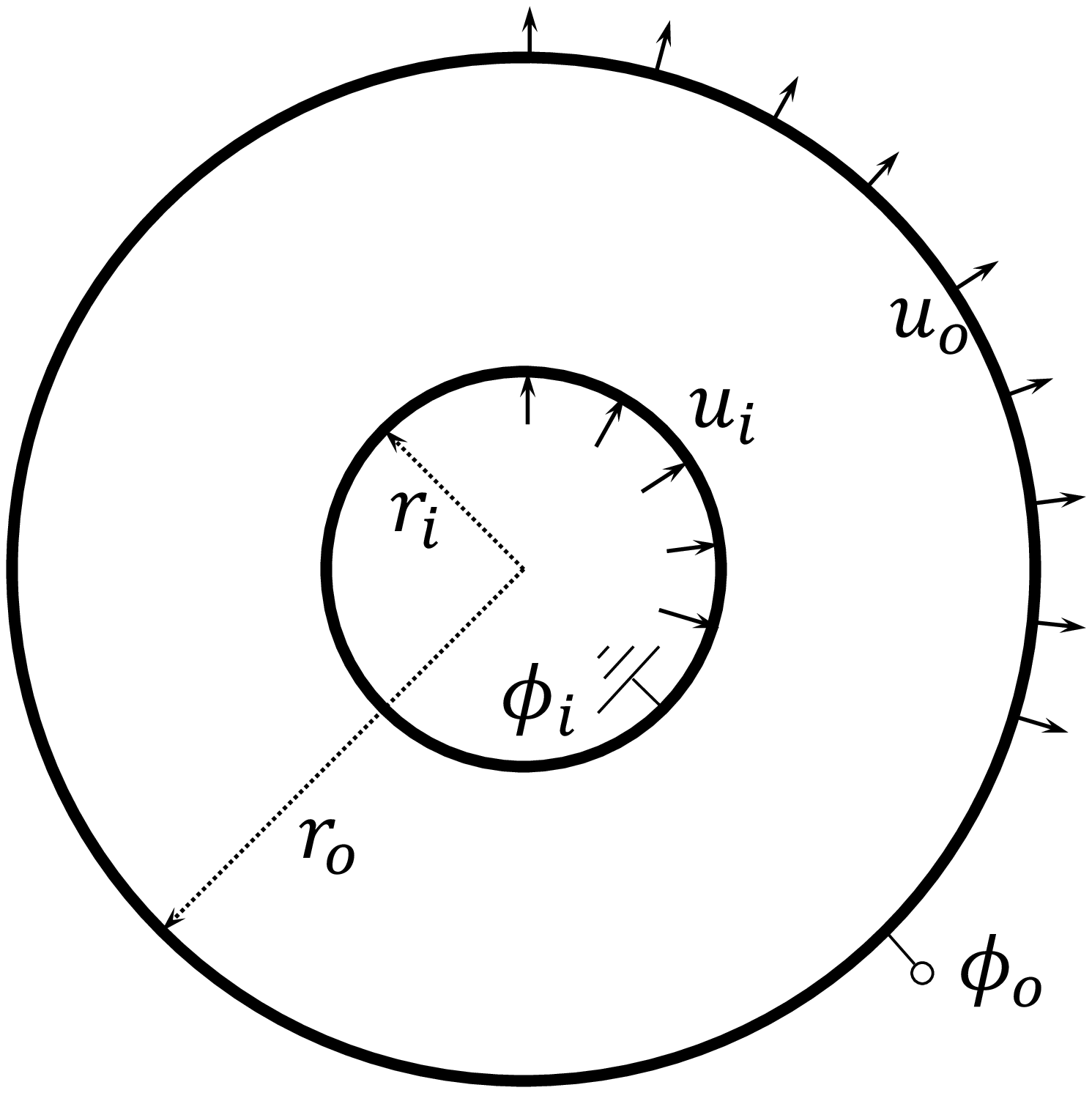}}  
    \subfigure[]{ 
    \label{fig:Ex1_Part} 
    \includegraphics[width=0.48\textwidth]{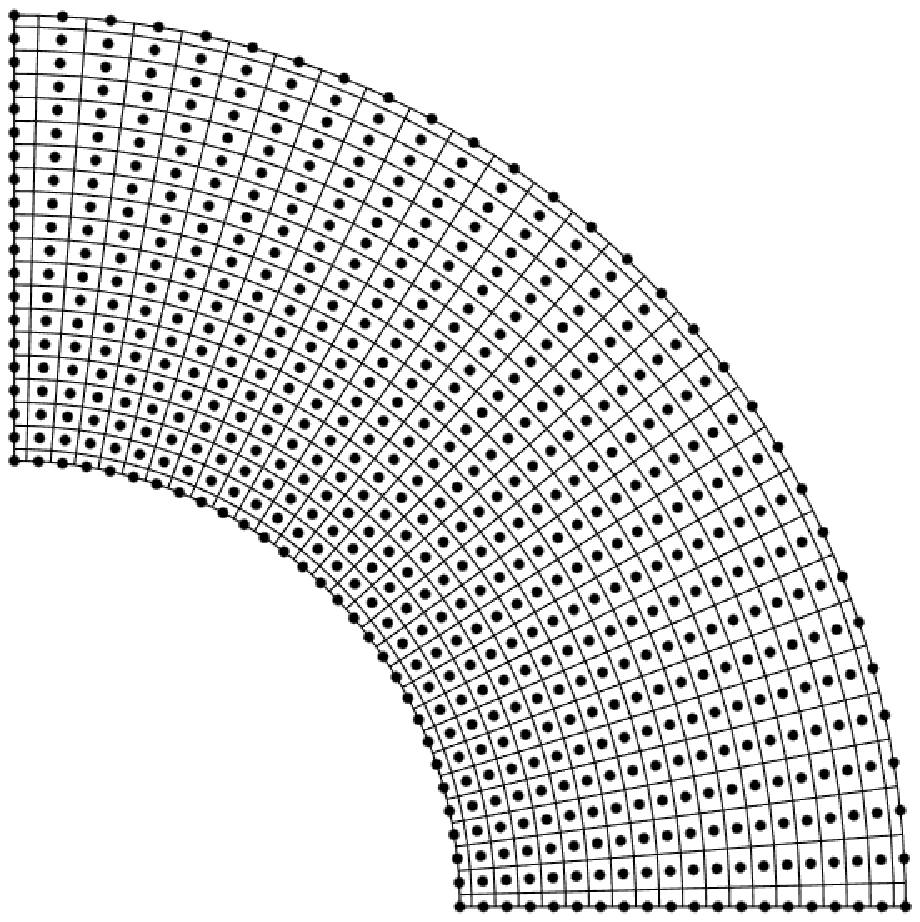}}  
  \caption{(a) An infinite length tube with an axisymmetric cross section. (b) The point distribution and domain partition in the FPM.} 
  \label{fig:Ex1_0} 
\end{figure}

We consider an isotropic flexoelectric material with cubic symmetry. The corresponding constitutive matrices are given in ~\ref{app:sym}, in which Young's modulus $E = 139~\mathrm{GPa}$, Poisson's ratio $\nu = 0.3$, internal material length $l = 2 \mathrm{\mu m}$, flexoelectric coefficients $\overline{\mu}_{11} = \overline{\mu}_{12} = \overline{\mu}_{44} = 1 \times 10^{-6}~\mathrm{C/m}$ and permittivity of the dielectric $\Lambda_{11} =  \Lambda_{33} = 1 \times 10^{-9}~\mathrm{F/m}$. The material is non-piezoelectric, i.e., $e_{31} = e_{33} = e_{15} = 0$. The geometric parameters are: $r_i = 10~\mathrm{\mu m}$, $r_o = 20~\mathrm{\mu m}$. And the boundary conditions are given as: $u_i = 0.045~\mathrm{\mu m}$, $u_o = 0.05~\mathrm{\mu m}$, $\phi_i = 0~\mathrm{V}$, $\phi_o = 1~\mathrm{V}$.

In using FPM, according to the symmetry, only a quarter of the entire domain is considered. Symmetric boundary conditions are applied. Fig.~\ref{fig:Ex1_Part} shows the problem domain and the 600 uniformly distributed points and subdomains used in the FPM. Both the primal and mixed FPMs based on full as well as reduced flexoelectric theories are employed. The computational parameters are given as (if applicable): $c_0 = \sqrt{5}$, $\eta_{11} = 1 \times 10^{10} E$, $\eta_{13} = 1 \times 10^{10} \Lambda_{11}$, $\eta_{21} = 2.0 E$, $\eta_{22} = \eta_{23} = \eta_{24} =0$.

\begin{figure}[htbp] 
  \centering 
    \subfigure[]{ 
    \label{fig:Ex1_01} 
    \includegraphics[width=0.48\textwidth]{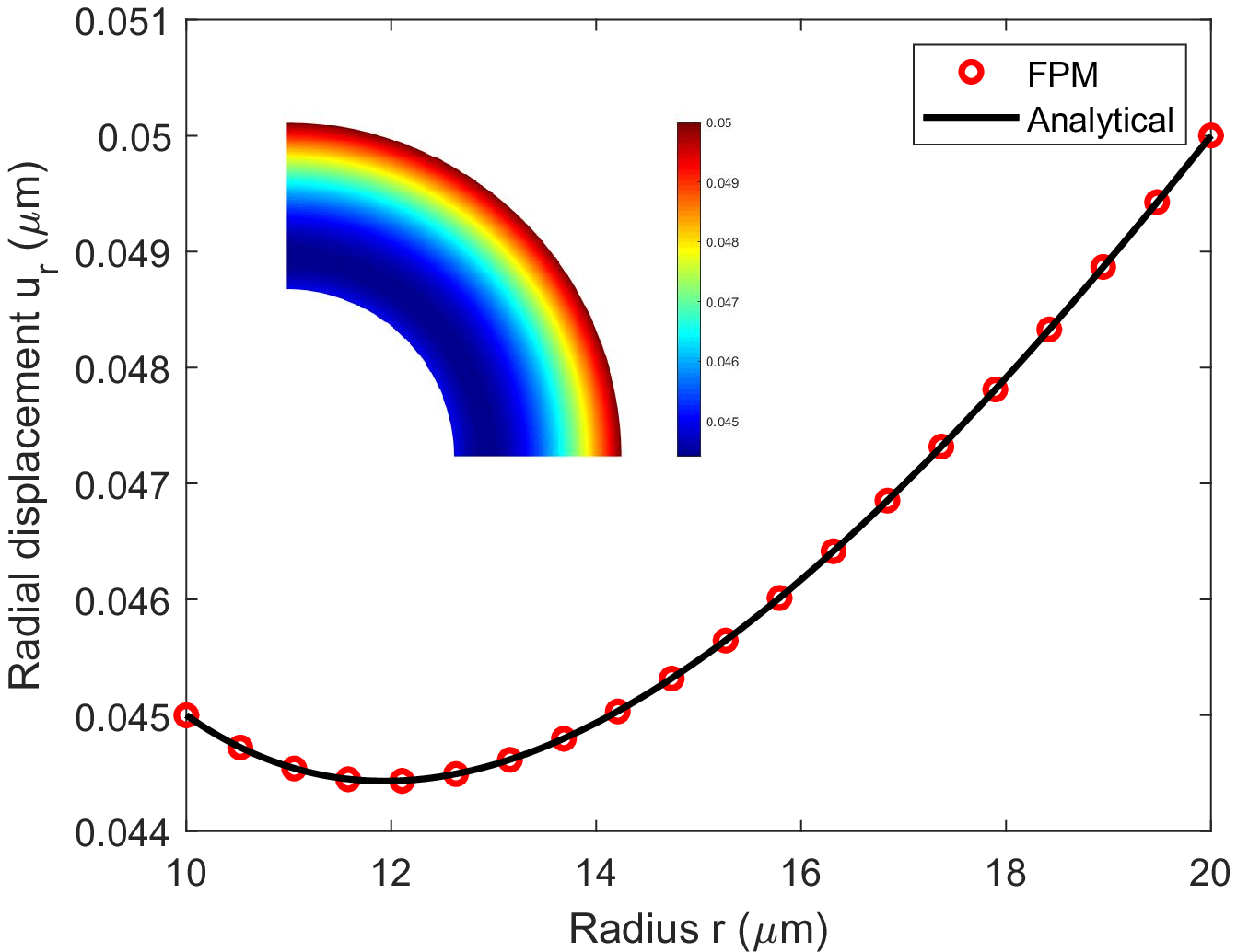}}  
    \subfigure[]{ 
    \label{fig:Ex1_02} 
    \includegraphics[width=0.48\textwidth]{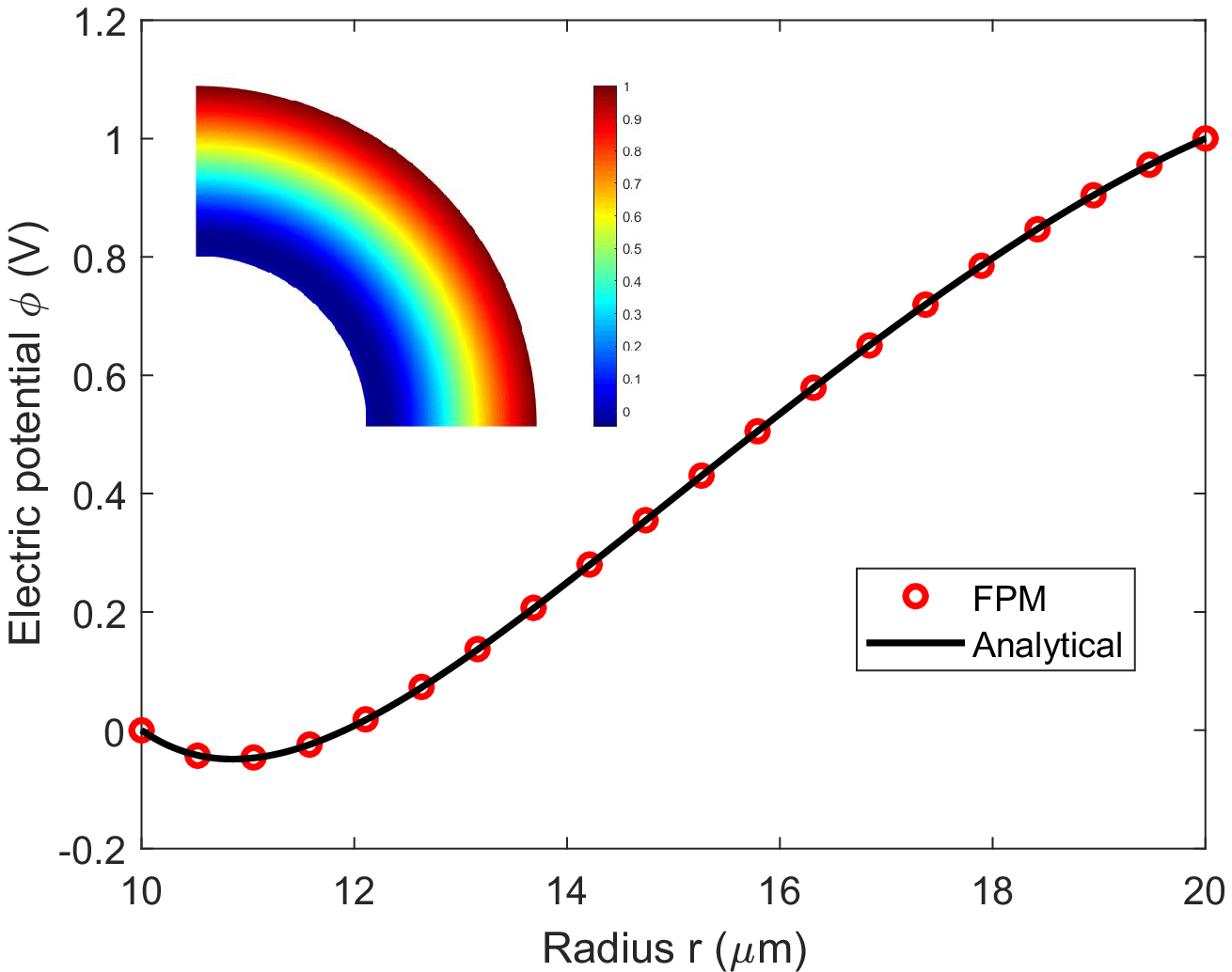}}  
    \subfigure[]{ 
    \label{fig:Ex1_03} 
    \includegraphics[width=0.48\textwidth]{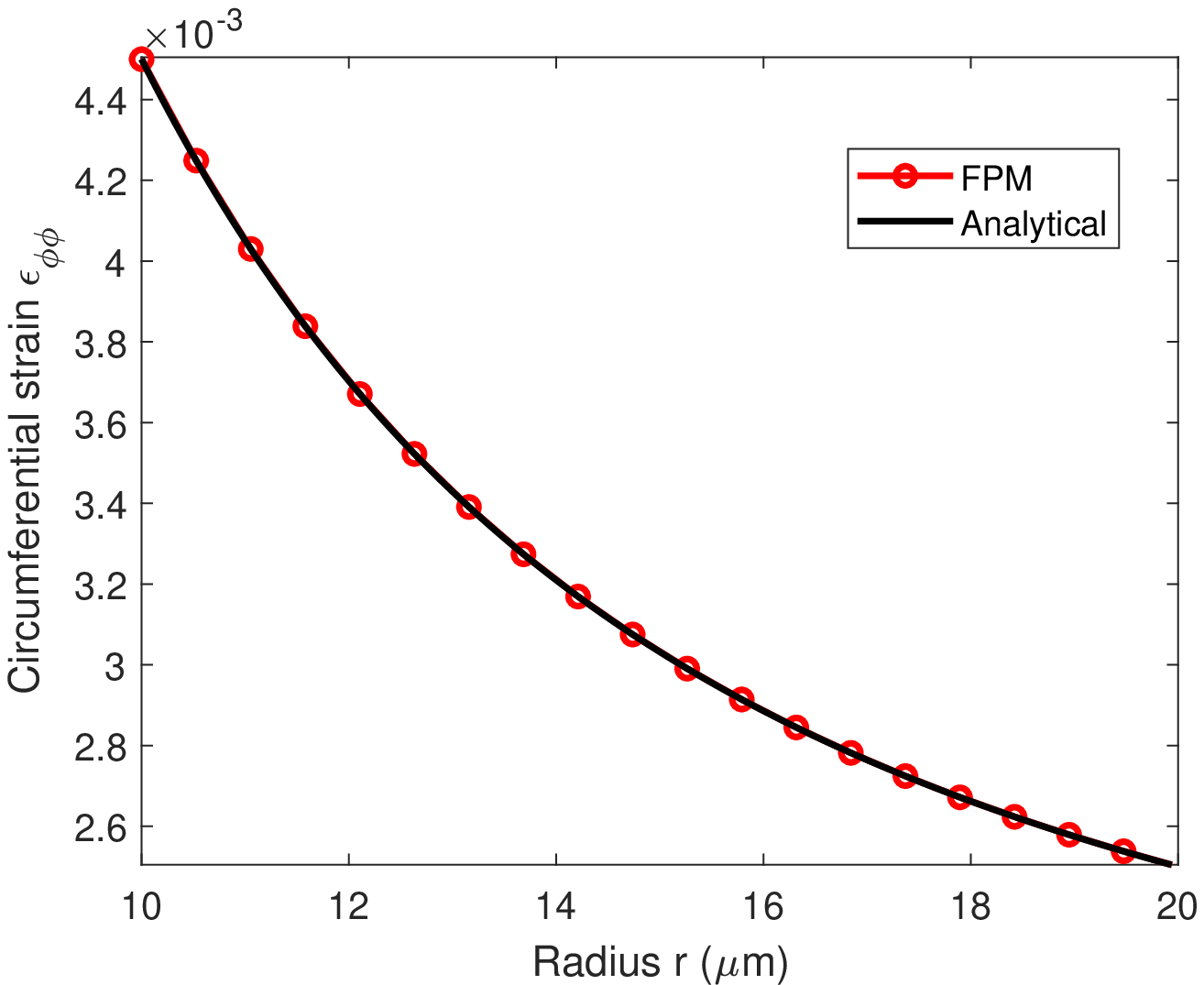}}  
    \subfigure[]{ 
    \label{fig:Ex1_04} 
    \includegraphics[width=0.48\textwidth]{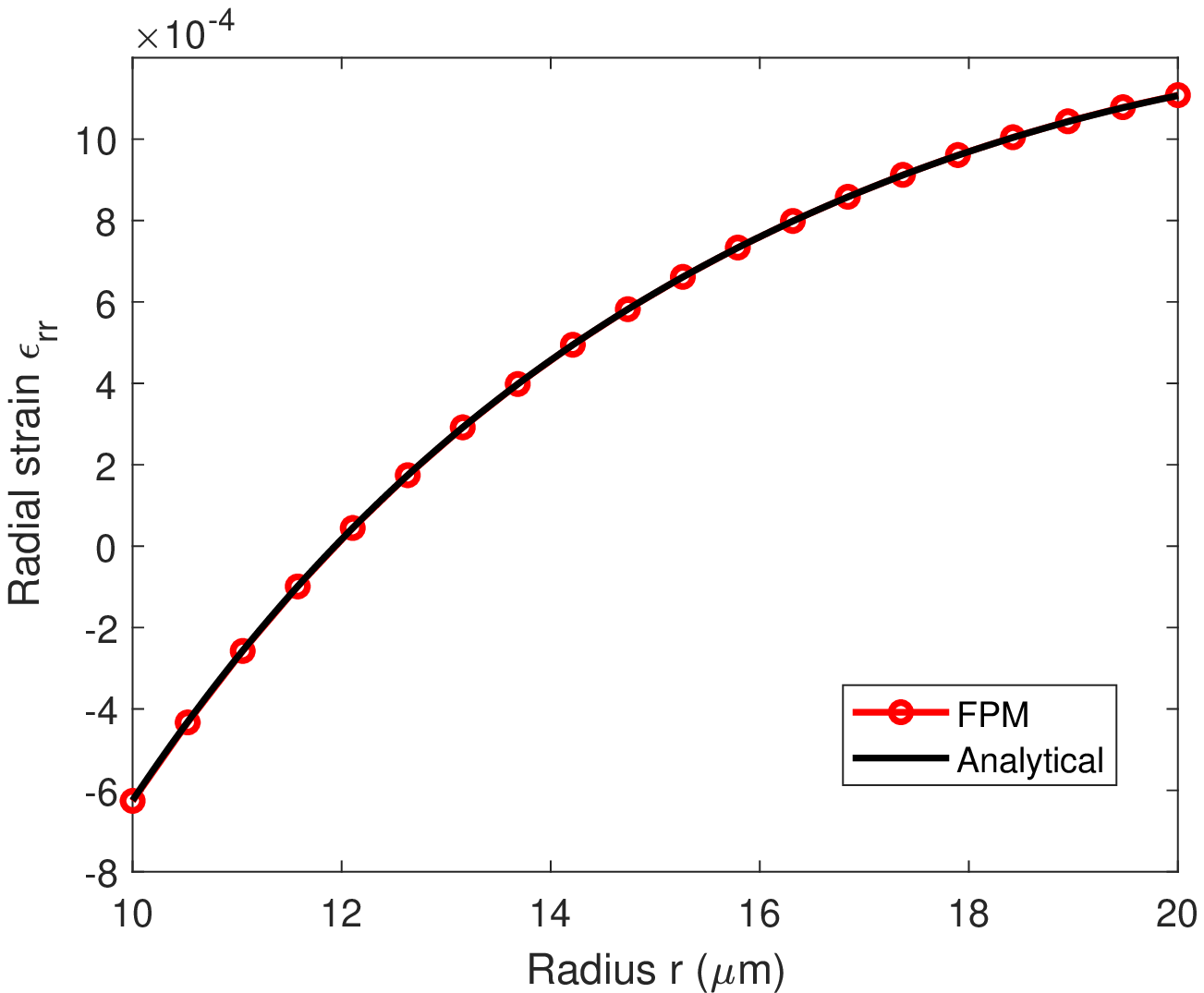}}  
  \caption{The computed solution by primal FPM based on the reduced flexoelectric theory. (a) Distribution of radial displacement. (b) Distribution of electric potential. (c) Distribution of circumferential strain. (d) Distribution of radial strain.} 
  \label{fig:Ex1_Solu} 
\end{figure}

The computed solutions are is presented in Fig.~\ref{fig:Ex1_Solu}. In this example, the solutions based on the full and reduced theories are approximately the same, with a relative error less than $10^{-5}$. Hence only the solutions using the reduced theory are shown here. The numerical results exhibit excellent agreement with the analytical solution based on the reduced flexoelectric theory given by \citet{Mao2014}. And it is also consistent with the solutions achieved by mixed FEM shown in \cite{Deng2017} and \cite{Zhuang2020}.

In this example, the relative errors for the displacement ($\mathbf{u}$), mechanical strain ($\boldsymbol{\varepsilon}$), strain gradient ($\boldsymbol{\kappa}$), electric potential ($\phi$), electric field ($\mathbf{E}$) and electric field gradient ($\mathbf{V}$) obtained by both the primal and mixed methods are shown in Table.~\ref{table:Ex1}. As can be seen, while the primal and mixed FPM with the same number of DoFs may show similar accuracy for the mechanical solutions, the mixed FPM can achieve much better results in the electrical field, especially for high-order variables.

\begin{figure}[htbp] 
  \centering 
    \subfigure[]{ 
    \label{fig:Ex1_Error_01} 
    \includegraphics[width=0.48\textwidth]{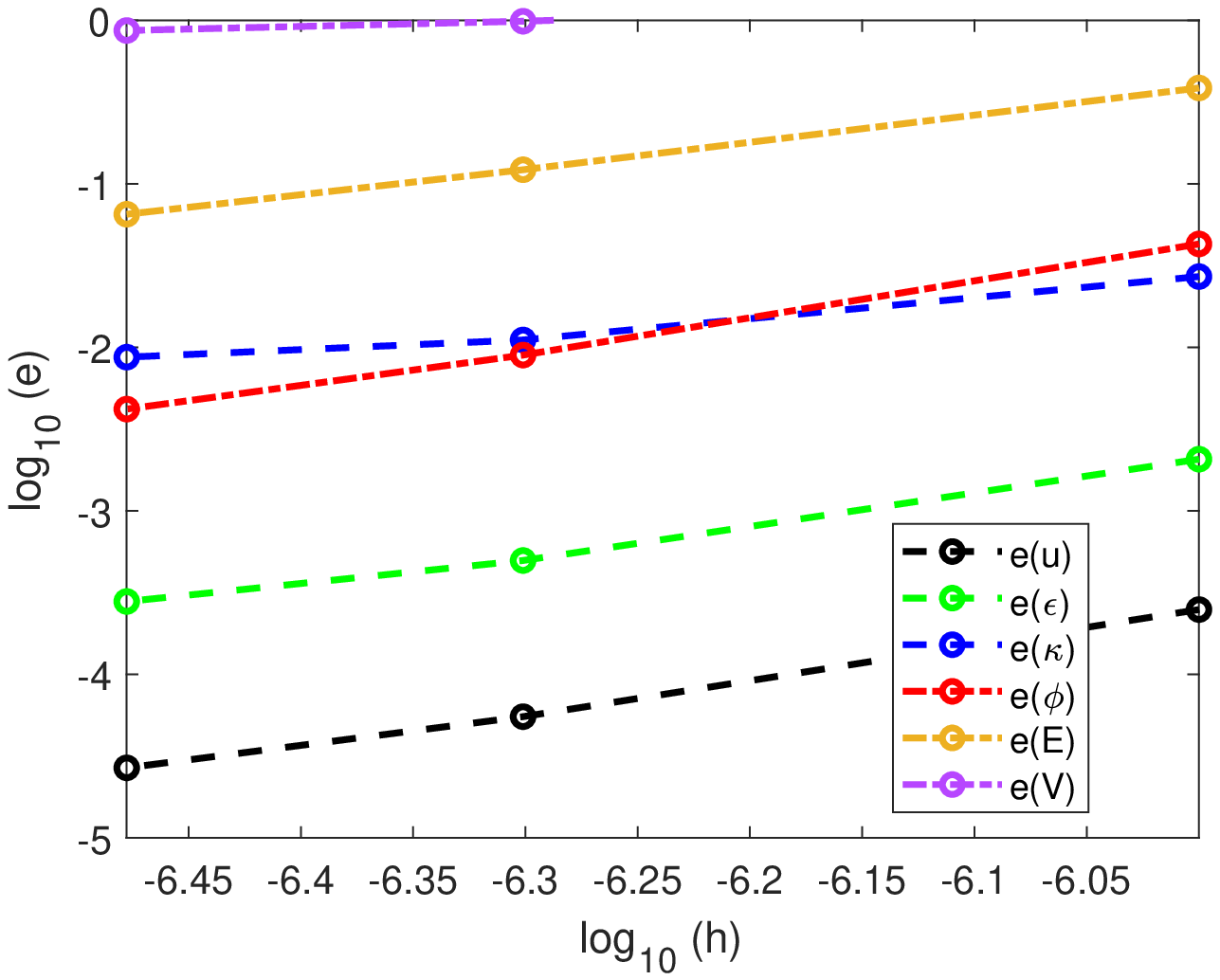}}  
    \subfigure[]{ 
    \label{fig:Ex1_Error_02} 
    \includegraphics[width=0.48\textwidth]{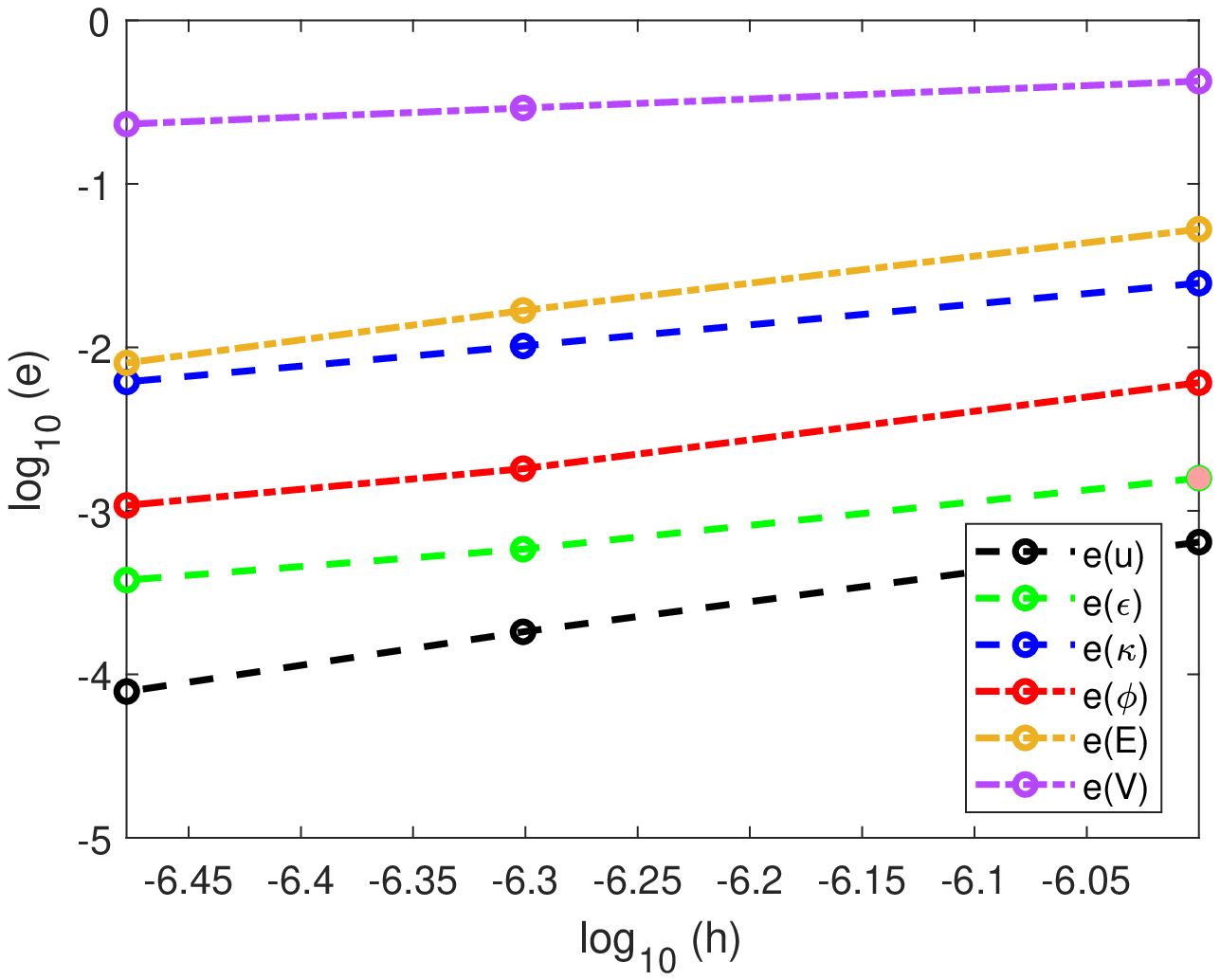}}  
  \caption{Relative errors and convergence rates for (a) primal FPM. (b) mixed FPM} 
  \label{fig:Ex1_Error} 
\end{figure}

In order to demonstrate the convergence of the FPM, Fig.~\ref{fig:Ex1_Error} shows the relative errors for the displacement ($\mathbf{u}$), mechanical strain ($\boldsymbol{\varepsilon}$), strain gradient ($\boldsymbol{\kappa}$), electric potential ($\phi$), electric field ($\mathbf{E}$) and electric field gradient ($\mathbf{V}$) obtained with 50, 200, 800 and 1800 uniformly distributed points respectively, and their dependence on the average distance between two neighboring points ($h$). As can be seen, for both primal and mixed FPM, the convergence rates for displacement, electric potential, as well as the strain and electric field, are all equal to 2. A linear convergence rate is expected for the strain and electric field gradients. However, the primal FPM is not recommended to estimate the electric field gradient because of its lower accuracy.

\begin{table}[htbp]
\caption{Relative errors of the primal and mixed FPM in solving Ex.1.}
\centering
{
\begin{tabular}{ | c | c | c | c | c | c | c |  }
\hline 
Relative error & $e \left( \mathbf{u} \right) $ & $e \left(\boldsymbol{\varepsilon} \right) $ & $e \left(\boldsymbol{\kappa} \right) $  & $e \left( \phi \right) $ & $e \left( \mathbf{E} \right) $ & $e \left( \mathbf{V} \right) $ \\
\hline
Primal FPM & $5.5 \times 10^{-5}$ & $5.0 \times 10^{-4}$ & $1.1 \times 10^{-2}$ & $9.0 \times 10^{-3}$ & $1.2 \times 10^{-1}$ & $9.2 \times 10^{-1}$ \\
\hline
Mixed FPM & $1.2 \times 10^{-4}$ & $5.8 \times 10^{-4}$ & $9.1 \times 10^{-3}$ & $1.8 \times 10^{-3}$ & $1.1 \times 10^{-2}$ & $2.5 \times 10^{-1}$ \\
\hline
\end{tabular}}
\label{table:Ex1}
\end{table}

\subsection{A 2D block}

\begin{figure}[htbp] 
  \centering 
    \subfigure[]{ 
    \label{fig:Ex2_Schem_01} 
    \includegraphics[width=0.48\textwidth]{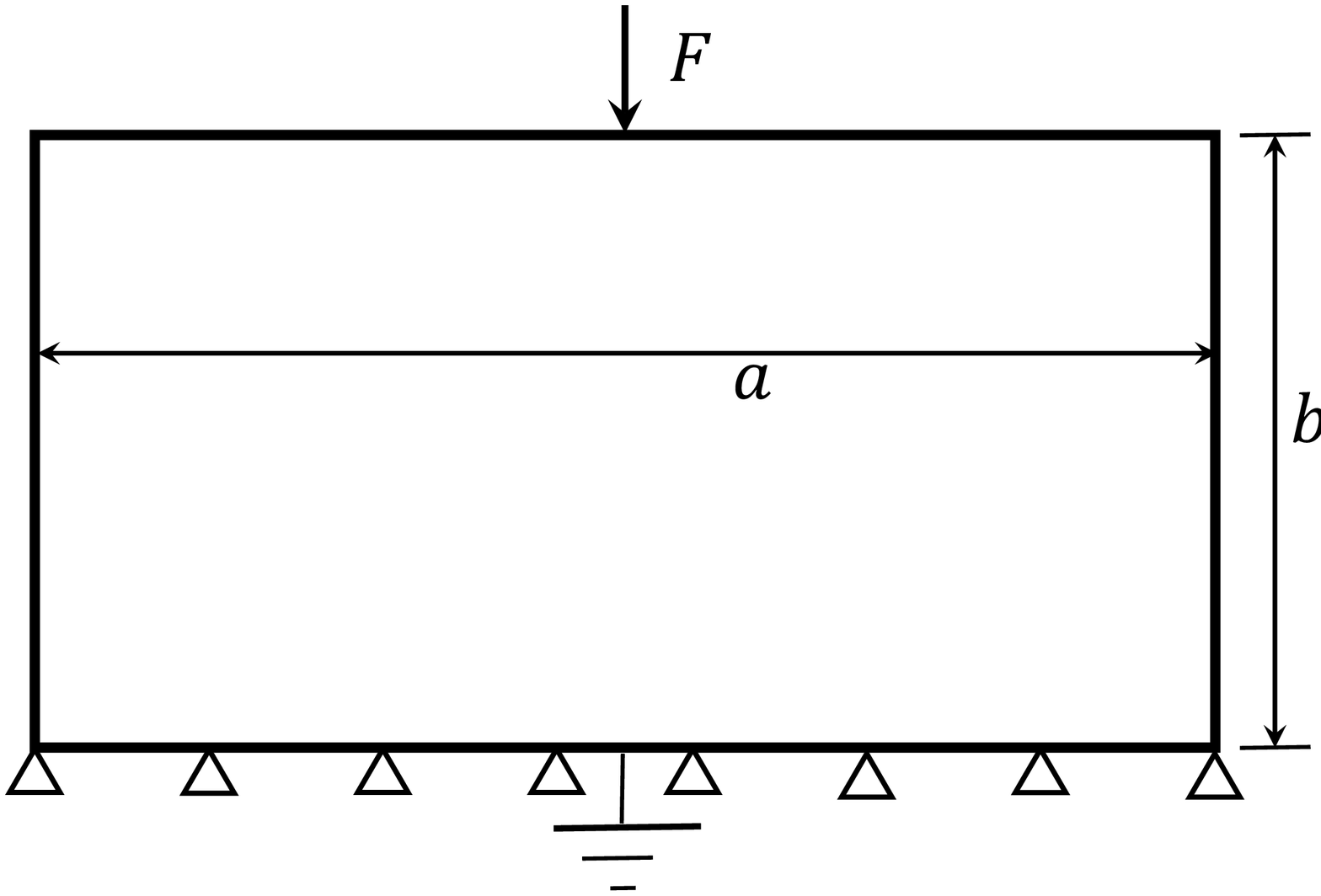}}  
    \subfigure[]{ 
    \label{fig:Ex2_Schem_02} 
    \includegraphics[width=0.48\textwidth]{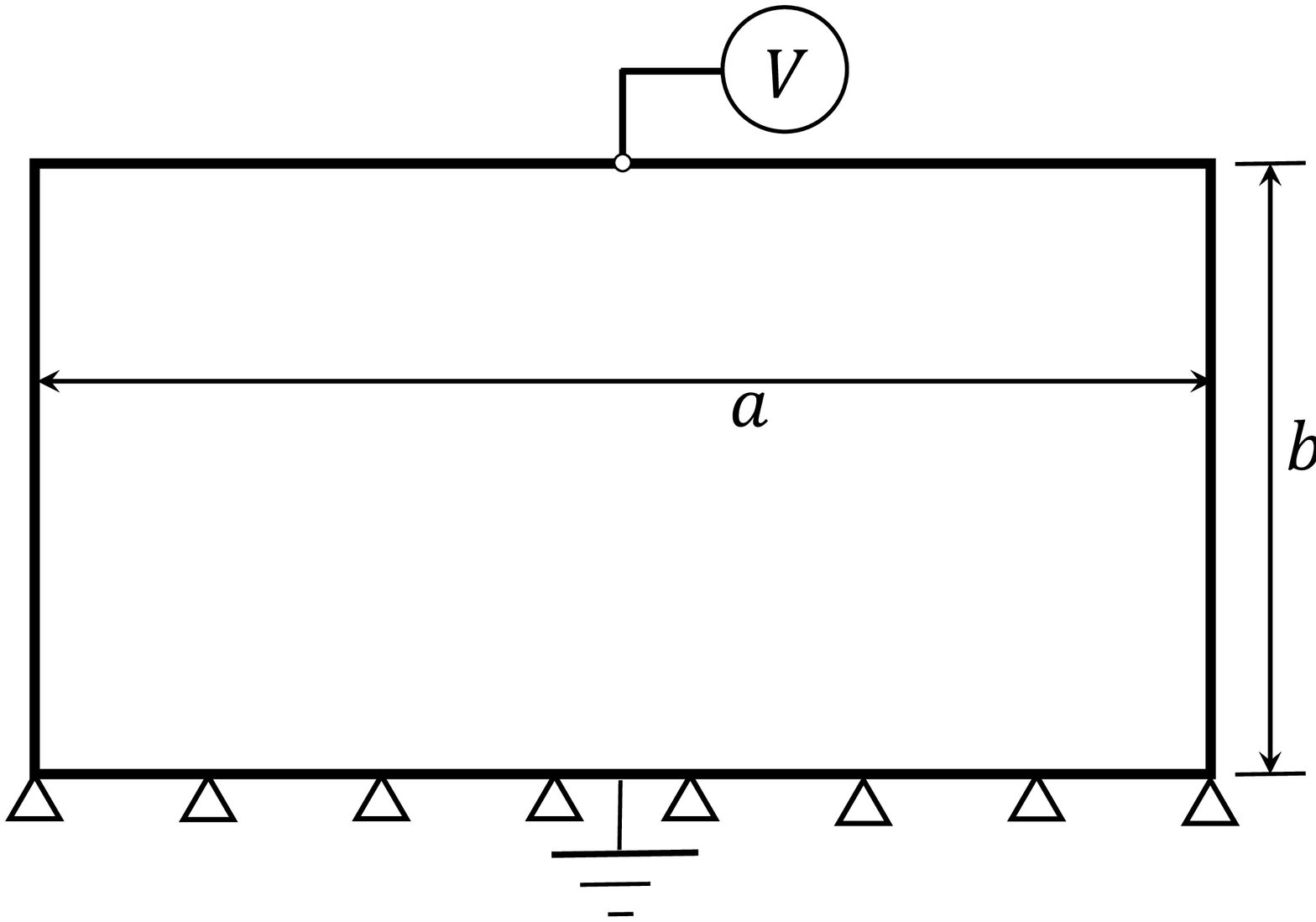}}  
  \caption{(a) 2D block subjected to a concentrated load. (a) 2D block subjected to a concentrated voltage.} 
  \label{fig:Ex2_Schem} 
\end{figure}

Next, we consider the 2D plane strain model of a block subjected to a concentrated load or voltage (see Fig.~\ref{fig:Ex2_Schem}). The example is motivated by the flexoelectric effect observed in some atomic force microscope experiments \cite{Zhuang2020}.  Significant strain gradient and electric field are expected.

\begin{figure}[htbp] 
  \centering 
    \subfigure[]{ 
    \label{fig:Ex2_01} 
    \includegraphics[width=0.48\textwidth]{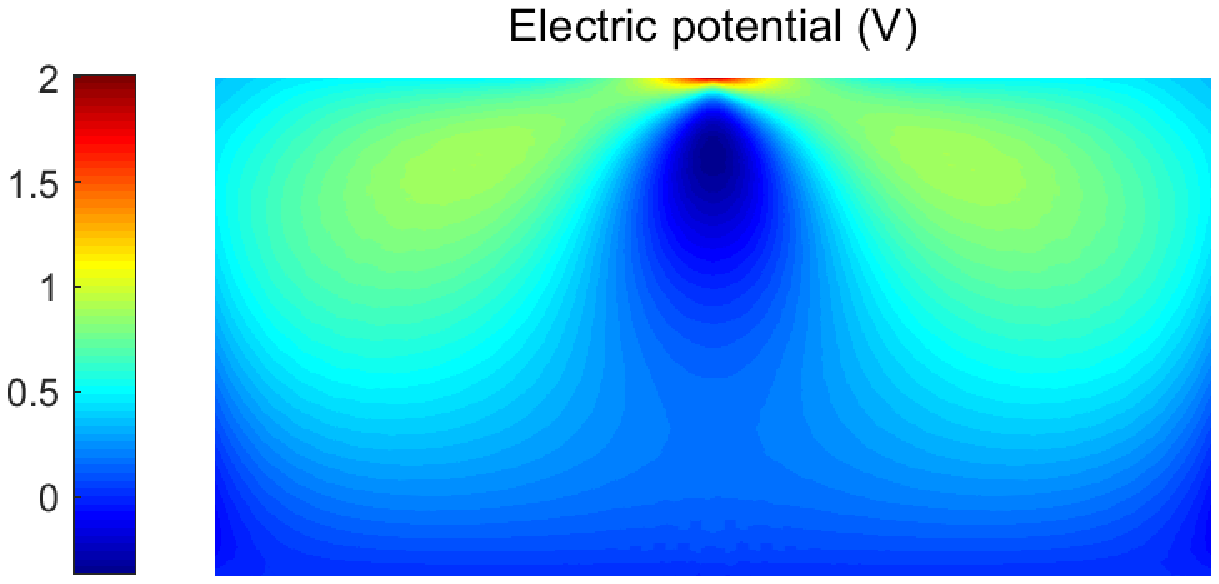}}  
    \subfigure[]{ 
    \label{fig:Ex2_02} 
    \includegraphics[width=0.48\textwidth]{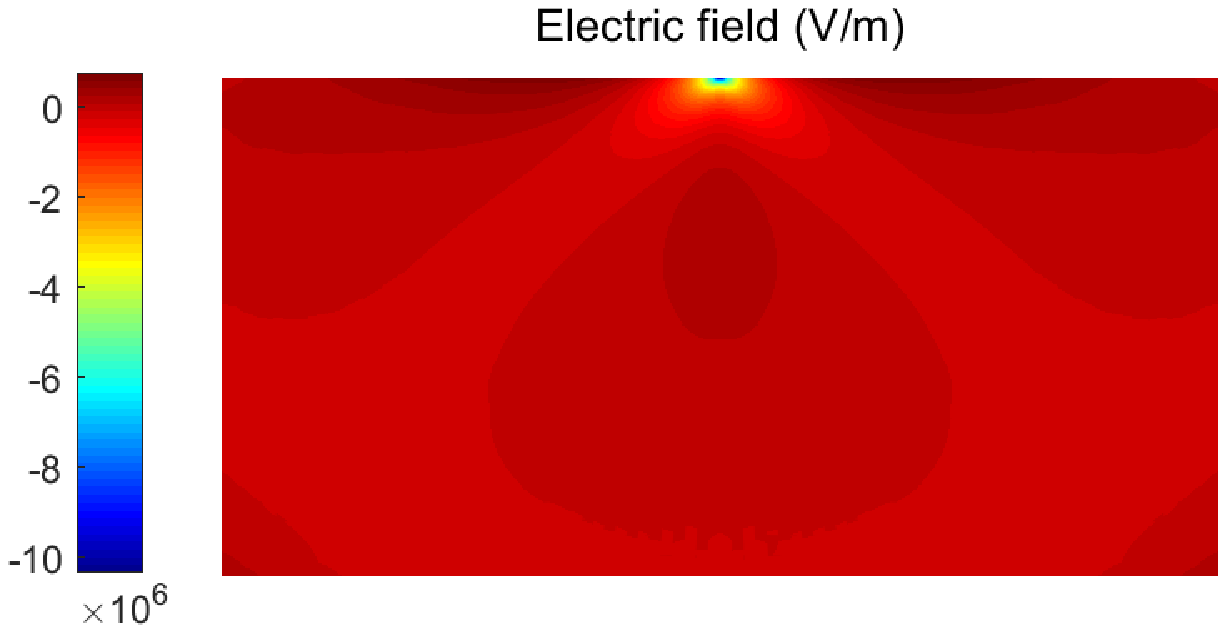}}  
  \caption{The computed solution of mixed FPM for a 2D block subjected to a concentrated load. (a) Distribution of the electric potential $\phi$. (b) Distribution of the electric field $E_2$.} 
  \label{fig:Ex2_Solu_1} 
\end{figure}

\begin{figure}[htbp] 
  \centering 
    \subfigure[]{ 
    \label{fig:Ex2_03} 
    \includegraphics[width=0.48\textwidth]{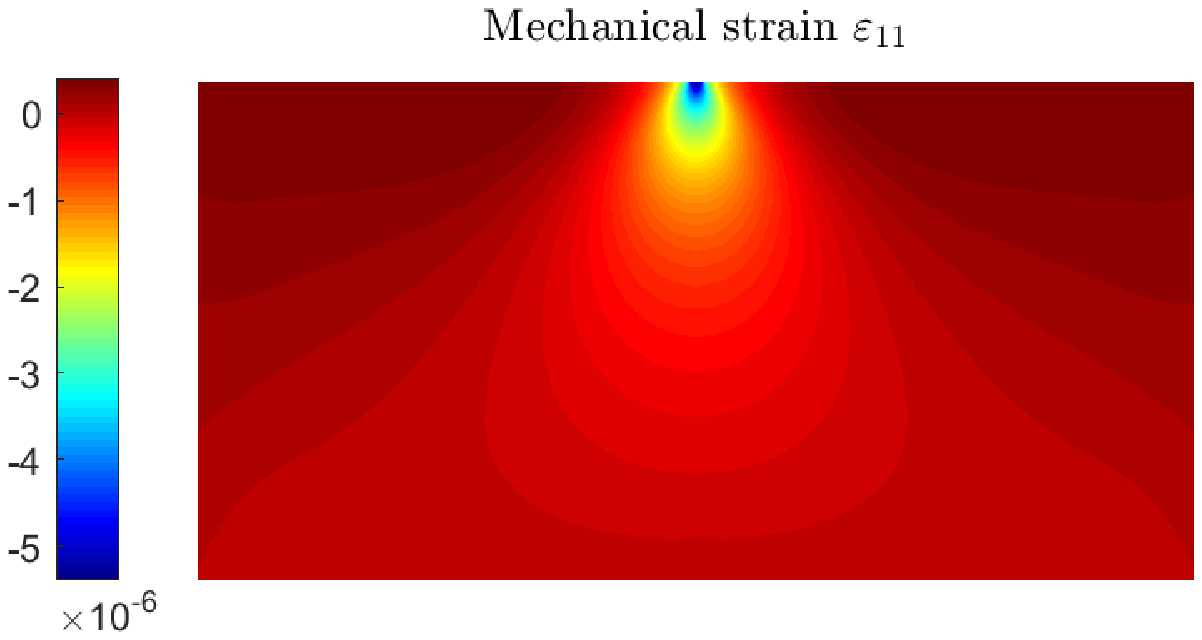}}  
    \subfigure[]{ 
    \label{fig:Ex2_04} 
    \includegraphics[width=0.48\textwidth]{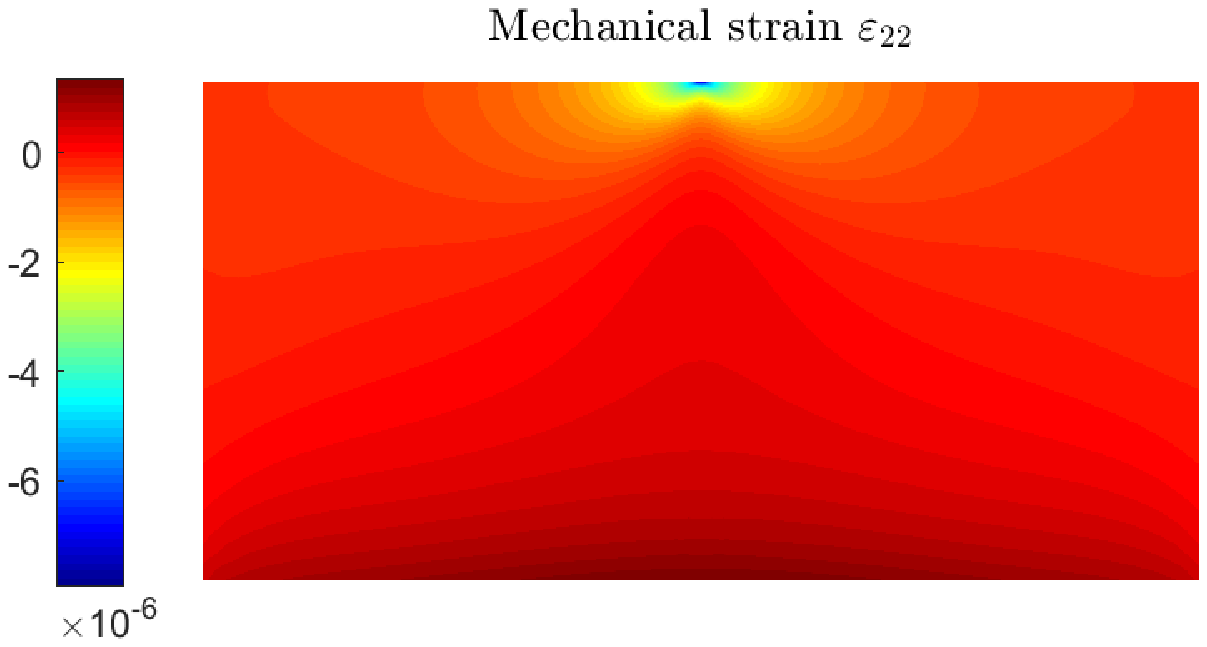}}  
  \caption{The computed solution of mixed FPM for a 2D block subjected to a concentrated voltage. (a) Distribution of the mechanical strain $\varepsilon_{11}$. (b) Distribution of the mechanical strain $\varepsilon_{22}$.} 
  \label{fig:Ex2_Solu_2} 
\end{figure}

The material properties are the same as in the first example. The geometric parameters are: $a = 20~\mathrm{\mu m}$, $b = 10~\mathrm{\mu m}$. The concentrated load $F = 100~\mathrm{\mu N}$ or voltage $V = 5~\mathrm{V}$ is applied uniformly in an area of $200~\mathrm{nm}$ width, to avoid the singularity. Based on the reduced flexoelectric theory, both the primal and mixed FPM are applied. The ABAQUS preprocessing module helps to generate the domain partition. The elements are converted into subdomains in the FPM, and the corresponding Points are placed at the centeroid of each subdomain. With 3200 points being employed, the FPM solutions are shown in Fig.~\ref{fig:Ex2_Solu_1} and \ref{fig:Ex2_Solu_2}. The solutions achieved by the primal and mixed FPMs agree well with each other. For simplicity, only the mixed FPM solution is shown at here.

The computational parameters are: $c_0 = \sqrt{20}$, $\eta_{11} = 1 \times 10^{10} E$, $\eta_{13} = 1 \times 10^{10} \Lambda_{11}$, $\eta_{21} = 1.0 E$, $\eta_{22} = \eta_{23} = \eta_{24} =0$. Fig.~\ref{fig:Ex2_Schem} shows the electrical response of the 2D block subjected to the concentrated load. Compared with the other areas, a much larger electric potential and electric field can be found under the loading point. Without any external electric load, the material exhibits a remarkable variation of electric field as a result of the flexoelectric phenomena. On the other hand, if the block is subjected to a concentrated voltage, a significant mechanical strain takes place (see Fig.~\ref{fig:Ex2_Solu_2}), which implies a local deformation caused by the external electric field. The same example is also considered in \cite{Deng2017} using the mixed FEM. And our solution shows qualitative consistency with the previous study.

\subsection{A cantilever Beam}

In the third example, a 2D cantilever beam is considered. The beam is subjected to a point load at the tip. Two kinds of electrical boundary conditions are imposed. First, the ``open circuit’' means that the electrode at the free end of the beam is grounded, i.e., $\widetilde{\phi} = 0$, while all the other boundaries are surface charge-free. Otherwise, the top and the bottom electrodes retain an electric potential difference, that is, the electric potential on the top surface are fixed to be zeros, whereas the bottom electrode is induced to have an unknown electric potential $V$. This is called a ``closed circuit'’ boundary condition. The material considered in this example is single barium titanate ($\mathrm{BaTiO_3}$) crystal, and the corresponding material parameters are: $E = 100~\mathrm{GPa}$, Poisson's ratio $\nu = 0.37$, internal material length $l = 0$, $e_{31} = -4.4~\mathrm{C/m^2}$, $e_{33} = e_{15} = 0$, $\overline{\mu}_{12} = 1 \times 10^{-6}~\mathrm{C/m}$, $\overline{\mu}_{11} = \overline{\mu}_{44} = 0$, $\Lambda_{11} = 11 \times 10^{-9}~\mathrm{F/m}$ and $ \Lambda_{33} = 12.48 \times 10^{-9}~\mathrm{F/m}$. Note that the material is anisotropic and piezoelectric, with no strain gradient effect. The poling direction is parallel with $y$-axis.

\begin{figure}[htbp] 
  \centering 
    \subfigure[]{ 
    \label{fig:Ex3_01} 
    \includegraphics[width=0.48\textwidth]{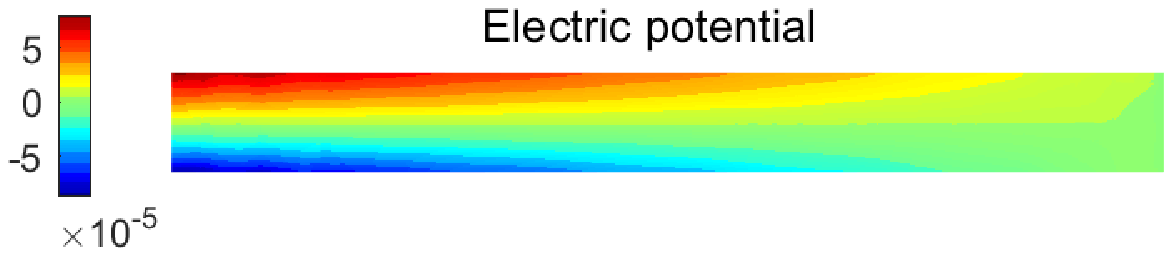}}  
    \subfigure[]{ 
    \label{fig:Ex3_02} 
    \includegraphics[width=0.48\textwidth]{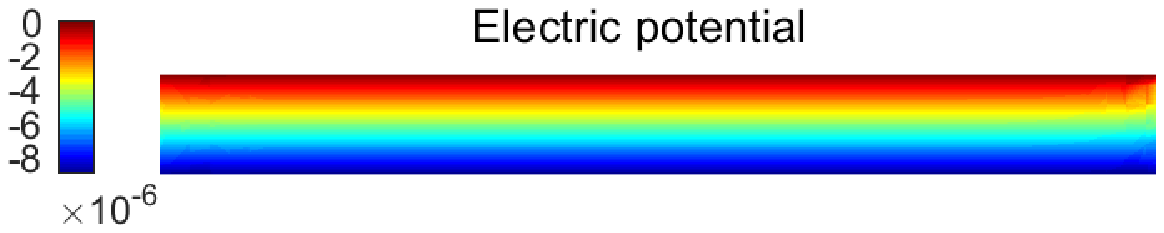}}  
  \caption{Distribution of the electric potential $\phi$. (a) Open circuit. (b) Closed circuit..} 
  \label{fig:Ex3_Solu_1} 
\end{figure}

For a given geometry, $h = 0.12468~\mathrm{\mu m}$, $L = 10 h$, when the point load $F = 100~\mathrm{\mu N}$, the FPM solution of the electric potential $\phi$ is presented in Fig.~\ref{fig:Ex3_Solu_1}. We used the primal FPM with 306 uniformly distributed points. The computational parameters are given as: $c_0 = \sqrt{20}$, $\eta_{11} = 1 \times 10^{10} E$, $\eta_{13} = 1 \times 10^{10} \Lambda_{11}$, $\eta_{21} = \eta_{22} =\eta_{23}=0$. The result shows good agreement with the previous studies based on mixed FEM \cite{Zhuang2020} and a Local Maximum Entropy (LME)-based meshless method \cite{Abdollahi2014}. It should be pointed out that in the primal FPM, though the trial function is discontinuous, its value in each subdomain is still dependent on the values of the neighboring points, and thus keeps a weak continuity. As a result, even if  we do not enforce continuity across the internal boundaries, i.e., $\eta_{21} = \eta_{22} =\eta_{23} = \eta_{24} =0$, the primal FPM is still a consistent algorithm. Yet an appropriate penalty parameter can help to improve the accuracy and stability of the solution, especially for the high-order variables.

\begin{figure}[htbp] 
  \centering 
    \includegraphics[width=0.55\textwidth]{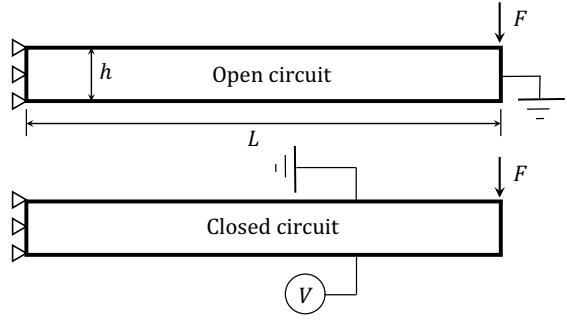}
  \caption{Cantilever beams subjected to a concentrated load and open/closed circuit.} 
  \label{fig:Ex3_Schem} 
\end{figure}

\begin{figure}[htbp] 
  \centering 
    \includegraphics[width=0.8\textwidth]{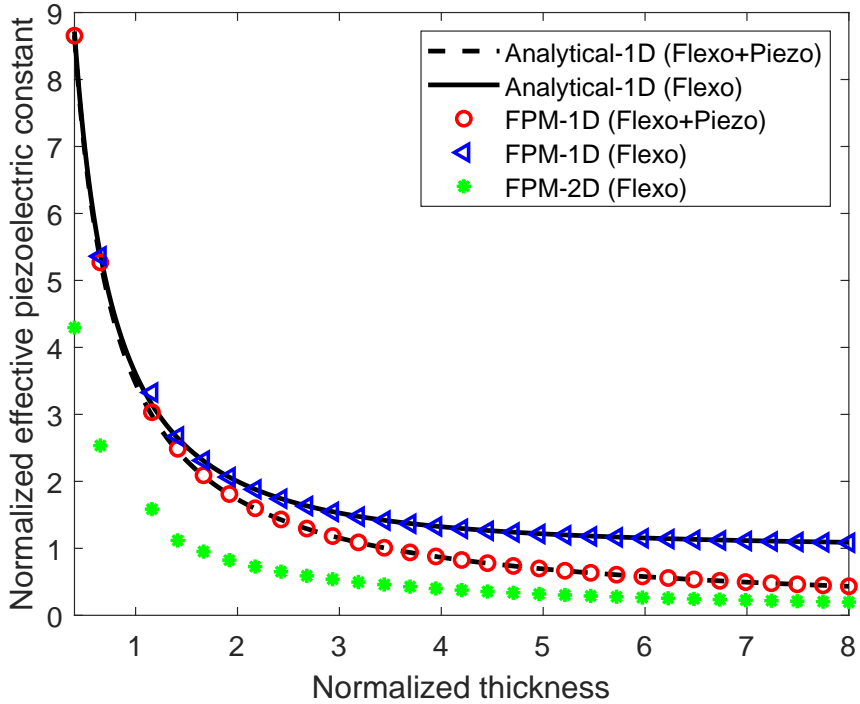}
  \caption{The computed solution: effective normalized piezoelectric constant.} 
  \label{fig:Ex3_Solu_2} 
\end{figure}

Furthermore, we define an electromechanical coupling factor $k_{eff}$ as the ratio between the electrostatic energy and elastic energy:
\begin{align}
\begin{split}
k_{eff} = \sqrt{ \frac{\frac{1}{2} \int_\Omega \boldsymbol{\varepsilon}^\mathrm{T} \mathbf{C}_{\sigma \varepsilon} \boldsymbol{\varepsilon} \mathrm{d} \Omega }{\frac{1}{2}  \int_\Omega \mathbf{E}^\mathrm{T} \boldsymbol{\Lambda}  \mathbf{E} \mathrm{d} \Omega} }.
\end{split}
\end{align}
And a normalized piezoelectric constant $\overline{e}$ is defined as the ratio between $k_{eff}$ obtained with and without considering the flexoelectric effect:
\begin{align}
\begin{split}
\overline{e} = \frac{k_{eff}}{k_{piezo}},
\end{split}
\end{align}
where $ k_{piezo}$ is the value of $k_{eff}$ when the flexoelectric effect is eliminated. Thus $\overline{e}$ is given as an estimate of the performance of the flexoelectric phenomenon. Fig.~\ref{fig:Ex3_Solu_2} shows the variation of $\overline{e}$ for the beam subjected to open circuit condition under different $h$. The normalized beam height $\overline{h} = -e_{31} h / \overline{\mu}_{12}$. A simplified 1D condition is also considered, in which only $\Lambda_{33}$ and $\overline{\mu}_{12}$ are non-zero. The results are compared with an analytical solution \cite{Majdoub2009}:
\begin{align}
\begin{split}
k_{eff} = \frac{\chi}{1+\chi} \sqrt{\frac{e_{31}^2 + 12 \left( \overline{\mu}_{12} / h \right)^2}{\Lambda_{33} E} }, \; \text{where} \; \chi = \frac{\Lambda_{33}}{\epsilon_0} - 1.
\end{split}
\end{align}
As can be seen, the influence of flexoelectricity rises dramatically when the beam height decreases. This  results from the enlargement of the strain gradient. The FPM results have shown excellent agreement with the analytical solution, as well as previous works \cite{Abdollahi2014}.

\subsection{A truncated pyramid}

\begin{figure}[htbp] 
  \centering 
    \subfigure[]{ 
    \label{fig:Ex4_Schem_01} 
    \includegraphics[width=0.48\textwidth]{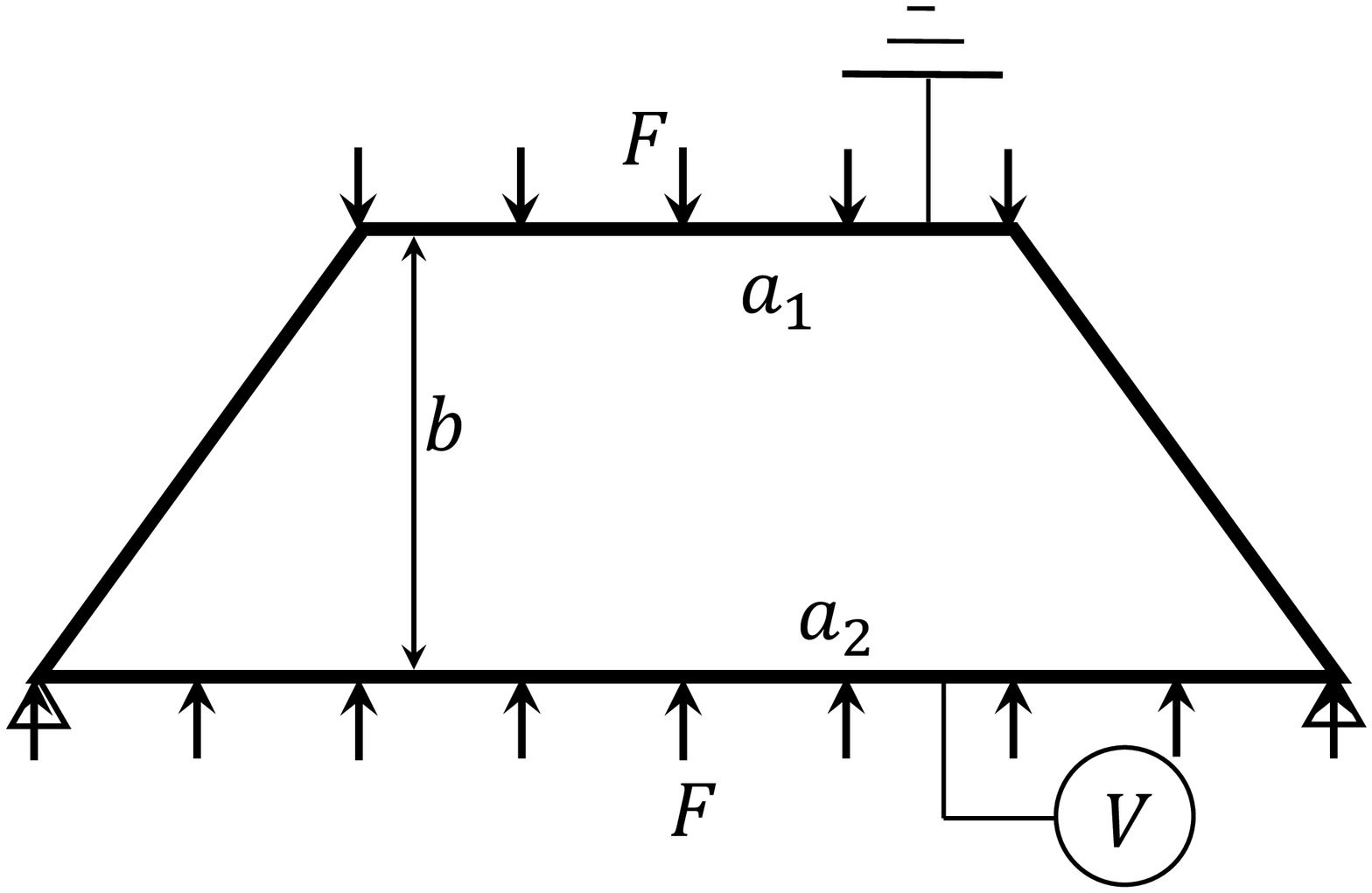}}  
    \subfigure[]{ 
    \label{fig:Ex4_Schem_02} 
    \includegraphics[width=0.48\textwidth]{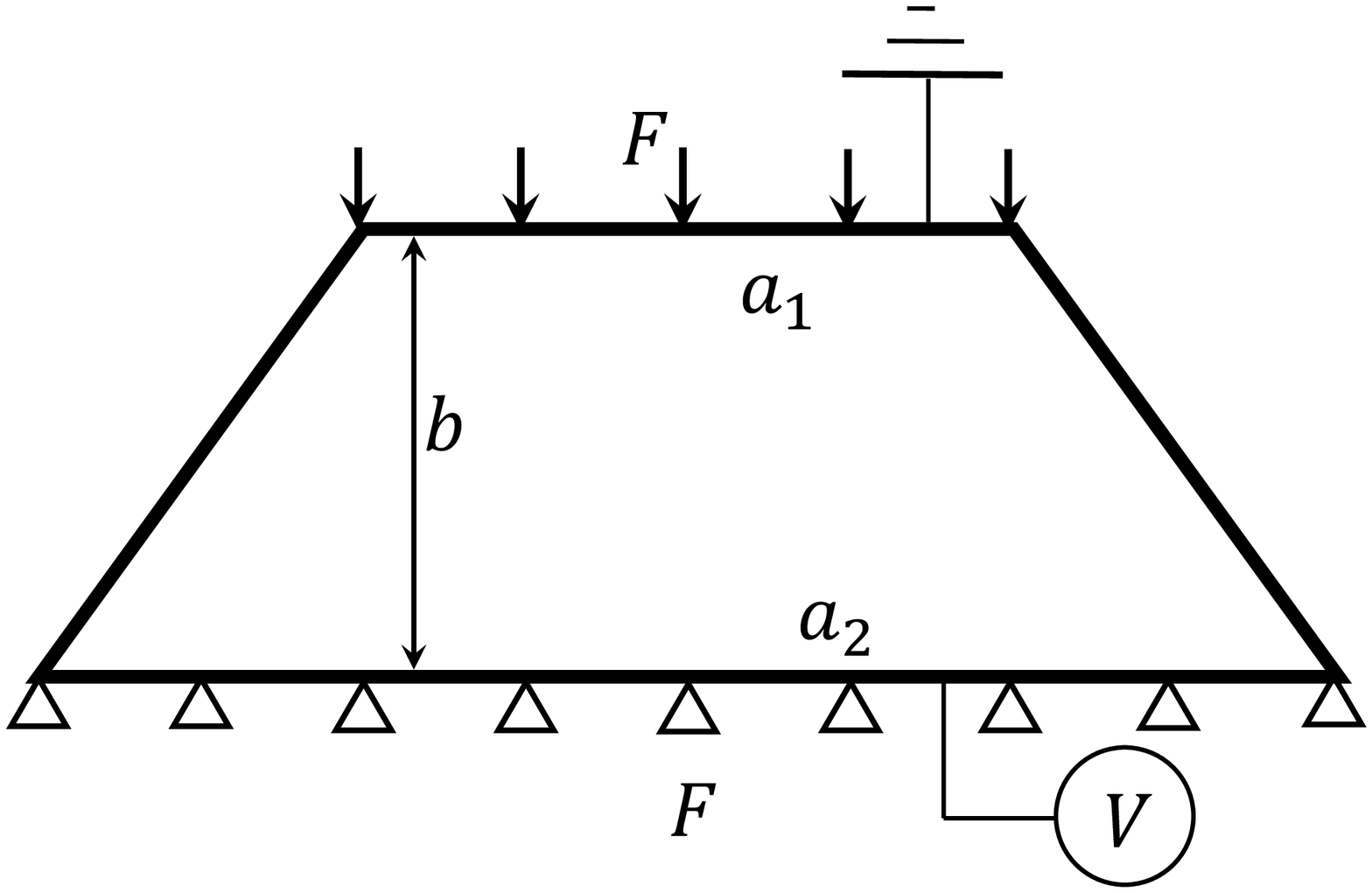}}  
  \caption{Truncated pyramid subjected to a uniformly distributed load with (a) flexible support. (b) rigid support.} 
  \label{fig:Ex4_Schem} 
\end{figure}

In this section, we study  the response of a truncated pyramid in plane strain subjected to uniformly distributed load with flexible and rigid supports (see Fig.~\ref{fig:Ex4_Schem}). The top surface of the truncated pyramid is grounded, i.e., $\widetilde{\phi} = 0$, while the bottom surface is attached to an equipotential electrode with an unknown $V$. The electrical boundary condition is the same as the ``closed circuit'’ in the cantilever beam example. Two kinds of mechanical boundary conditions are considered here. First, the `‘flexible support’' (see Fig.~\ref {fig:Ex4_Schem_01}) represents a simply supported condition and the load $F$ is applied uniformly on both the top and bottom surfaces. Alternatively, the truncated pyramid can also be rigidly supported. Thus the vertical displacement of the bottom surface is fixed to be zero, and a non-uniform reaction force will be induced on the bottom surface.

\begin{figure}[htbp] 
  \centering 
    \subfigure[]{ 
    \label{fig:Ex4_01} 
    \includegraphics[width=0.48\textwidth]{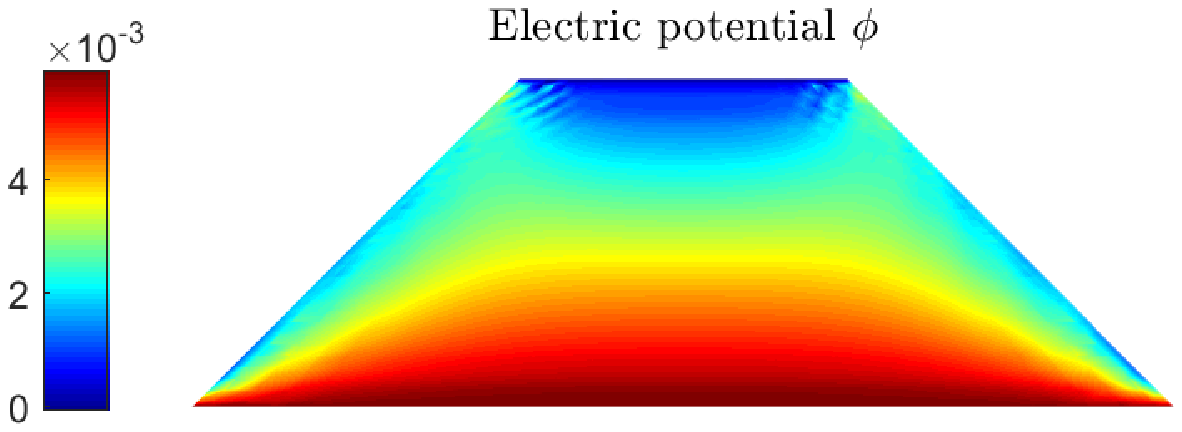}}  
    \subfigure[]{ 
    \label{fig:Ex4_02} 
    \includegraphics[width=0.48\textwidth]{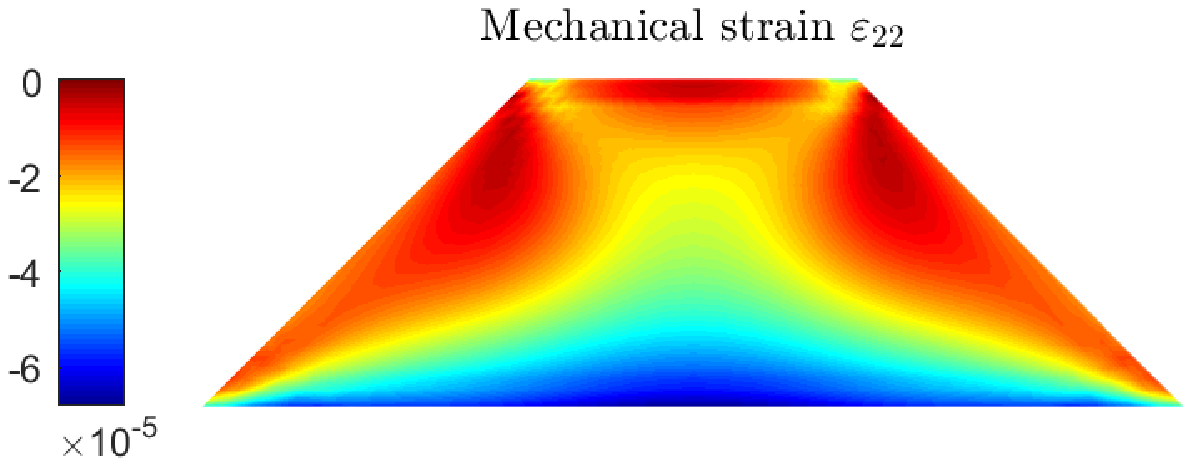}}  
  \caption{The computed solution of primal FPM for a truncated pyramid with flexible support. (a) Distribution of the electric potential $\phi$. (b) Distribution of the mechanical strain $\varepsilon_{22}$.} 
  \label{fig:Ex4_Solu_1} 
\end{figure}

\begin{figure}[htbp] 
  \centering 
    \subfigure[]{ 
    \label{fig:Ex4_03} 
    \includegraphics[width=0.48\textwidth]{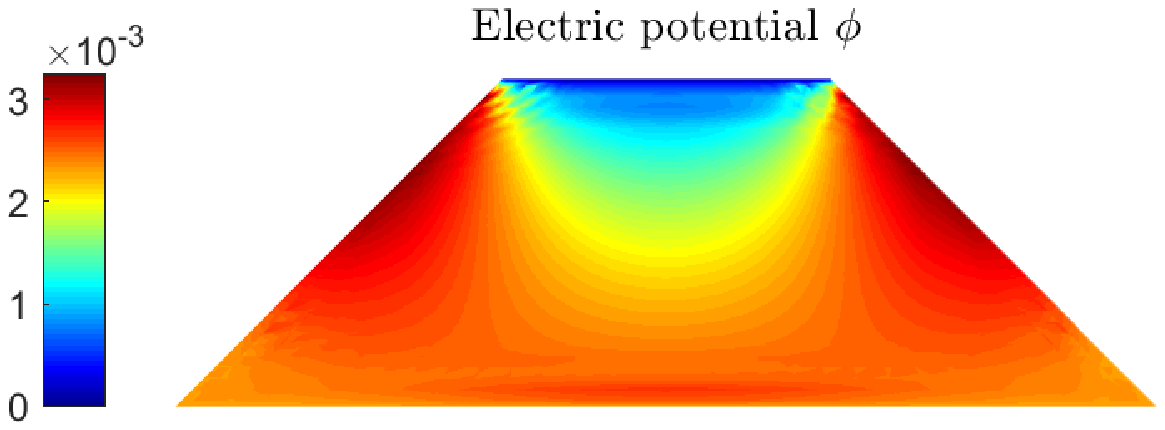}}  
    \subfigure[]{ 
    \label{fig:Ex4_04} 
    \includegraphics[width=0.48\textwidth]{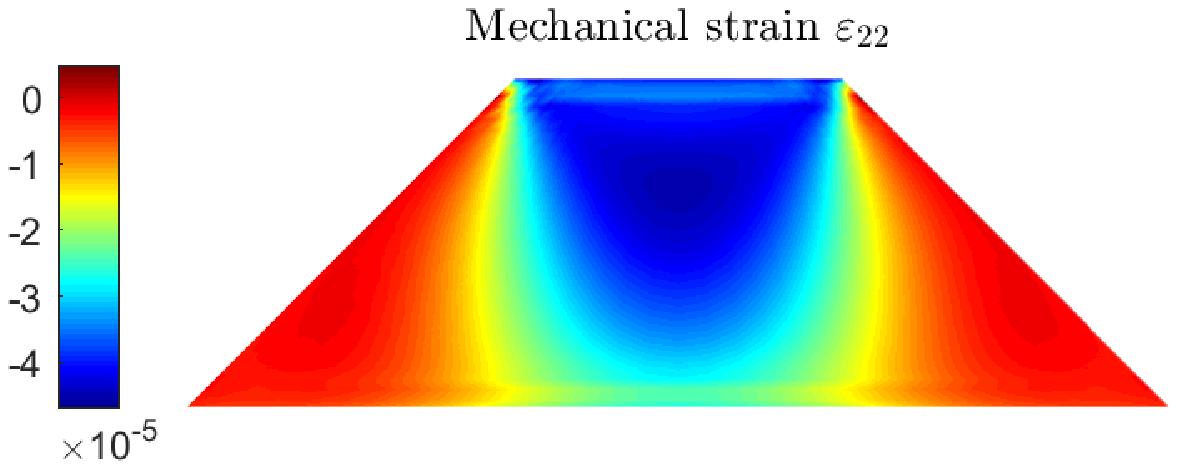}}  
  \caption{The computed solution of primal FPM for a truncated pyramid with rigid support. (a) Distribution of the electric potential $\phi$. (b) Distribution of the mechanical strain $\varepsilon_{22}$.} 
  \label{fig:Ex4_Solu_2} 
\end{figure}

The material is the same as in the previous example ($\mathrm{BaTiO_3}$). And the poling direction is still parallel to the $y$-axis. The geometric parameters are given as: $a_1 = 750~\mathrm{\mu m}$, $a_2 = 2250~\mathrm{\mu m}$, $h = 750~\mathrm{\mu m}$. The load $F = 450~\mathrm{kN}$. With 306 uniform points and quadrilateral partitions are used, the numerical solution of the electric potential $\phi$ and mechanical strain $\varepsilon_{22}$ for both flexible and rigid boundary conditions are presented in Fig.~\ref{fig:Ex4_Solu_1} and \ref{fig:Ex4_Solu_2} respectively. The primal FPM with reduced flexoelectric theory is employed. The computational parameters $c_0 = \sqrt{20}$, $\eta_{11} = 1 \times 10^{10} E$, $\eta_{13} = 1 \times 10^{10} \Lambda_{11}$, $\eta_{21} = 1.0 E$, $\eta_{22} =0$, $\eta_{23} =\Lambda_{33}$. The deformation of the truncated pyramid with flexible support shows a bending component. Thus, the structure exhibits a highly inhomogeneous mechanical strain field, which results in a remarkable strain gradient and an inhomogeneous electric field, especially near the pyramid’s corners (see Fig.~\ref{fig:Ex4_Solu_1}), whereas the rigid support prevents the bending behavior of the truncated pyramid. As a result, the variation of mechanical strain is relatively smooth, and the electric potential induced at the bottom surface is much smaller than in the flexible case (see Fig.~\ref{fig:Ex4_Solu_2}). The FPM solutions are in good agreement with the results by \citet{Abdollahi2014} and \citet{Zhuang2020}.

\subsection{Elliptical hole in a plate}

\begin{figure}[htbp] 
  \centering 
    \subfigure[]{ 
    \label{fig:Ex5_Schem_01} 
    \includegraphics[width=0.48\textwidth]{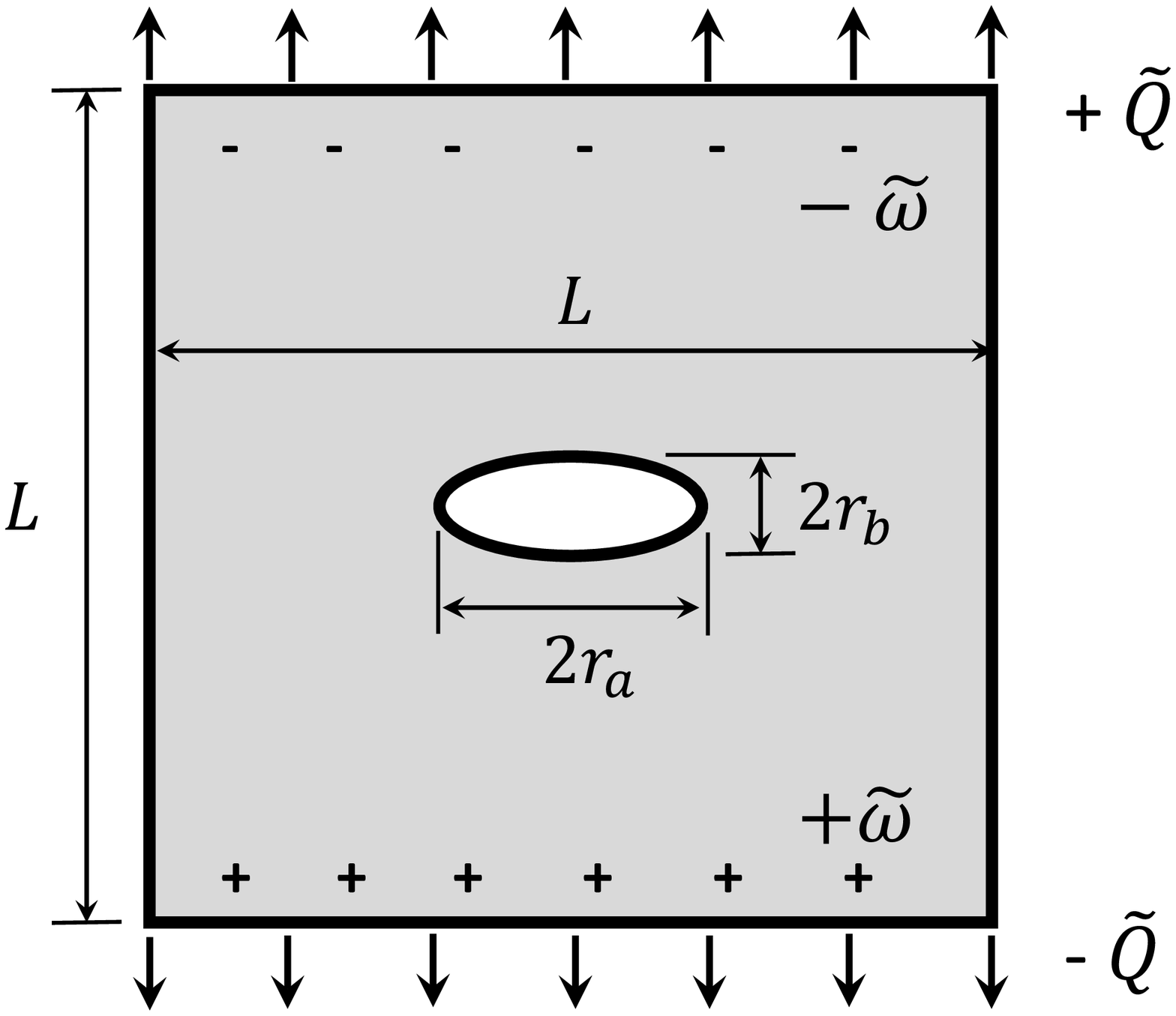}}  
    \subfigure[]{ 
    \label{fig:Ex5_Schem_02} 
    \includegraphics[width=0.48\textwidth]{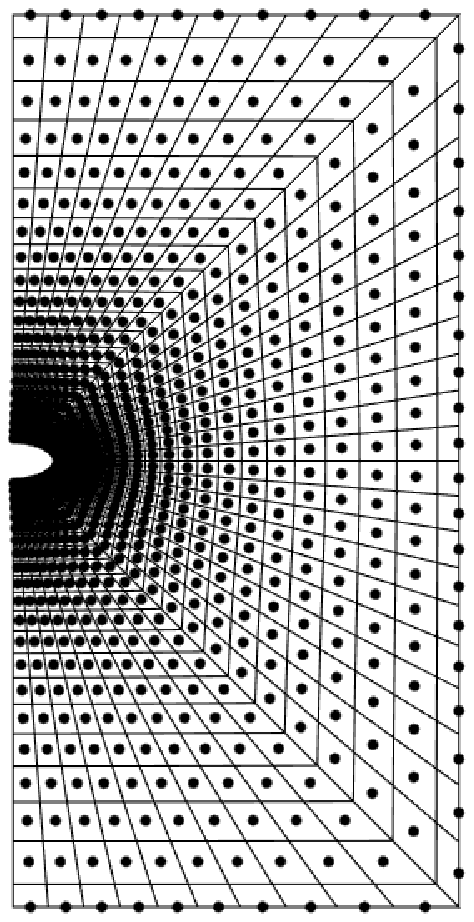}}  
  \caption{(a) A plate with an elliptical hole. (b) The point distribution and domain partition in the FPM.} 
  \label{fig:Ex1_Schem} 
\end{figure}

In the Ex.~5, we consider the flexoelectric behavior in a more complicated geometry. As shown in Fig.~\ref{fig:Ex5_Schem_01}, a cylindrical hole is in a plane strain tension field and uniform electric field. In this example, we consider a $L \times L$ square domain where $L = 20 r_a$. The plane is subjected to a uniformly distributed tensile stress on the top and bottom surfaces, i.e., $\sigma_{22} = \pm \widetilde{Q} $. And the electric field is created by the opposite surface charges $\pm \widetilde{\omega}$ on the two surfaces. The surface of the cylindrical hole is traction and charge-free.

The geometric parameters: $r_a = 10~\mathrm{ \mu m}$, $r_b = r_a / 2$. The material parameters: $E = 139~\mathrm{GPa}$, Poisson's ratio $\nu = 0.3$, $l = r_a / 3 =  3.33 \mathrm{\mu m}$, $\overline{\mu}_{11} = -6.3 \times 10^{-5}~\mathrm{C/m}$, $\overline{\mu}_{12} = 5.2 \times 10^{-6}~\mathrm{C/m}$, $\overline{\mu}_{44} = -3.4 \times 10^{-5}~\mathrm{C/m}$, $\Lambda_{11} =  \Lambda_{33} = 4.9 \times 10^{-9}~\mathrm{F/m}$, $e_{31} = e_{33} = e_{15} = 0$. The external loads: $\widetilde{Q} = E/200 = 695~\mathrm{MPa}$, $\widetilde{\omega} = 0.0837~\mathrm{C/m^2}$.

\begin{figure}[htbp] 
  \centering 
    \subfigure[]{ 
    \label{fig:Ex5_01} 
    \includegraphics[width=0.48\textwidth]{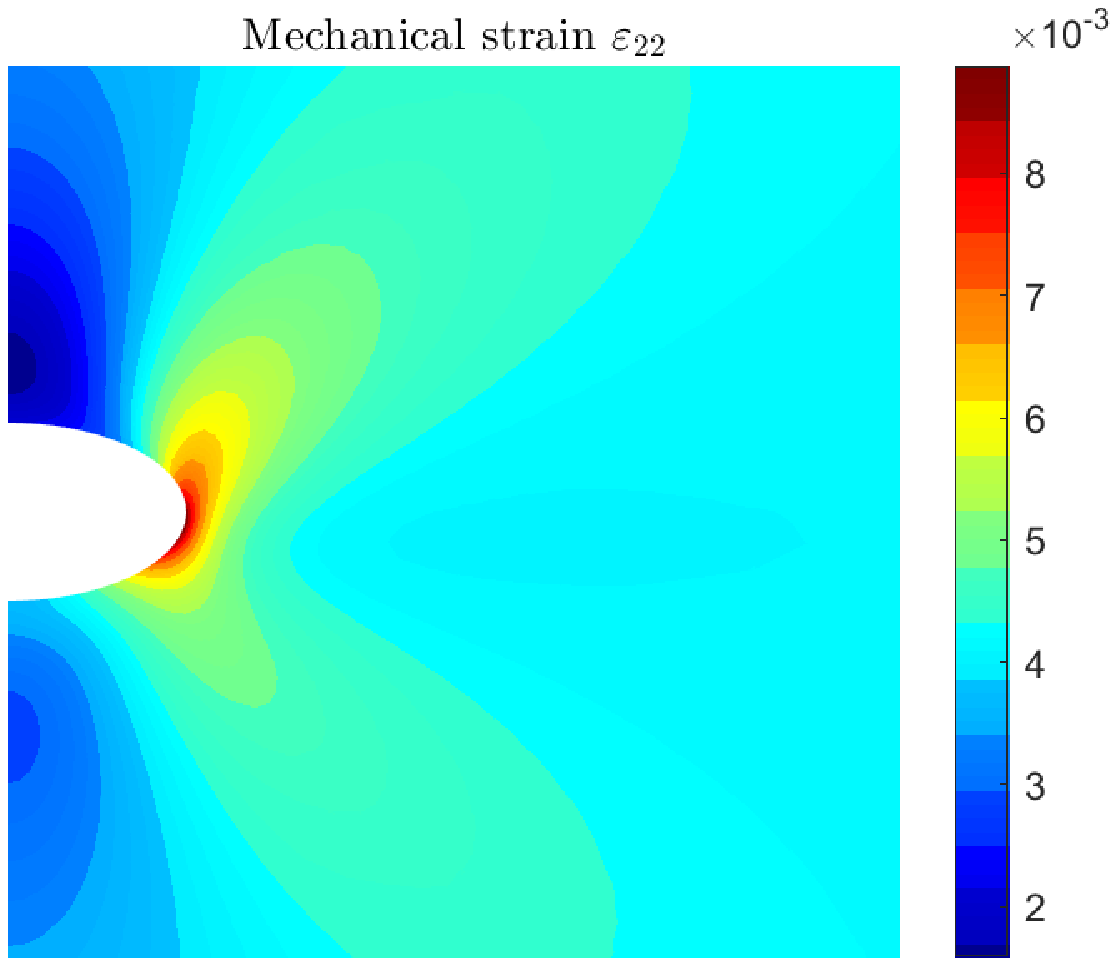}}  
    \subfigure[]{ 
    \label{fig:Ex5_02} 
    \includegraphics[width=0.48\textwidth]{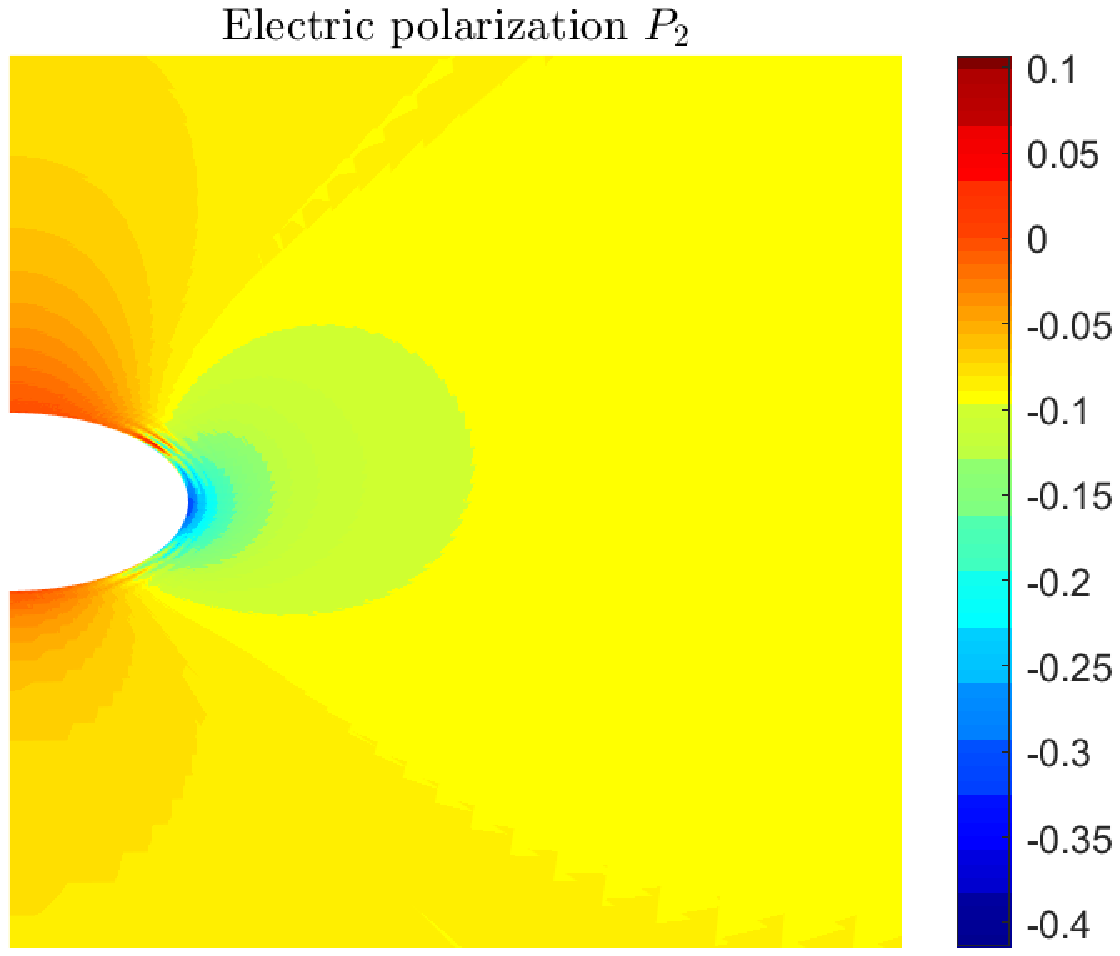}}  
  \caption{The computed solution by primal FPM for Ex.~5. (a) Distribution of the the mechanical strain $\varepsilon_{22}$. (b) Distribution of the electric polarization $P_2$.} 
  \label{fig:Ex5_Solu} 
\end{figure}

According to the symmetry, we can analyze only one half of the domain with symmetric boundary conditions. The primal FPM is employed with 3360 points. As the variation of displacement and electric polarization is expected to be more complicated, more points are distributed in the vicinity of the cylindrical hole surface (see Fig.~\ref{fig:Ex5_Schem_02}). The computational parameters: $c_0 = \sqrt{20}$, $\eta_{21} = 2.0 E$, $\eta_{22} =1.0 E$, $\eta_{23} = 0.1 \Lambda_{33}$. The computed solutions of the mechanical strain $\varepsilon_{22}$ and electric polarization $P_2$ are presented in Fig.~\ref{fig:Ex5_Solu}. The results are consistent with the work of \citet{Mao2016}. As can be seen, a concentration of the mechanical strain and electric polarization occurs in vicinity of the ``tip'’ of the hole. Besides, the system presents an asymmetric behavior along the $y$-axis, even though all the external loads are symmetric, and the material is isotropic. This implies that the direction of the polarization field is important in the analysis of flexoelectric responses.

\subsection{Influence of the electroelastic stress: A 1D example}

\begin{figure}[htbp]
    \centering
    \begin{minipage}[t]{0.48\textwidth}
        \centering
        \includegraphics[width=1\textwidth]{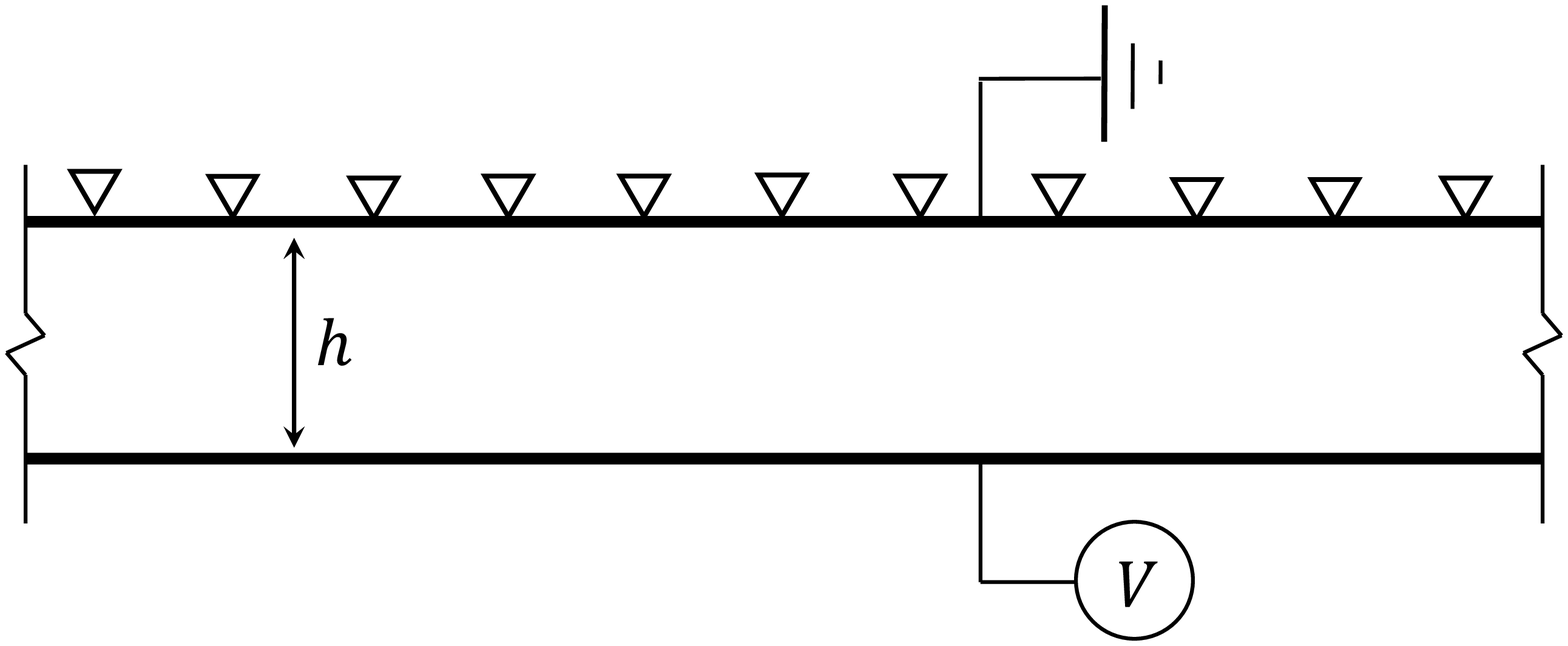}
        \caption{An infinite dielectric layer.}
        \label{fig:Ex6_Schem}
    \end{minipage}
    \begin{minipage}[t]{0.48\textwidth}
        \centering
        \includegraphics[width=1\textwidth]{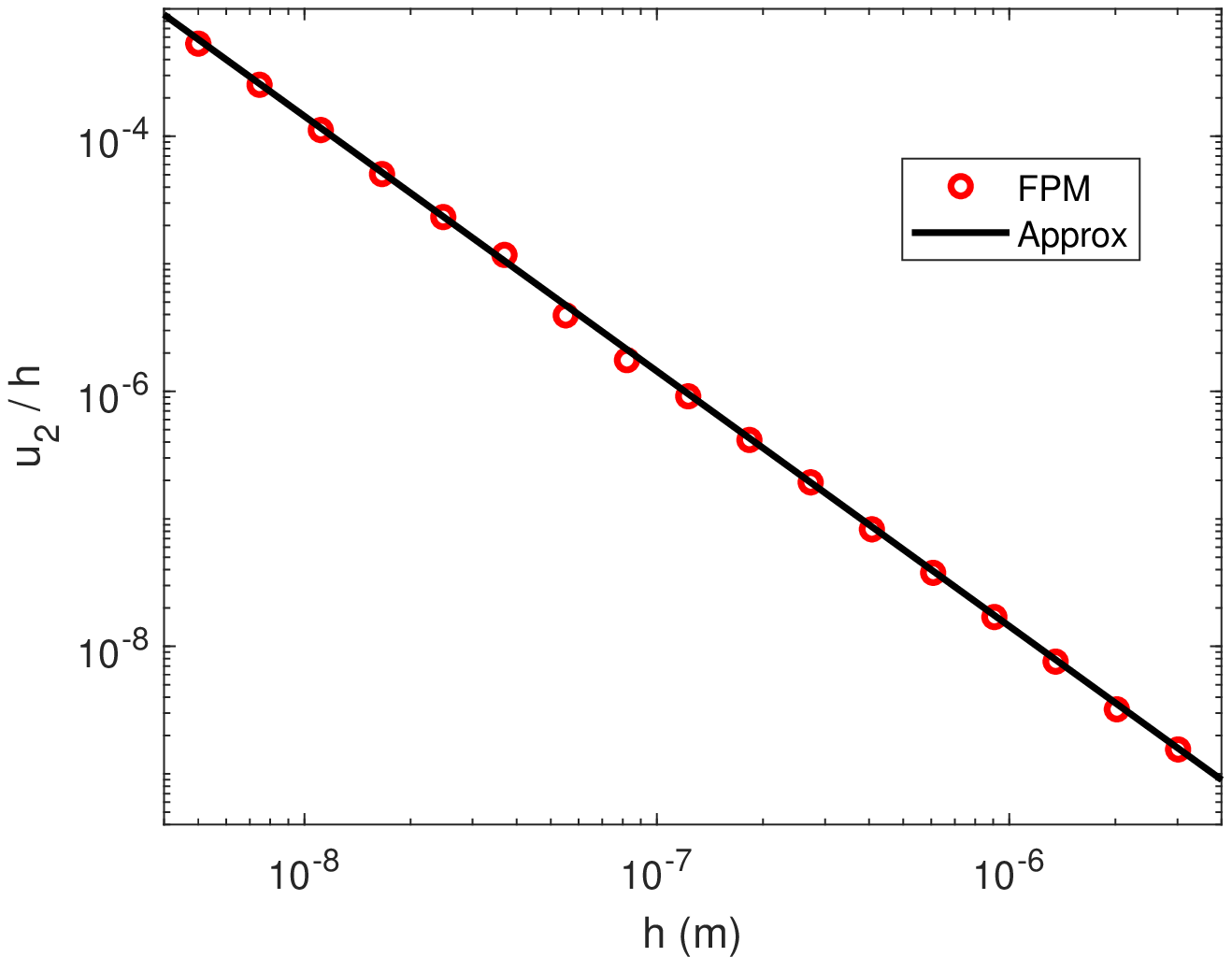}
        \caption{The vertical displacement under varying thickness $h$.}
        \label{fig:Ex6_01}
    \end{minipage}
\end{figure}

All the previous examples are analyzed in the reduced flexoelectric theory, i.e., the electrostatic stress $\boldsymbol{\sigma}^{ES}$ is omitted. However, some previous studies have pointed out that the electrostatic stress could be very strong at the nano-scale \cite{Hu2010}. Here we introduce a simple 1D example to illustrate the influence of the electrostatic stress. An infinite dielectric layer with thickness $h$ is considered. The material is isotropic and centrosymmetric cubic. The top surface is fixed and grounded, whereas the bottom surface is traction-free but subjected to a given electric potential $V$. Thus the field will lead to one-dimensional responses. In this example, the material properties are: $E = 139~\mathrm{GPa}$, Poisson's ratio $\nu = 0$, $l = 0$, $\overline{\mu}_{11} = \overline{\mu}_{12} = \overline{\mu}_{44} = 0$, $\Lambda_{33} = 1 \times 10^{-9}~\mathrm{F/m}$, $e_{31} = e_{33} = e_{15} = 0$, $\Phi_{33} = 1 \times 10^{-17}~\mathrm{F \cdot m}$, $b_{33} = 1 \times 10^{-4}~\mathrm{C/m}$. Note that the high-order electric effect is also taken into consideration in this example. The primal FPM based on full flexoelectric theory is used with 20 points distributed uniformly in $y$-direction. Since the FPM formula based on the full theory is nonlinear, the Newton-Raphson solution method is used. The relative tolerance of the displacement and electric potential vectors is set to be $1 \times 10^{-6}$. In this example, it usually takes 1-5 iterations to reach the actual solution. We especially concentrated on the solution of vertical displacement $u_2$ at the bottom surface. An approximate solution for this one-dimensional behavior is given in \cite{Hu2010}. When the thickness of the layer varies between nano and micro scales, the distribution of the displacement is roughly linear. Fig.~\ref{fig:Ex6_01} presents the normalized displacement at the bottom $u_2 / h$ under varying thickness $h$. Generally, when the thickness of the layer decreases, the normalized displacement at the bottom $u_2 / h$ increases inversely proportional to $h^2$. This result agrees well with the approximate equations \cite{Hu2010}. This example implies that the electrostatic stress, as well as the influence of the electric field gradient could be significant for problems at nano scale. Therefore, the full flexoelectric theory is highly recommended to be applied in analyzing such systems.

\subsection{Influence of the electroelastic stress: A hollow cylinder}

\begin{figure}[htbp] 
  \centering 
    \subfigure[]{ 
    \label{fig:Ex1_Comp_01} 
    \includegraphics[width=0.48\textwidth]{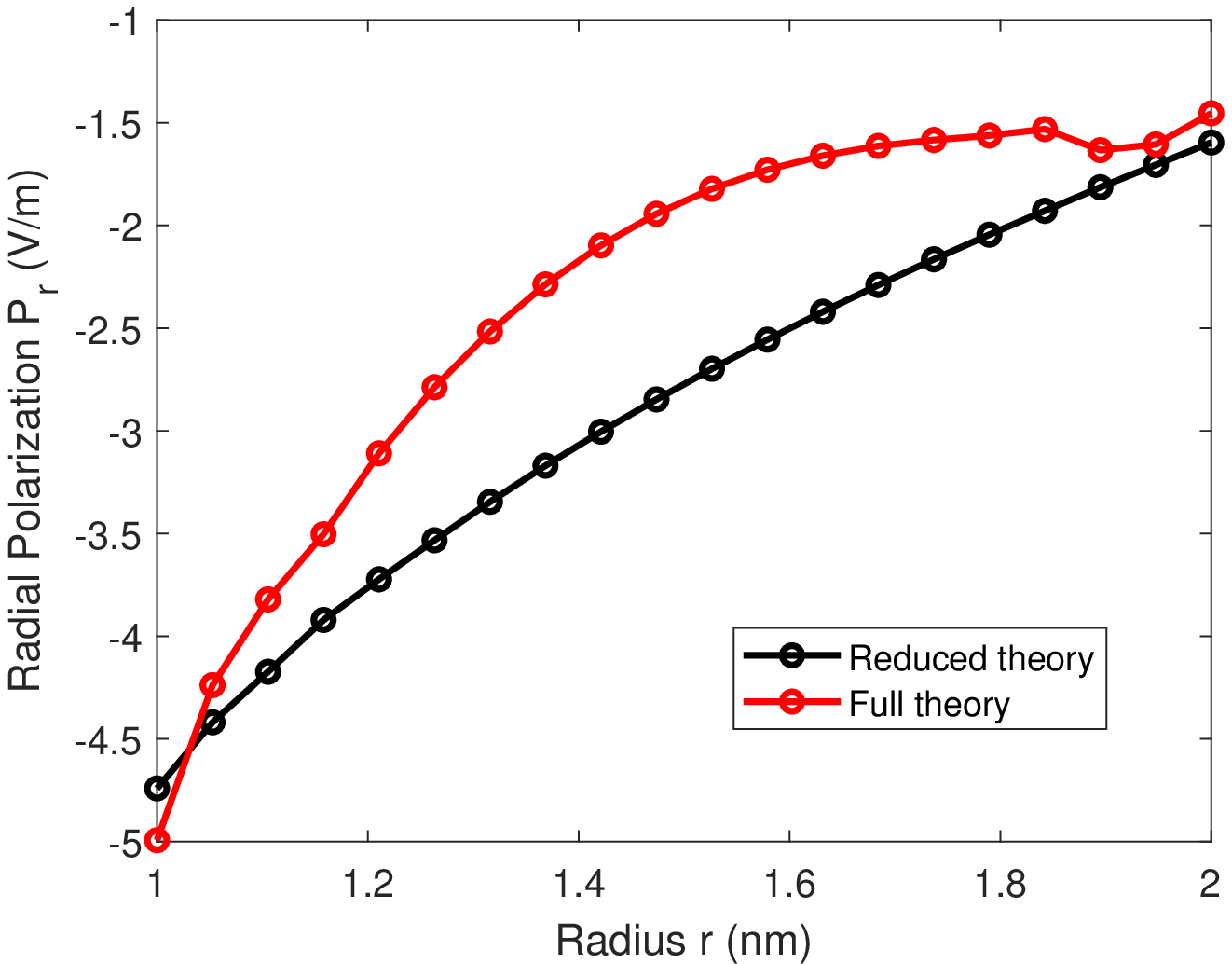}}  
    \subfigure[]{ 
    \label{fig:Ex1_Comp_02} 
    \includegraphics[width=0.48\textwidth]{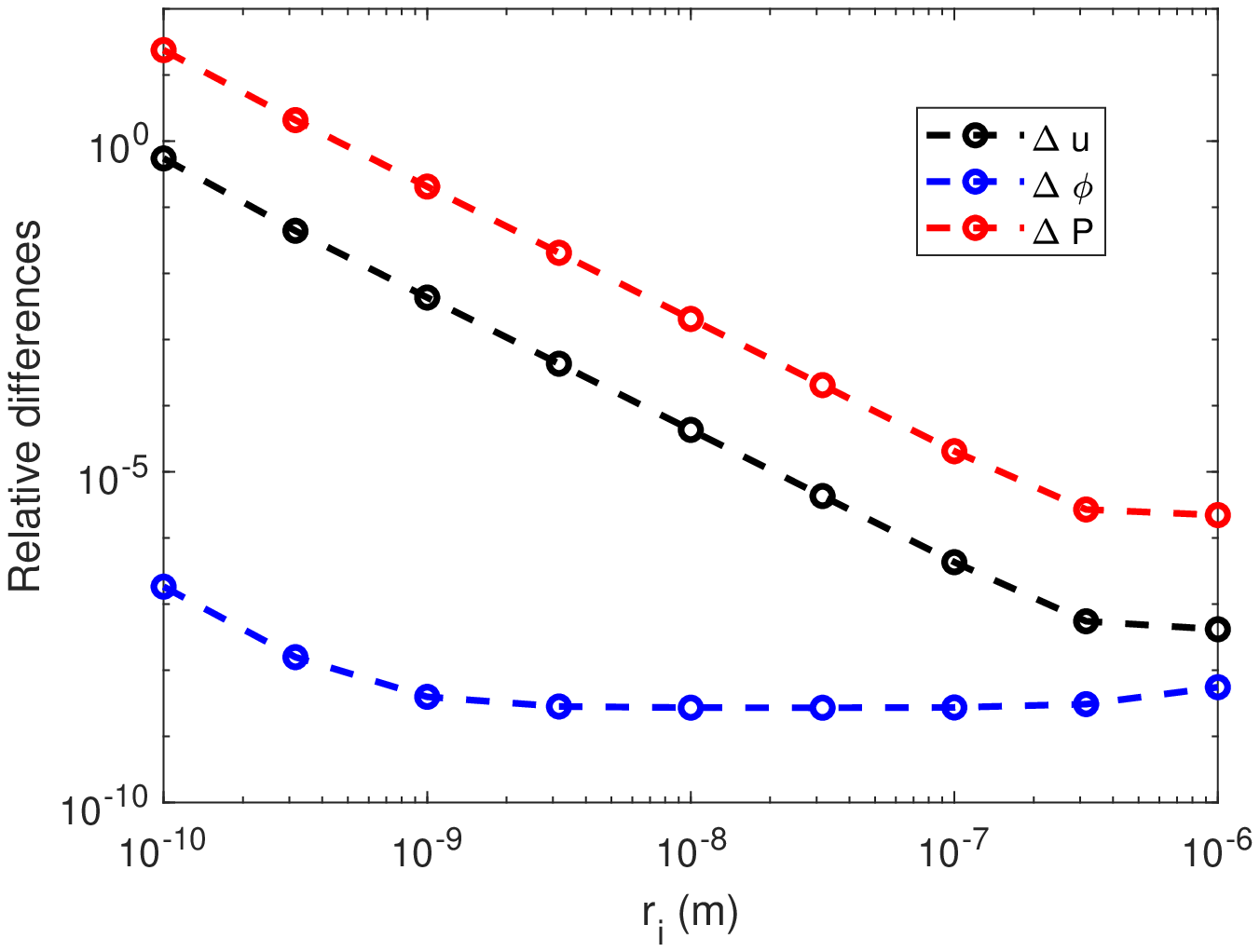}}  
  \caption{(a) The distribution of the radial electric polarization $P_r$ under full and reduced flexoelectric theories when $r_i = 1~\mathrm{nm}$. (b) The relative differences of displacement $\mathbf{u}$, electric potential $\phi$ and electric polarization $\mathbf{P}$ between the two theories under varying inner radius $r_i$.} 
  \label{fig:Ex1_Comp} 
\end{figure}

Furthermore, we reconsider the first example but miniaturize its geometric size. All the material properties and the given electric potentials remain the same. Yet the radii of the cylinder $r_i$ and $r_o$, as well as the fixed displacements $u_i$ and $u_o$, decrease analogously. In Ex.~1, when $r_i = 10~\mathrm{\mu m}$, as has been stated, the solutions based on full and reduced flexoelectric theories are approximately the same. However, when the inner radius $r_i$ decreases to $1~\mathrm{nm}$, as shown in Fig.~\ref{fig:Ex1_Comp_01}, a significant difference occurs between the two theories. This mainly results from the electrostatic behavior. The relative differences of the two theories are defined the same as Eqn.~\ref{eqn:error} where the exact and computed solutions are replaced by solutions based on the full and reduced theories respectively. Fig.~\ref{fig:Ex1_Comp_02} demonstrates the variation of the relative differences of the displacement, electric potential and electric polarizations between the two theories when the geometric size of the system changes. While the relative difference of the displacement field remains small, the differences of the electric potential and electric polarization become enormous when $r_i$ drops below $0.1~\mathrm{\mu m}$. In fact, the influence of the electrostatic stress can be dominant at the nano-scale. The same as in the previous 1D example, we can conclude that the full theory which takes the electrostatic stress and electric field gradient into consideration is recommended to be applied in analyzing nano-scale systems.

\section{Analysis of stationary cracks} \label{sec:crack_1}

In this section, a number of numerical examples with stationary cracks are considered. The material PZT-5H is under study in all the following examples. The material parameters are shown in \ref{app:PZT-5H}, in which $\alpha$ is a scaling factor introduced here to assess the influence of the strain and electric field gradients. The poling direction of the material is always parallel to $y$-axis.

\subsection{A square plate with a central crack}

\begin{figure}[htbp] 
  \centering 
    \subfigure[]{ 
    \label{fig:Ex21_Schem_01} 
    \includegraphics[width=0.48\textwidth]{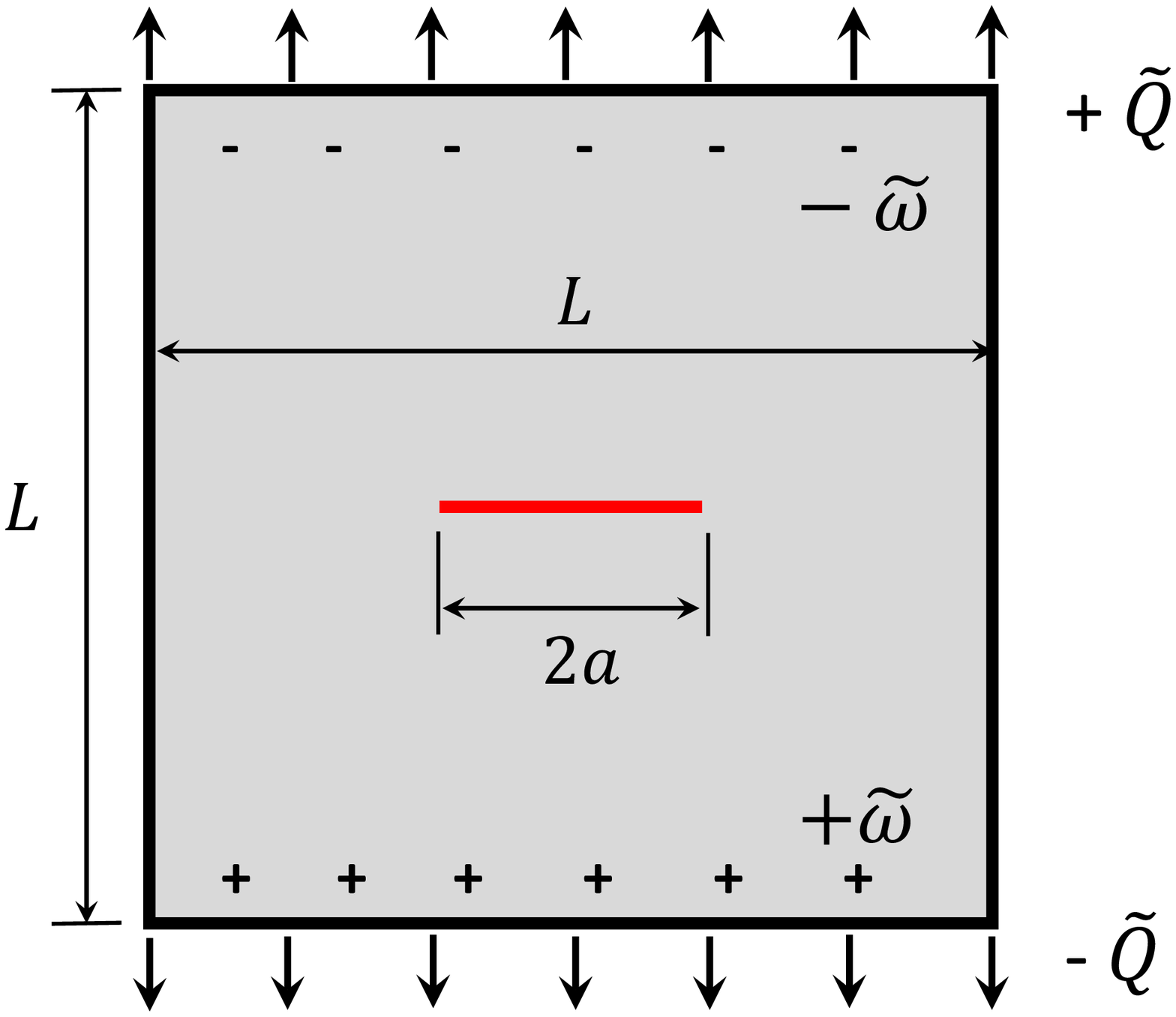}}  
    \subfigure[]{ 
    \label{fig:Ex21_Schem_02} 
    \includegraphics[width=0.48\textwidth]{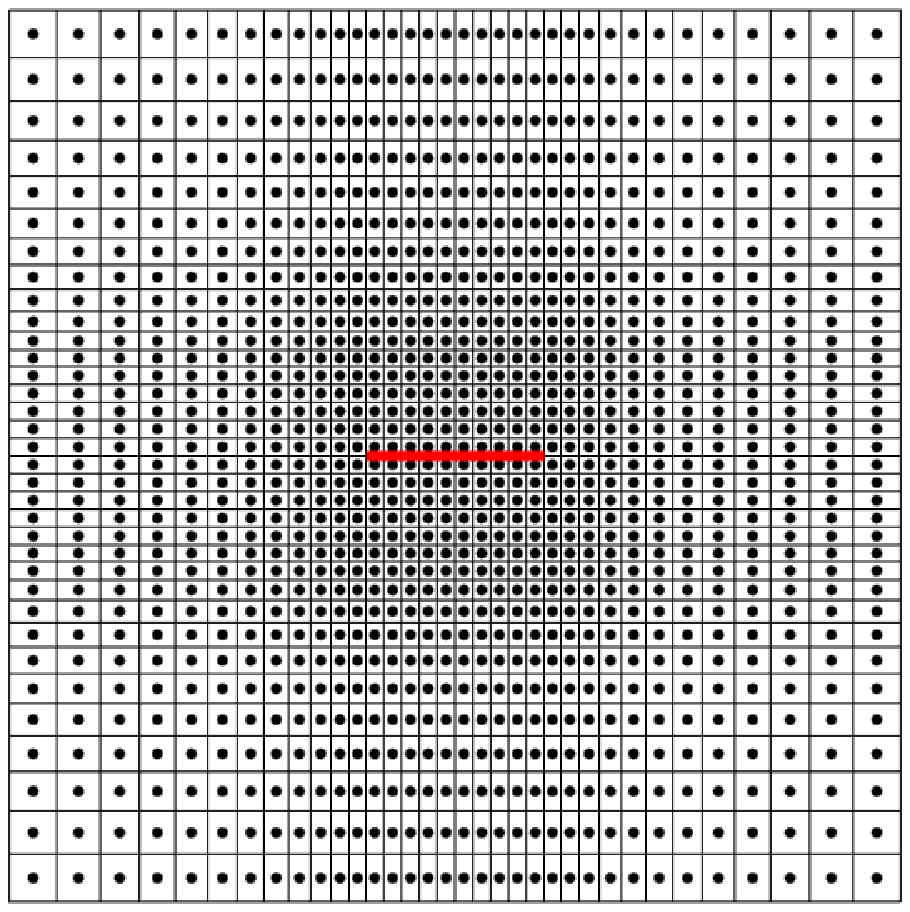}}  
  \caption{(a) A square plate with a central crack.  (b) The point distribution and domain partition in the FPM.} 
  \label{fig:Ex21_Schem} 
\end{figure}

First, a square plate with a central crack is considered (see Fig.~\ref{fig:Ex21_Schem_01}). The geometric parameters: $L = 1~\mathrm{\mu m}$, $a =  0.1~\mathrm{\mu m}$. The plate is subjected to a mechanical tensile stress $\pm \widetilde{Q} = 1.17~\mathrm{MPa}$. The external electrical loading is given as opposite surface charges $D_i n_i = \pm \widetilde{\omega}$ on the top and bottom surfaces. Three different loadings are considered: $D1 = -5 \times 10^{-4}~\mathrm{C/m^2}$, $D2 = 0$ and $D3 = 1 \times 10^{-3}~\mathrm{C/m^2}$. The crack is electrically impermeable.

As shown in Fig.~\ref{fig:Ex21_Schem_02}, a total of 2916 points and the same number of subdomains are used. For the internal boundaries on the crack, the interactions between its neighboring points are eliminated, and the cracked internal boundary will be converted in to traction-free boundaries. We analyze the problem with primal FPM based on full flexoelectric theory. Computational parameters: $c_0 = \sqrt{20}$, $\eta_{21} = 2.0 E$, $\eta_{22} =50 E$, $\eta_{23} = 2.0 \Lambda_{33}$, $\eta_{24} = 0$. The distributions of the displacement $u_2$ and electric potential $\phi$ on the upper crack surface under different electrical loadings ($D1$, $D2$ and $D3$) with $\alpha = 0, 2, 4$ are presented in Fig.~\ref{fig:Ex21_Solu}. As can be seen, the electrical loading has a considerable influence on both the displacement and electric potential results on the crack surface. The strain and electric field gradients (assessed by the size-factor $\alpha$), on the other hand, can decrease the displacement variation and enhance the electric potential in the same time.

\begin{figure}[htbp] 
  \centering 
    \subfigure[$\alpha = 0$]{ 
    \label{fig:Ex21_01_01} 
    \includegraphics[width=0.48\textwidth]{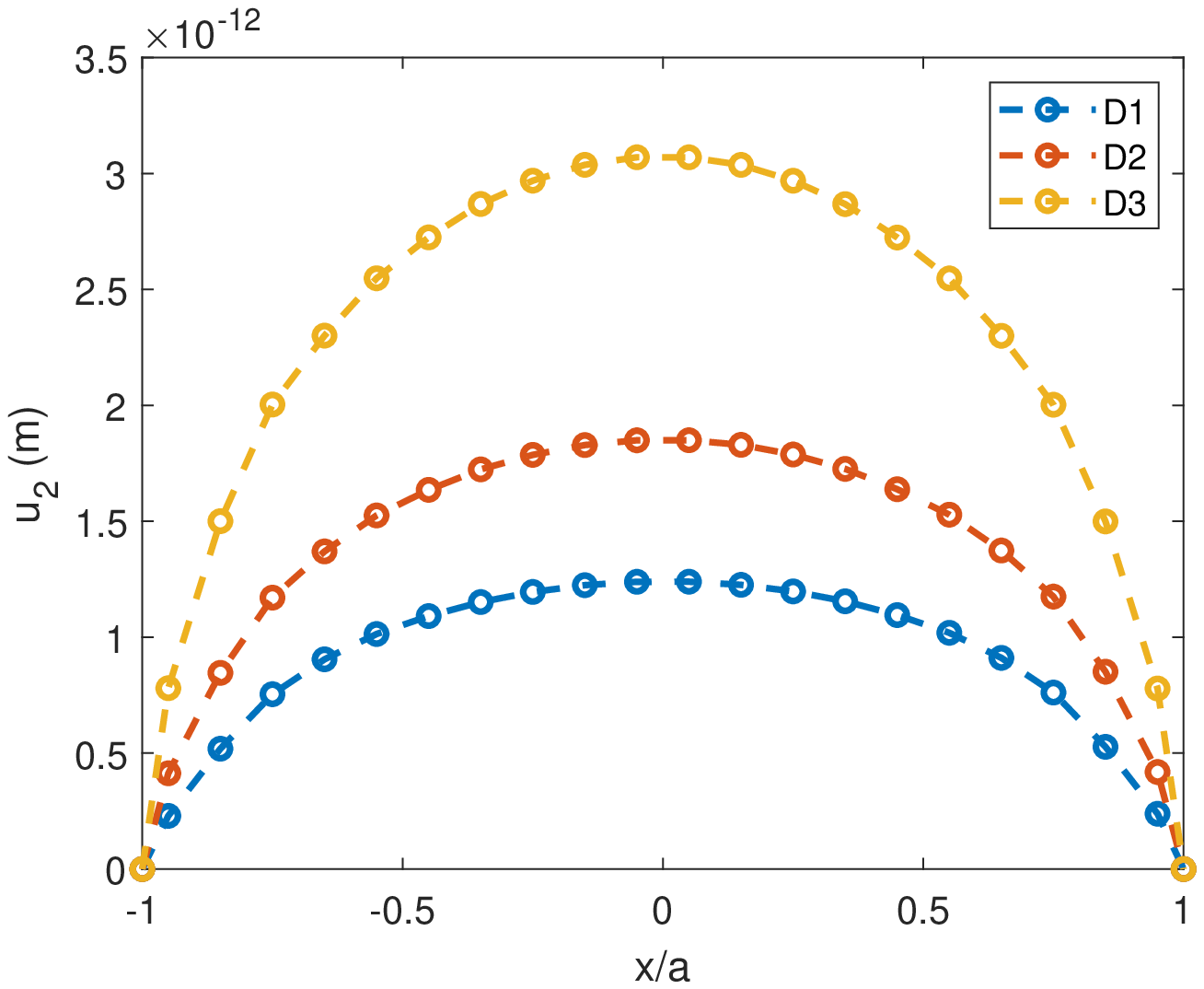}}  
    \subfigure[$\alpha = 0$]{ 
    \label{fig:Ex21_01_02} 
    \includegraphics[width=0.48\textwidth]{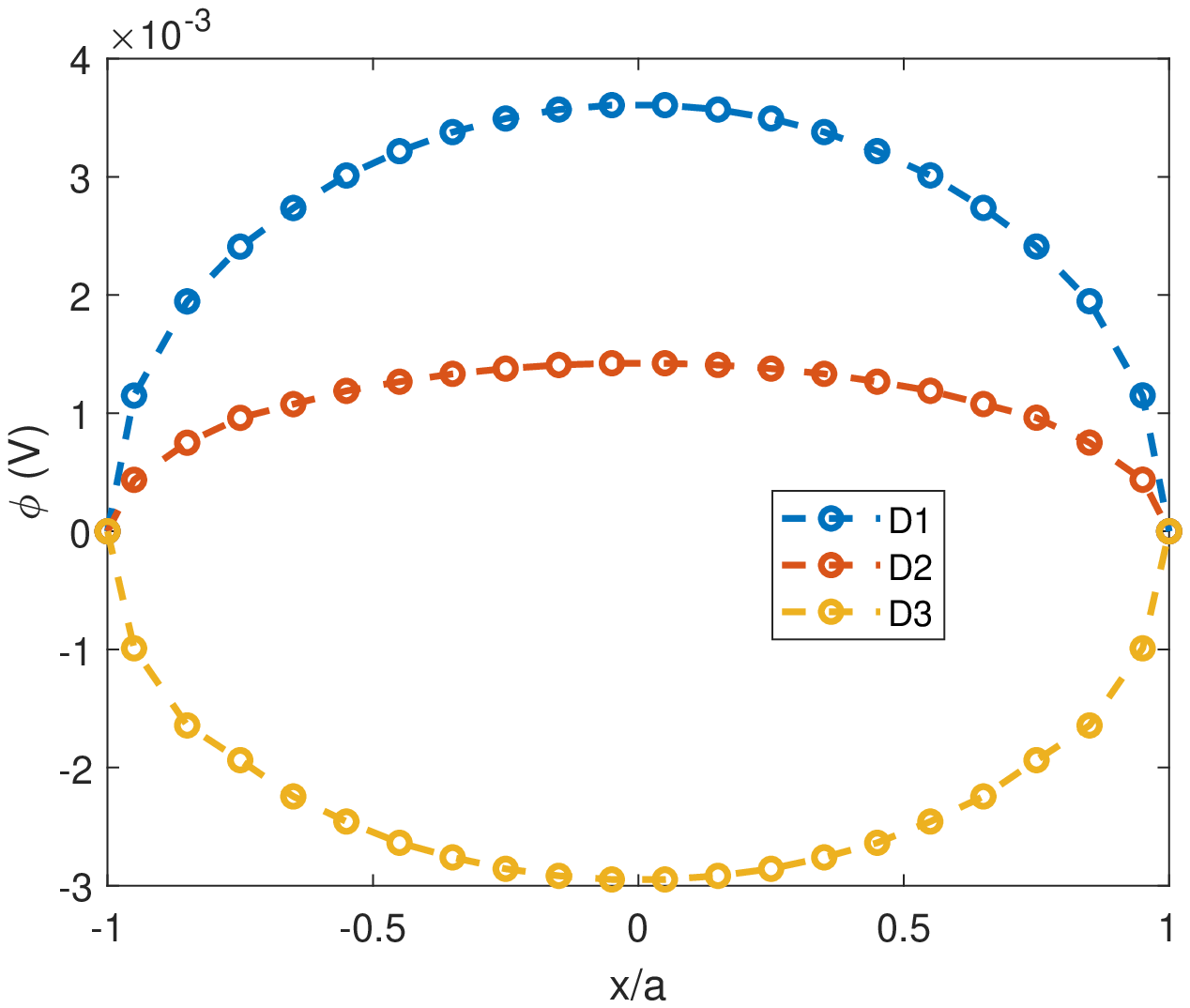}}  
        \subfigure[$\alpha = 2$]{ 
    \label{fig:Ex21_02_01} 
    \includegraphics[width=0.48\textwidth]{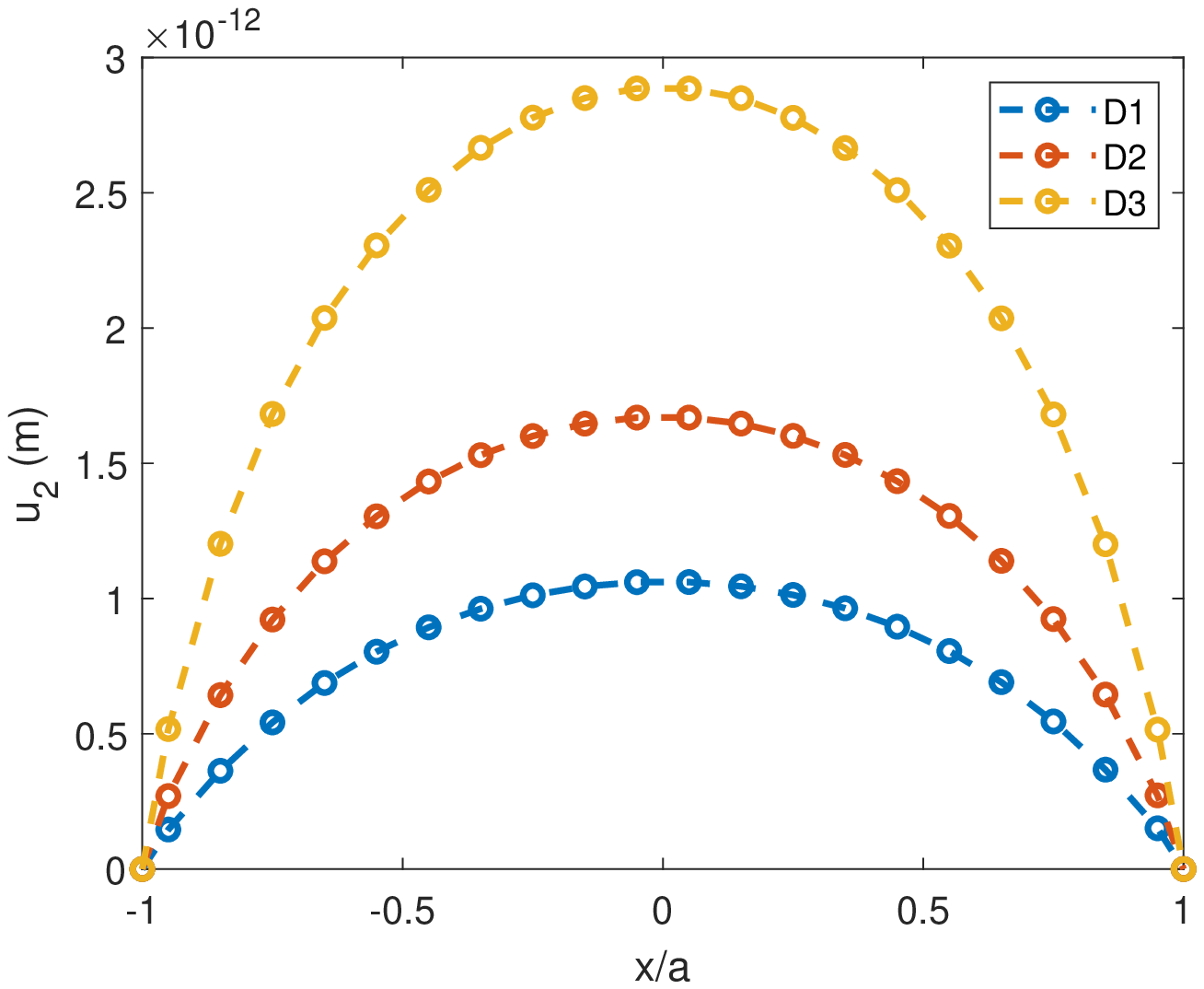}}  
    \subfigure[$\alpha = 2$]{ 
    \label{fig:Ex21_02_02} 
    \includegraphics[width=0.48\textwidth]{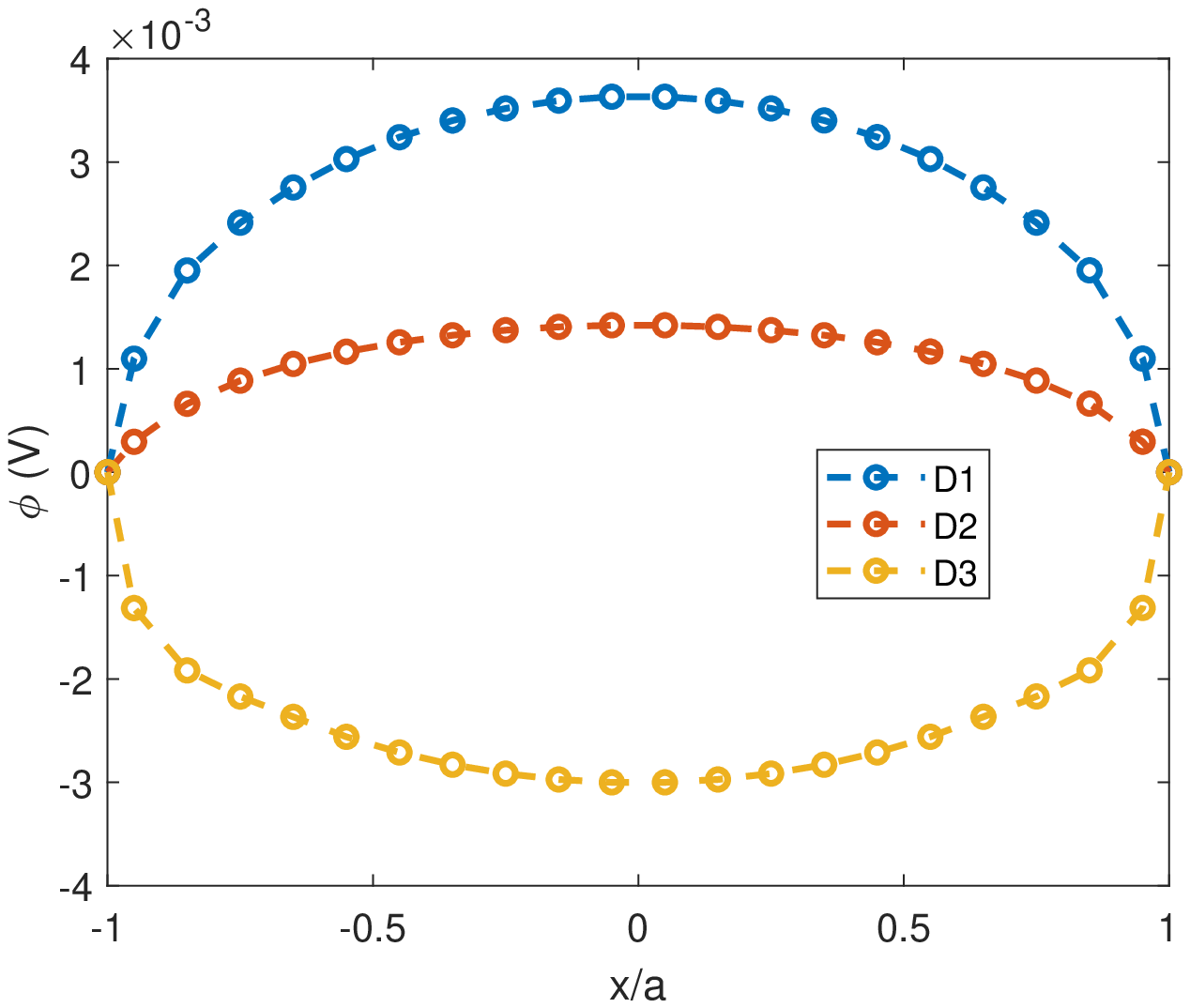}}  
        \subfigure[$\alpha = 4$]{ 
    \label{fig:Ex21_03_01} 
    \includegraphics[width=0.48\textwidth]{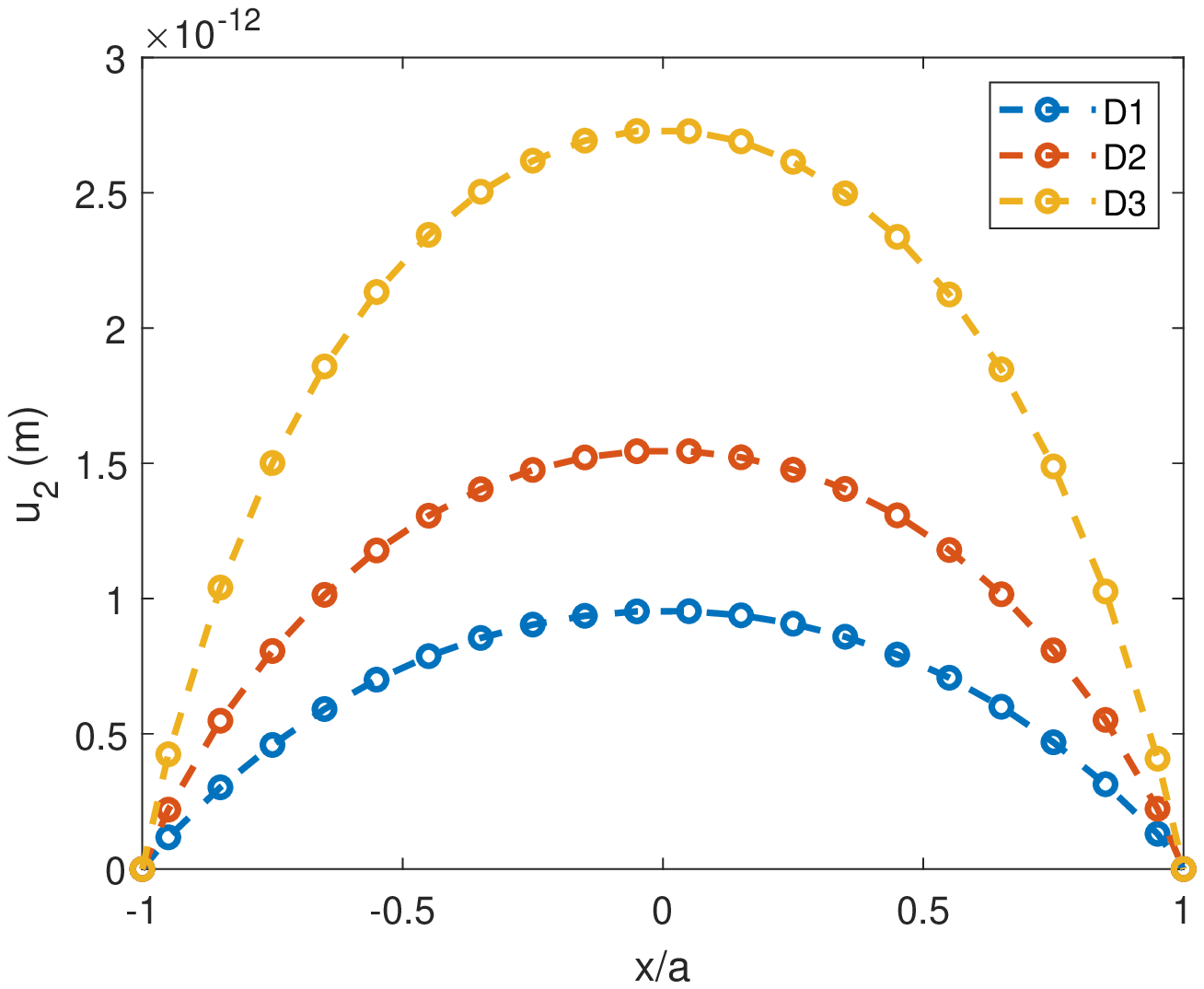}}  
    \subfigure[$\alpha = 4$]{ 
    \label{fig:Ex21_03_02} 
    \includegraphics[width=0.48\textwidth]{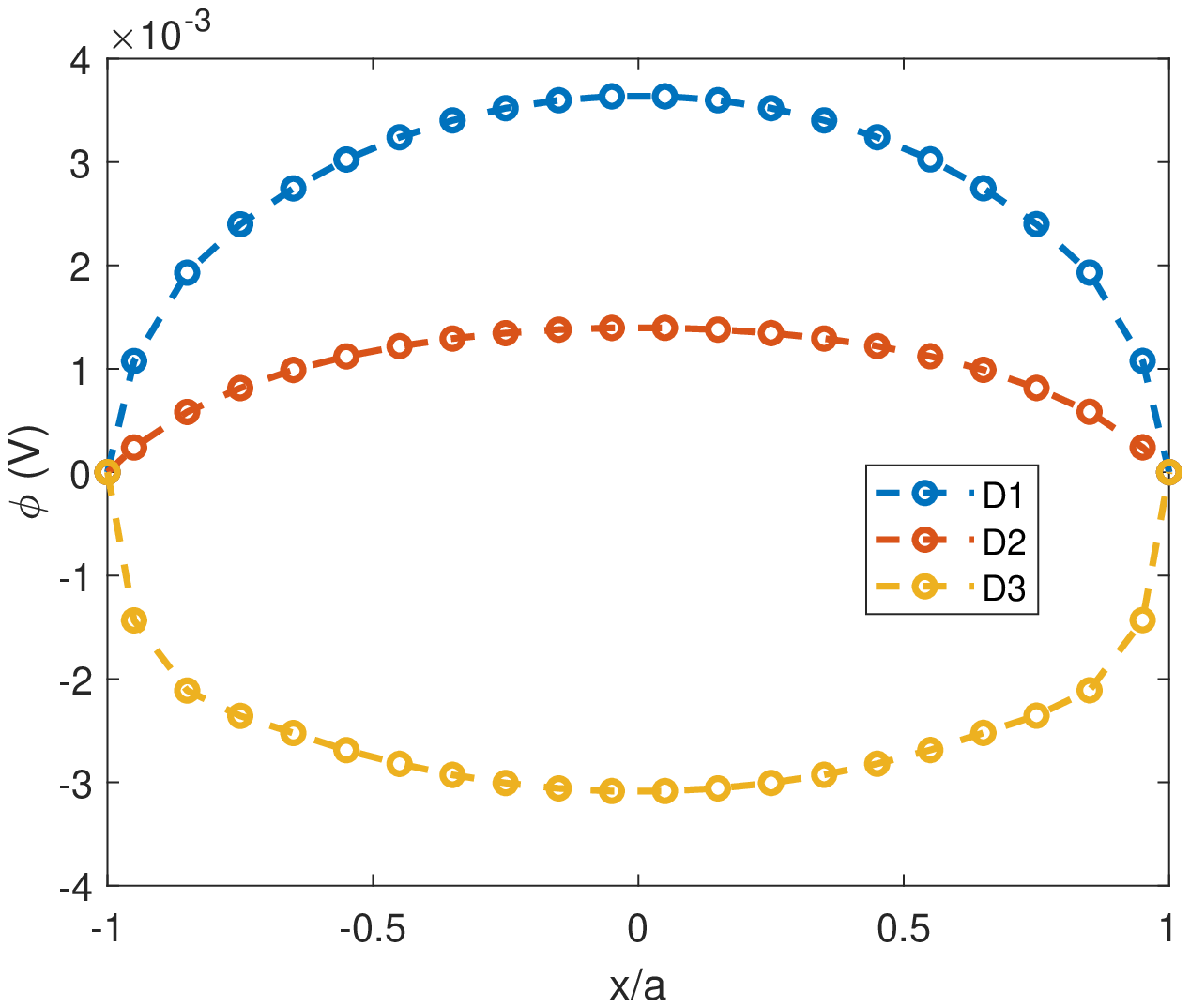}}  
  \caption{Distribution of displacement $u_2$ and electric potential $\phi$ on the upper crack-face.} 
  \label{fig:Ex21_Solu} 
\end{figure}

\begin{figure}[htbp] 
  \centering 
    \subfigure[D1]{ 
    \label{fig:Ex22_01_01} 
    \includegraphics[width=0.48\textwidth]{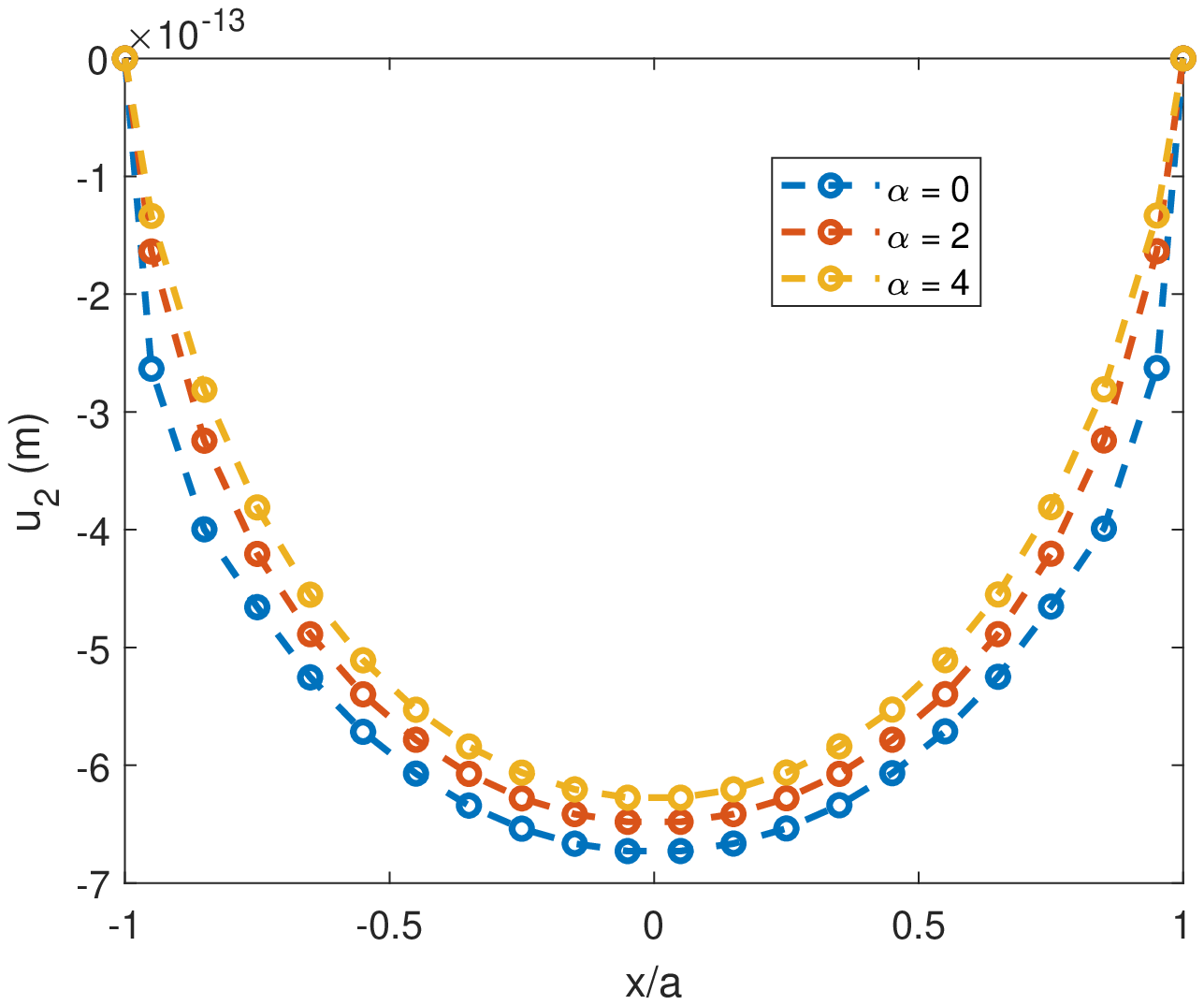}}  
    \subfigure[D1]{ 
    \label{fig:Ex22_01_02} 
    \includegraphics[width=0.48\textwidth]{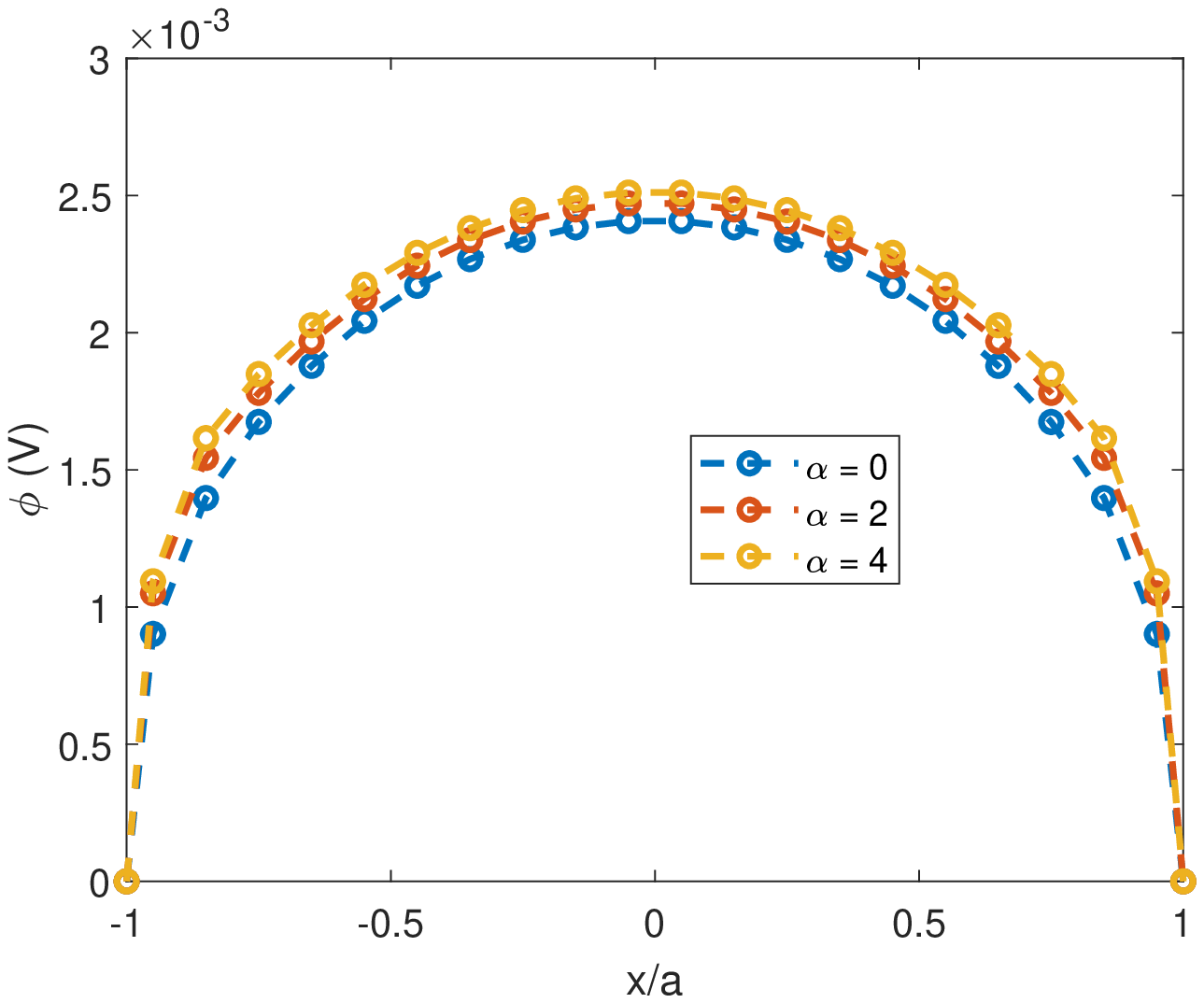}}  
        \subfigure[D3]{ 
    \label{fig:Ex22_02_01} 
    \includegraphics[width=0.48\textwidth]{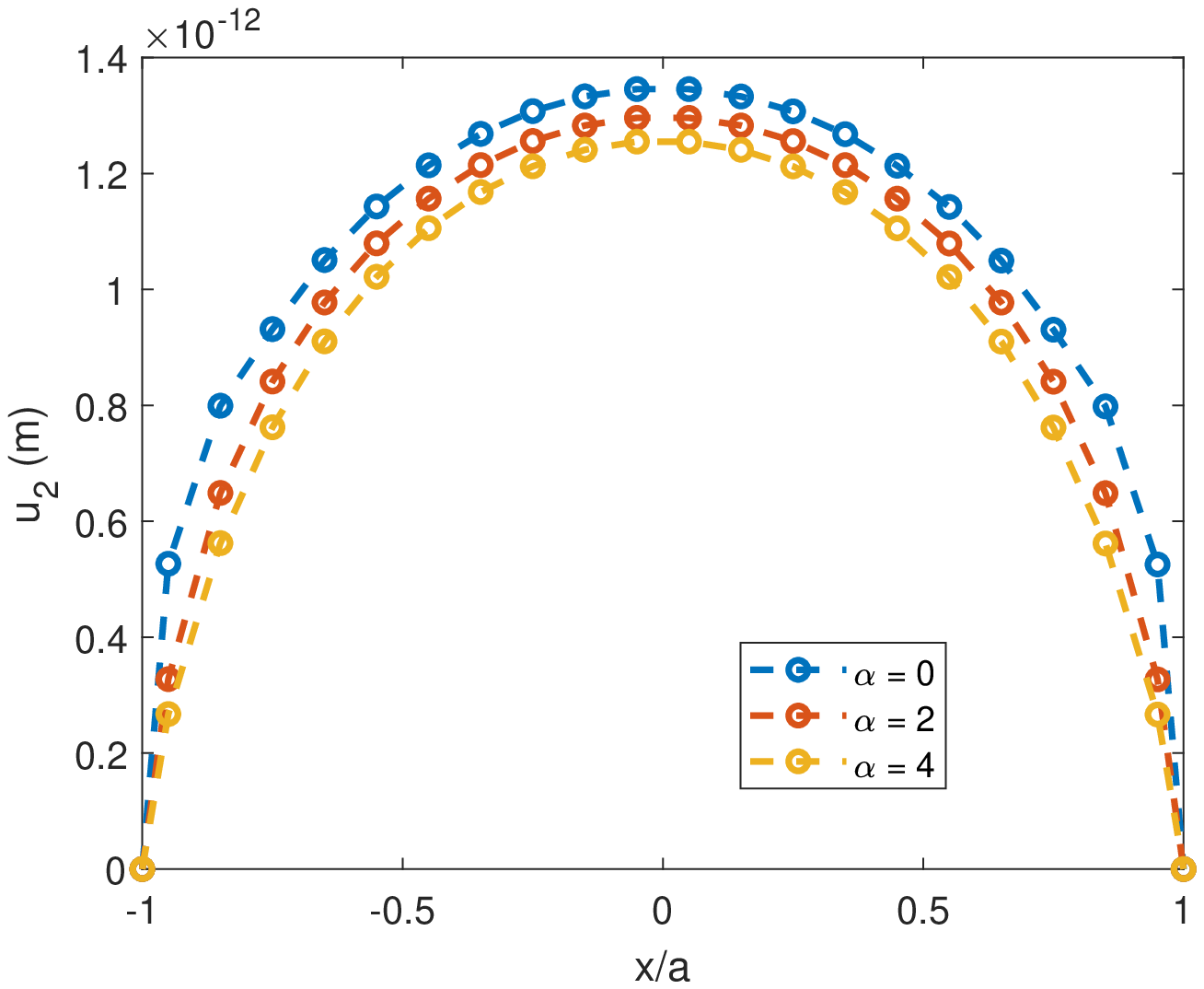}}  
    \subfigure[D3]{ 
    \label{fig:Ex22_02_02} 
    \includegraphics[width=0.48\textwidth]{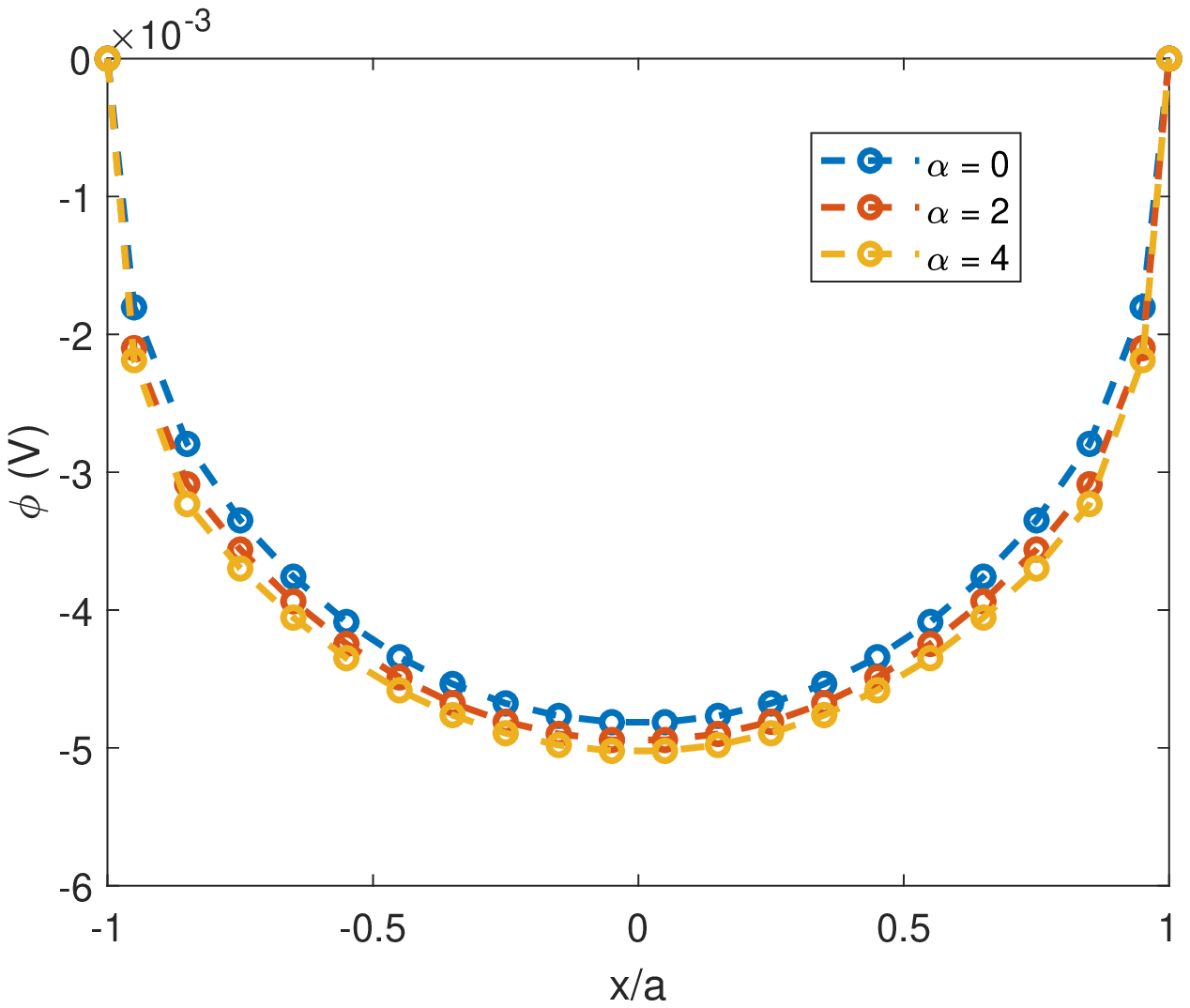}}   
  \caption{Distribution of displacement $u_2$ and electric potential $\phi$ on the upper crack-face. (a) (b) $D1 = -5 \times 10^{-4}~\mathrm{C/m^2}$. (c) (d) $D3 = 1 \times 10^{-3}~\mathrm{C/m^2}$} 
  \label{fig:Ex21_Solu} 
\end{figure}

Next, for the plate shown in Fig.~\ref{fig:Ex21_Schem_01}, we simply concentrate on the electrical loading, i.e., the tensile stress $\widetilde{Q}$ is absent. The solutions for the crack-opening-displacement $u_2$ and electric potential $\phi$ are shown in Fig.~\ref{fig:Ex21_Solu}. A change of sign in the electric loading will result in a change of sign in the induced displacement and electric potential on the crack. Note that the negative $u_2$ on the crack shown in Fig.~\ref{fig:Ex22_01_01} is not real, which implies that a surface contact approach should be included for further studies. When comparing the results of $\alpha = 0, 2$ and $4$, we can conclude that in this example, the strain and electric field gradient effects can reduce or ``smooth'’ the mechanical deformation, and enlarge the electric response on the crack.

\subsection{A plate with two symmetric cracks on the boundary of a central hole}

\begin{figure}[htbp] 
  \centering 
    \subfigure[]{ 
    \label{fig:Ex23_Schem_01} 
    \includegraphics[width=0.48\textwidth]{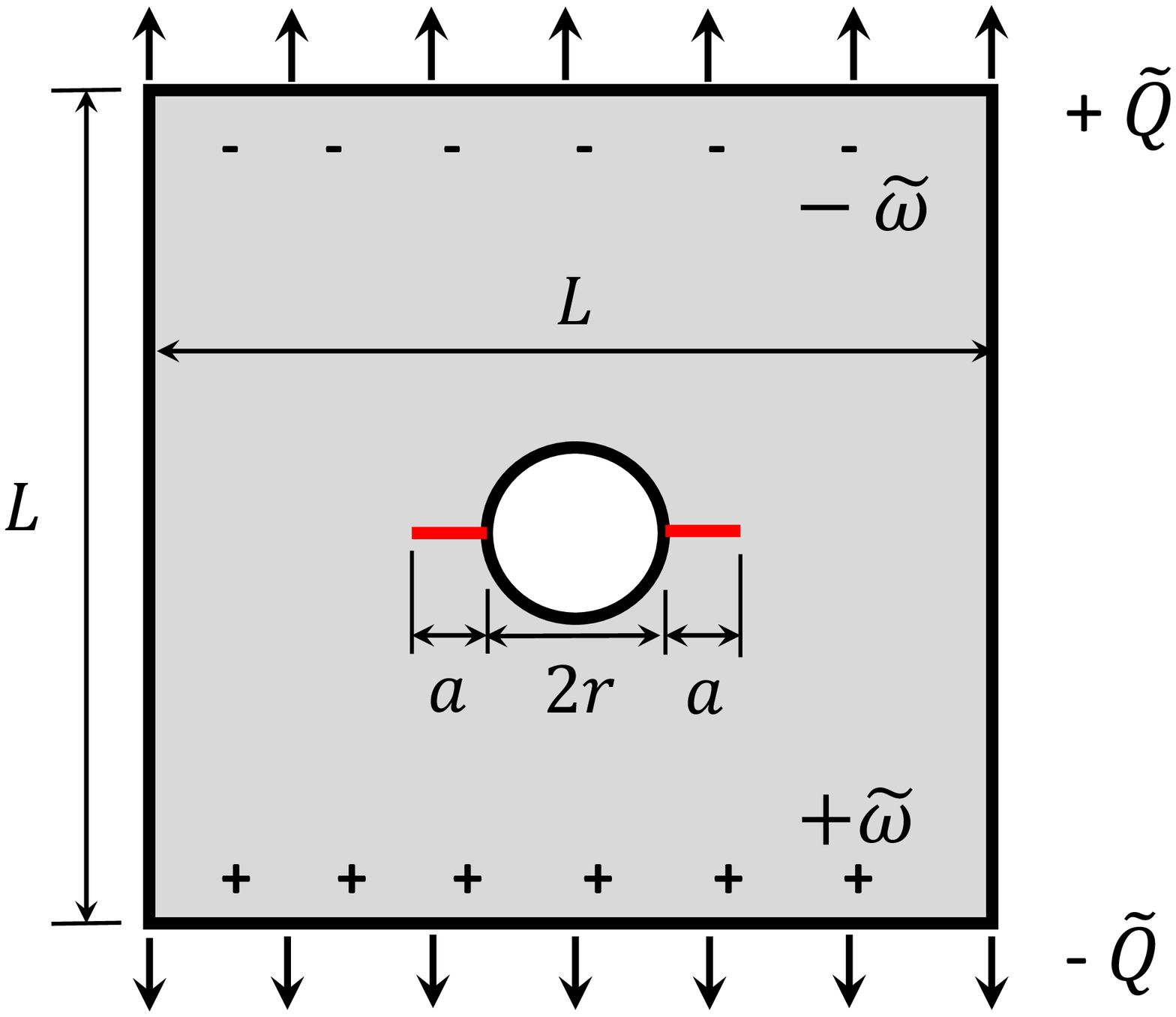}}  
    \subfigure[]{ 
    \label{fig:Ex23_Schem_02} 
    \includegraphics[width=0.48\textwidth]{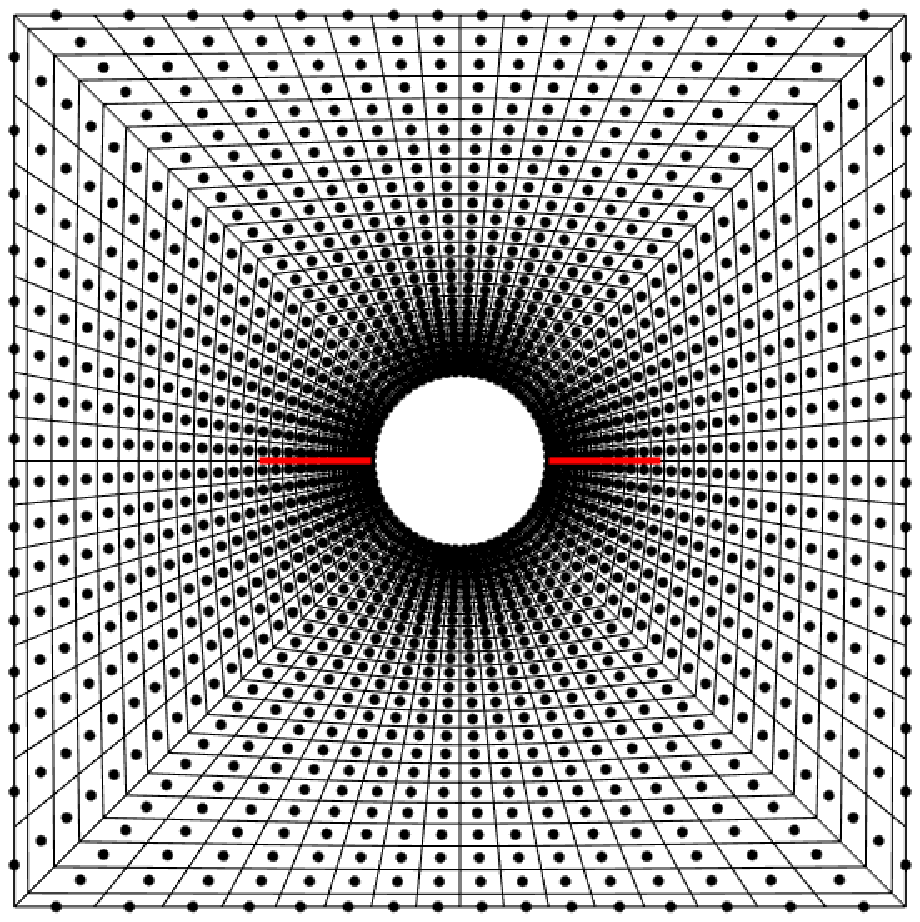}}  
  \caption{(a) A plate with two symmetric cracks on the boundary of a central hole.  (b) The point distribution and domain partition in the FPM.} 
  \label{fig:Ex23_Schem} 
\end{figure}

\begin{figure}[htbp] 
  \centering 
    \subfigure[D1]{ 
    \label{fig:Ex23_01_01} 
    \includegraphics[width=0.48\textwidth]{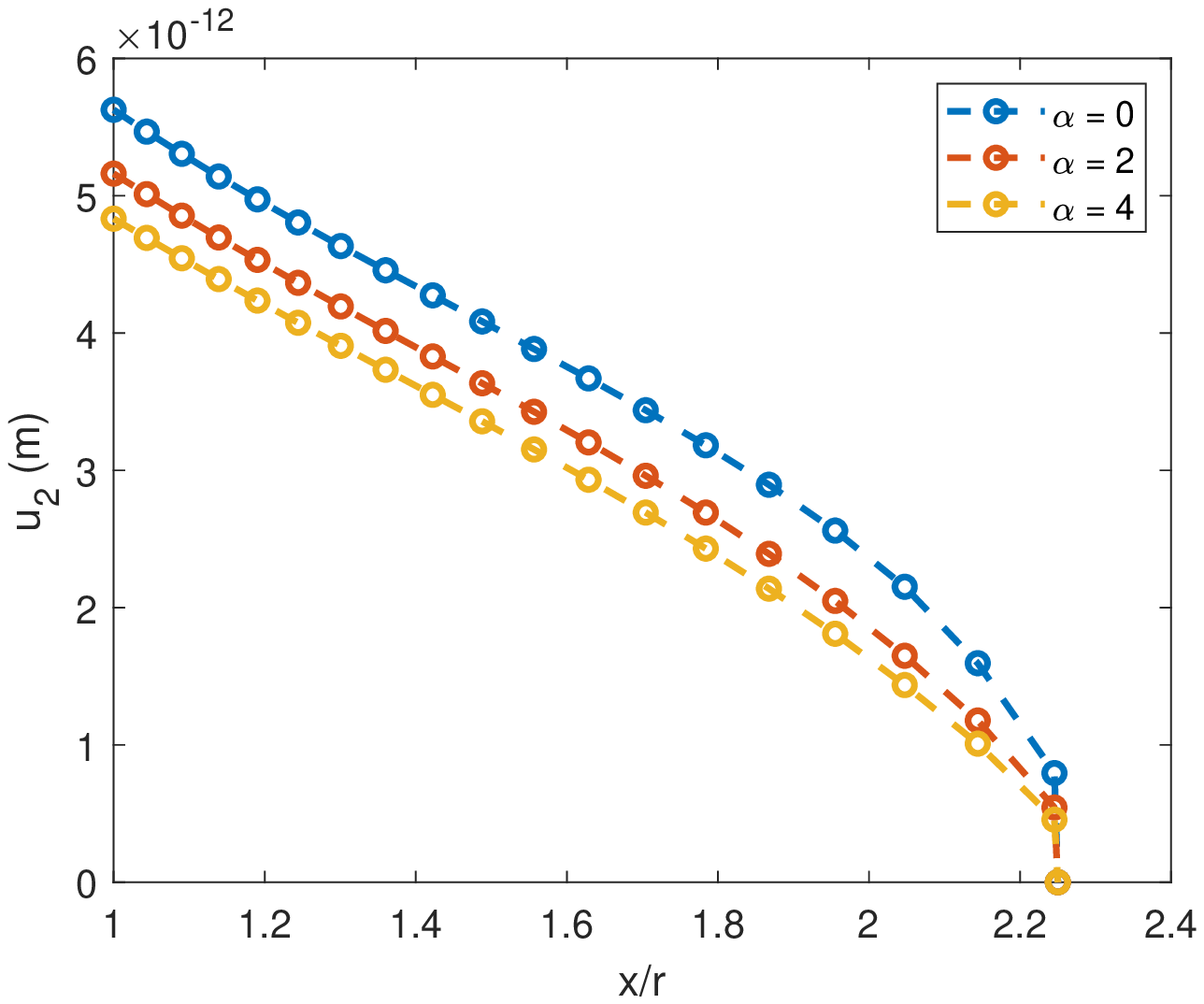}}  
    \subfigure[D1]{ 
    \label{fig:Ex23_01_02} 
    \includegraphics[width=0.48\textwidth]{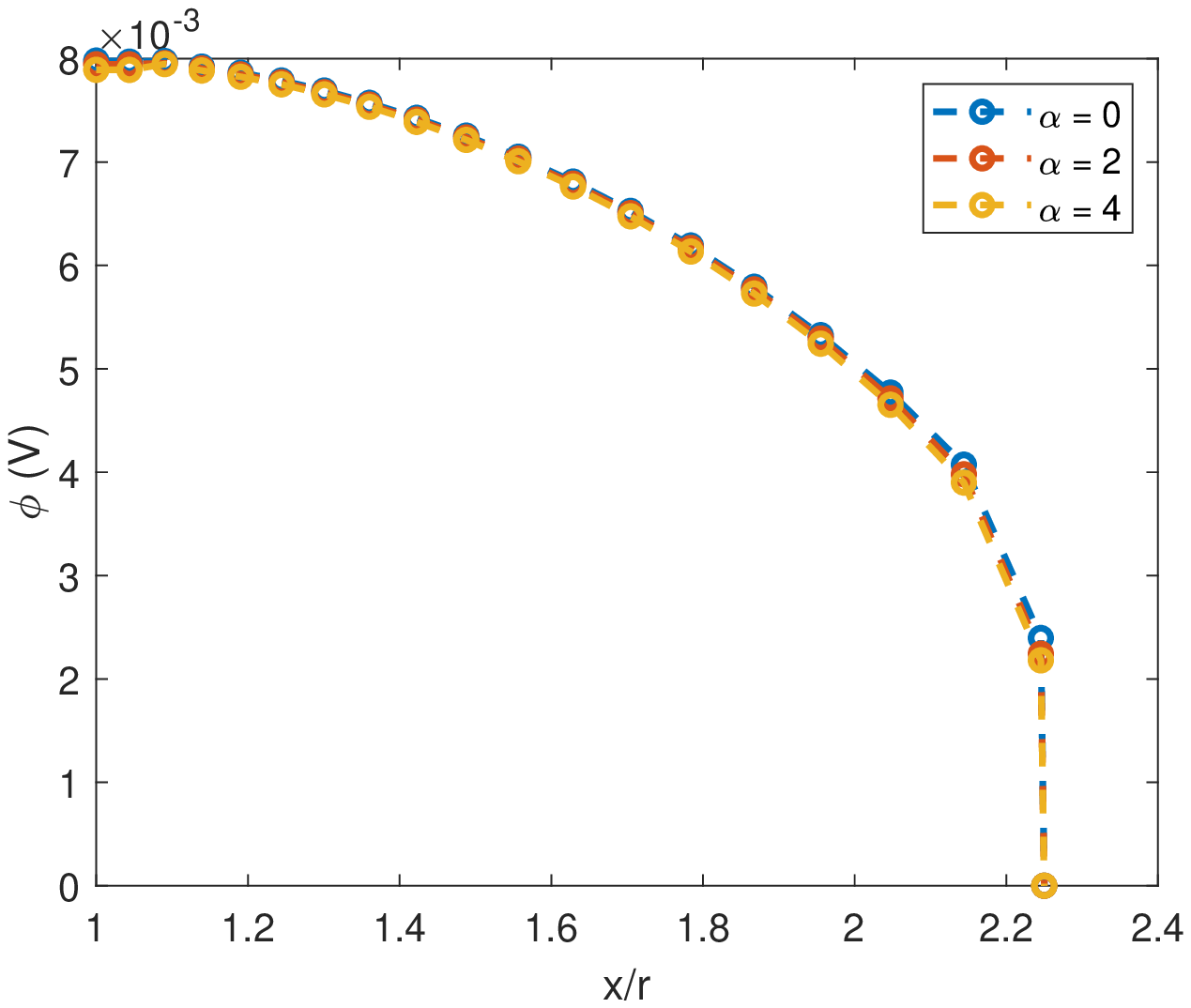}}  
        \subfigure[D3]{ 
    \label{fig:Ex23_02_01} 
    \includegraphics[width=0.48\textwidth]{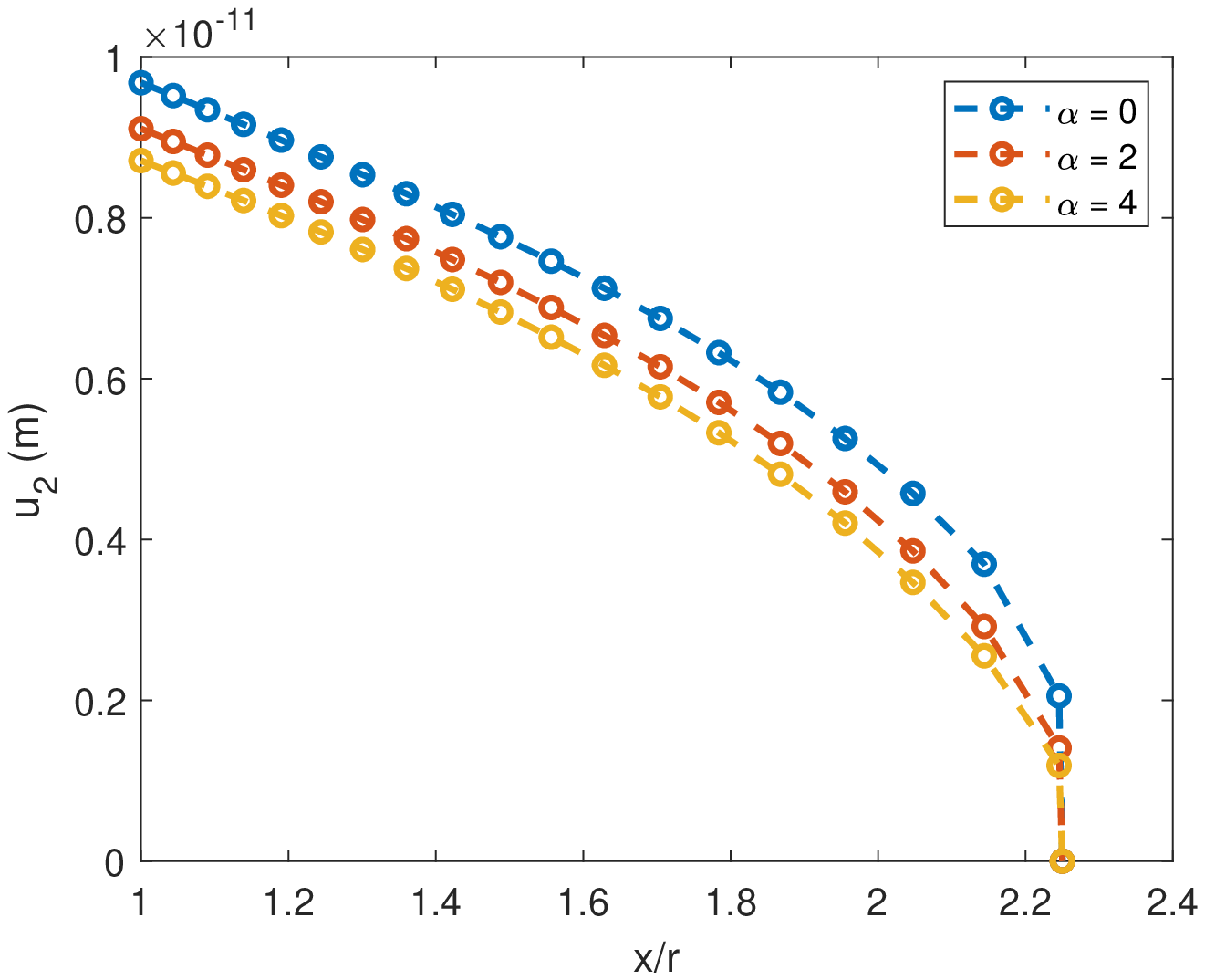}}  
    \subfigure[D3]{ 
    \label{fig:Ex23_02_02} 
    \includegraphics[width=0.48\textwidth]{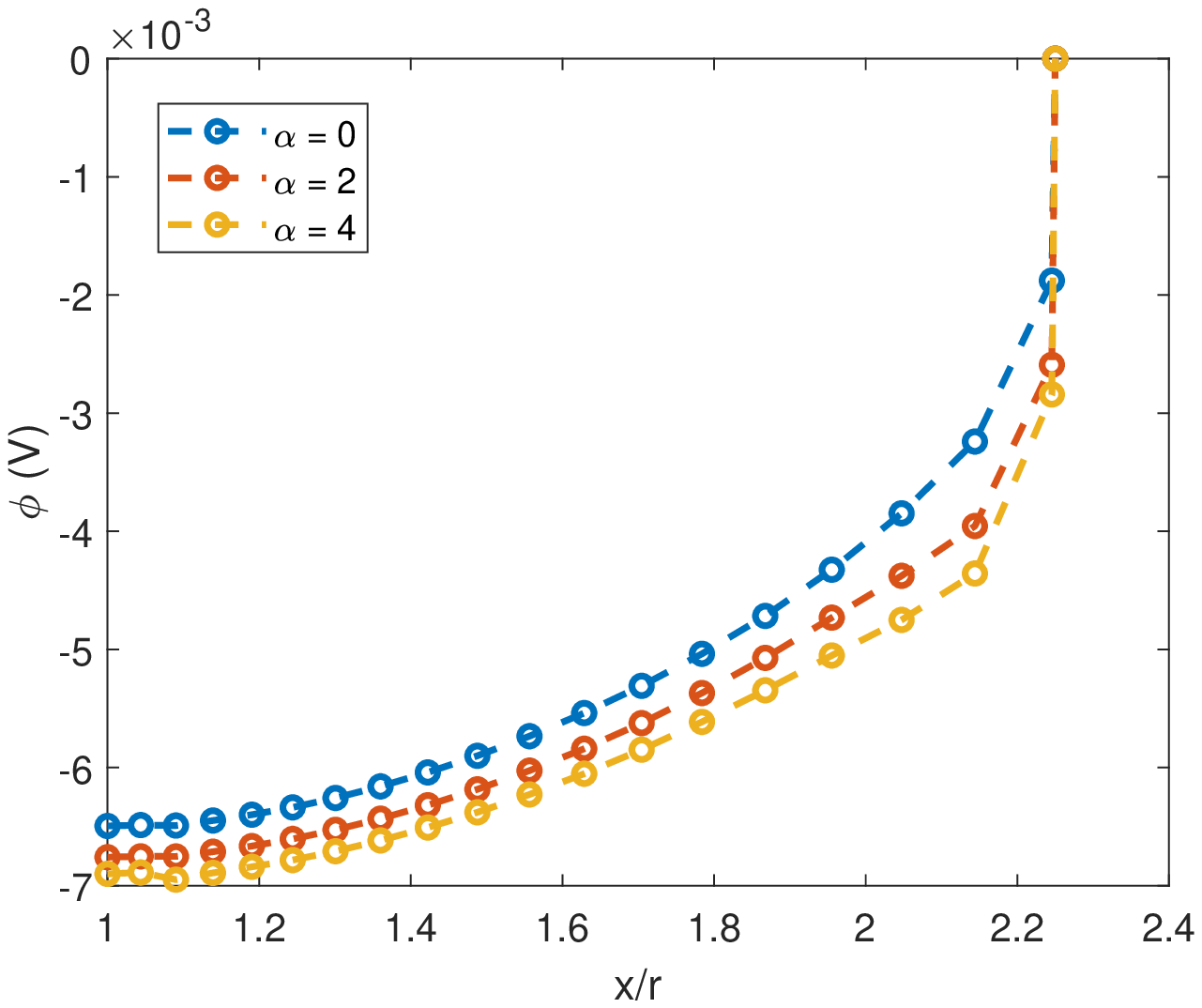}}   
  \caption{Distribution of displacement $u_2$ and electric potential $\phi$ on the upper crack-face. (a) (b) $D1 = -5 \times 10^{-4}~\mathrm{C/m^2}$. (c) (d) $D3 = 1 \times 10^{-3}~\mathrm{C/m^2}$} 
  \label{fig:Ex23_Solu} 
\end{figure}

The second example is shown in Fig.~\ref{fig:Ex23_Schem_01}. The plate with a circular hole at the center and two symmetric cracks on the boundary of the hole is subjected to the same mechanical and electrical loadings as the previous example, i.e., $\widetilde{Q} = 1.17~\mathrm{MPa}$ and $\widetilde{\omega} = D1, D2$ or $D3$. The hole surface is traction and charge free. And the cracks are electrically impermeable. The geometric parameters: $L = 1~\mathrm{\mu m}$, $r = 0.1~\mathrm{\mu m}$, $a = 0.125~\mathrm{\mu m}$. 3456 points are distributed in the domain and quadrilateral partition is applied (see Fig.~\ref{fig:Ex21_Schem_02}). The FPM computational parameters: $c_0 = \sqrt{20}$, $\eta_{21} = 2.0 E$, $\eta_{22} =50 E$, $\eta_{23} = 2.0 \Lambda_{33}$, $\eta_{24} = 0$. Fig.~\ref{fig:Ex23_Solu} exhibits the distribution of displacement $u_2$ and electric potential $\phi$ on the upper surface of the right crack under different electrical loadings and $\alpha$. Similar to the previous example, the presence of the strain and electric field gradients lead to a decrease of the crack-opening displacement. Yet their influence is limited. In this example, since the poling direction of the material is perpendicular to the crack, the behaviors of the two cracks are completely symmetric.

\subsection{A rectangular plate with an edge crack}

\begin{figure}[htbp] 
  \centering 
    \subfigure[]{ 
    \label{fig:Ex24_Schem_01} 
    \includegraphics[width=0.48\textwidth]{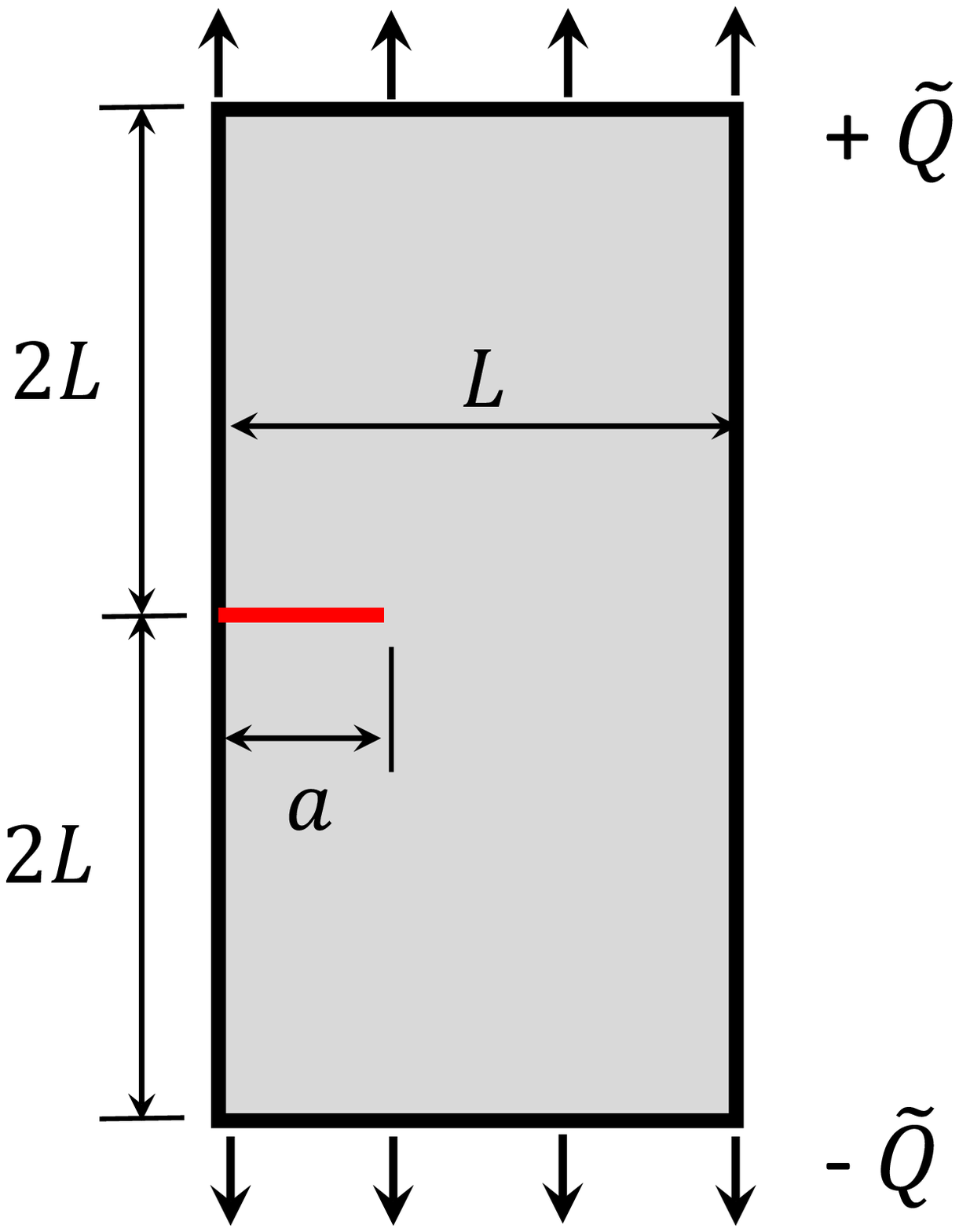}}  
    \subfigure[]{ 
    \label{fig:Ex24_Schem_02} 
    \includegraphics[width=0.48\textwidth]{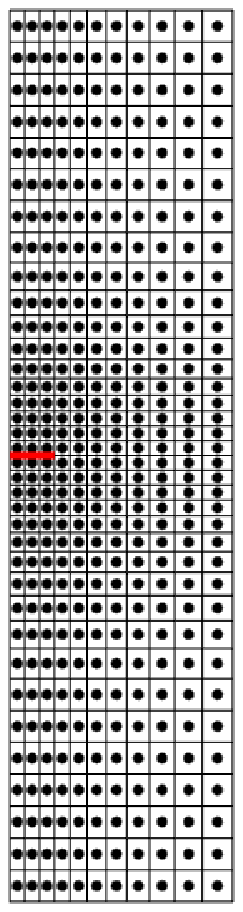}}  
  \caption{(a) A rectangular plate with an edge crack.  (b) The point distribution and domain partition in the FPM.} 
  \label{fig:Ex24_Schem} 
\end{figure}

At last, we consider a long rectangular plate with an edge crack. As shown in Fig.~\ref{fig:Ex24_Schem_01}, the plate is subjected to mechanical tensile stress $\widetilde{Q} = 1.17~\mathrm{MPa}$. The electrical loading is absent. The geometric parameters: $L = 1~\mathrm{\mu m}$, $a = 0.1~\mathrm{\mu m}$. In this example, two kinds of crack-face boundary conditions are considered (see Part. I for details). For the electrically impermeable cracks, the crack face is assumed to be surface charge free. Whereas for the electrically permeable cracks, the electric potential is continuous between the upper and lower crack faces.

\begin{figure}[htbp] 
  \centering 
    \subfigure[Electrically impermeable]{ 
    \label{fig:Ex24_01_01} 
    \includegraphics[width=0.48\textwidth]{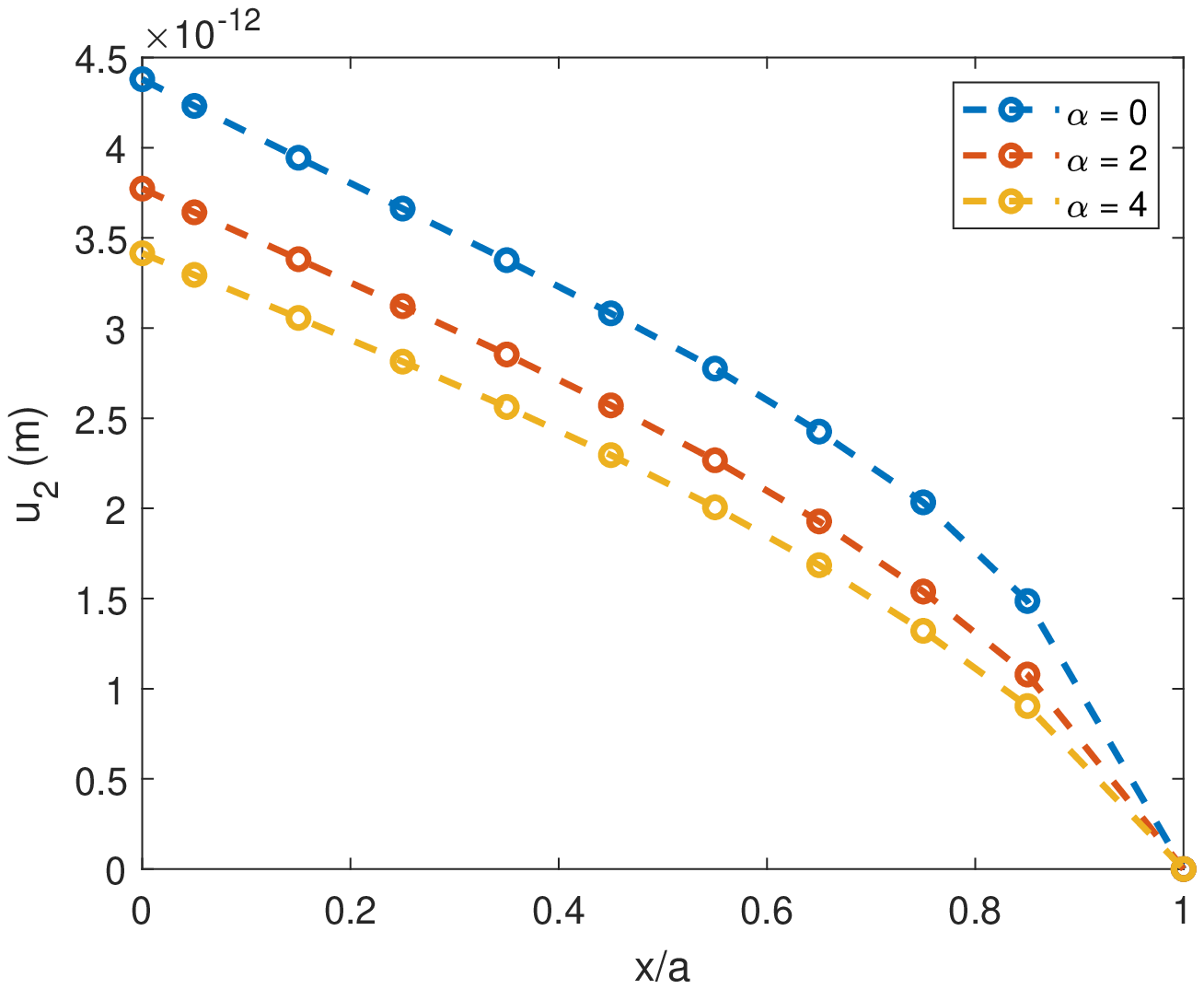}}  
    \subfigure[Electrically impermeable]{ 
    \label{fig:Ex24_01_02} 
    \includegraphics[width=0.48\textwidth]{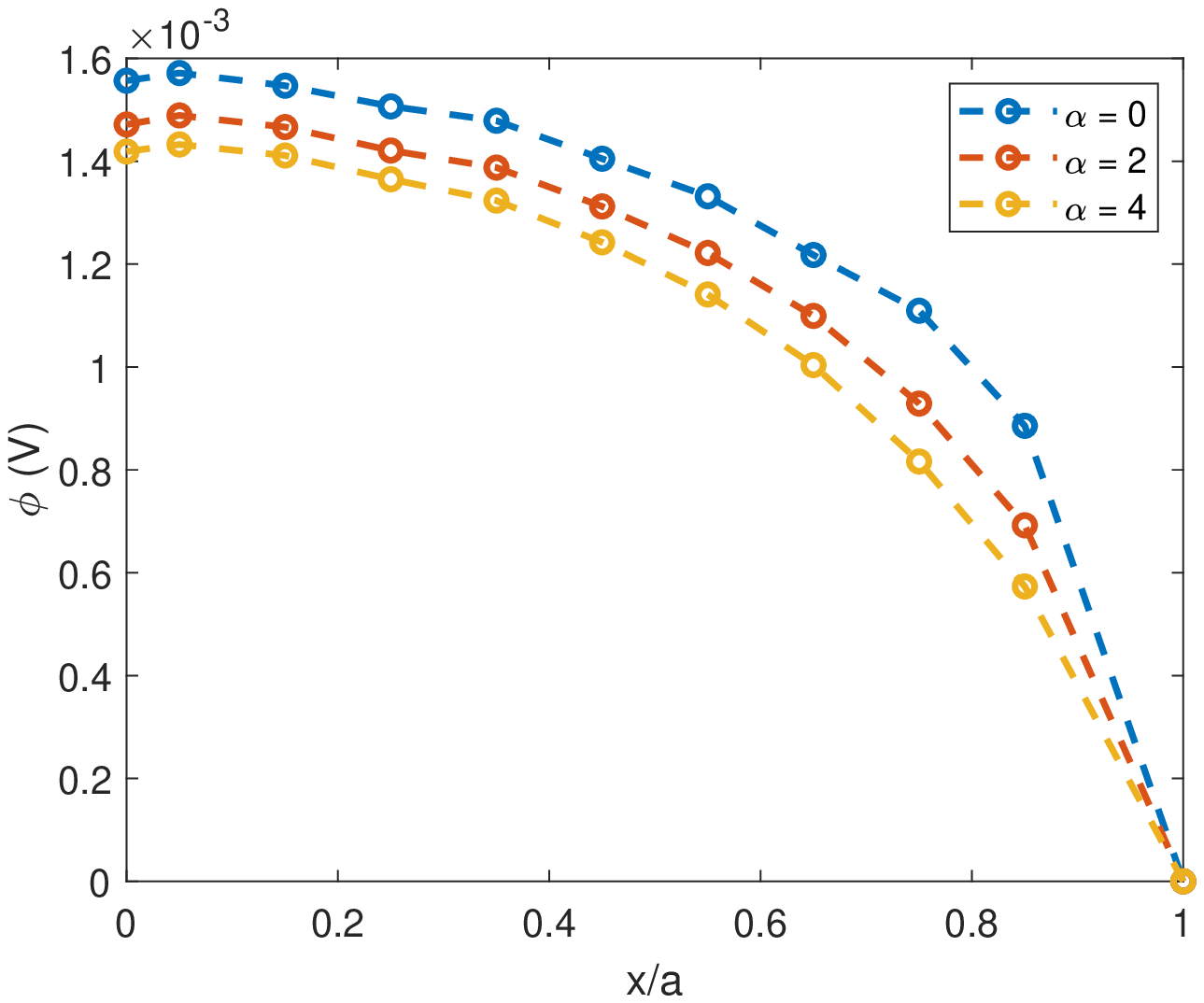}}  
        \subfigure[Electrically permeable]{ 
    \label{fig:Ex24_02_01} 
    \includegraphics[width=0.48\textwidth]{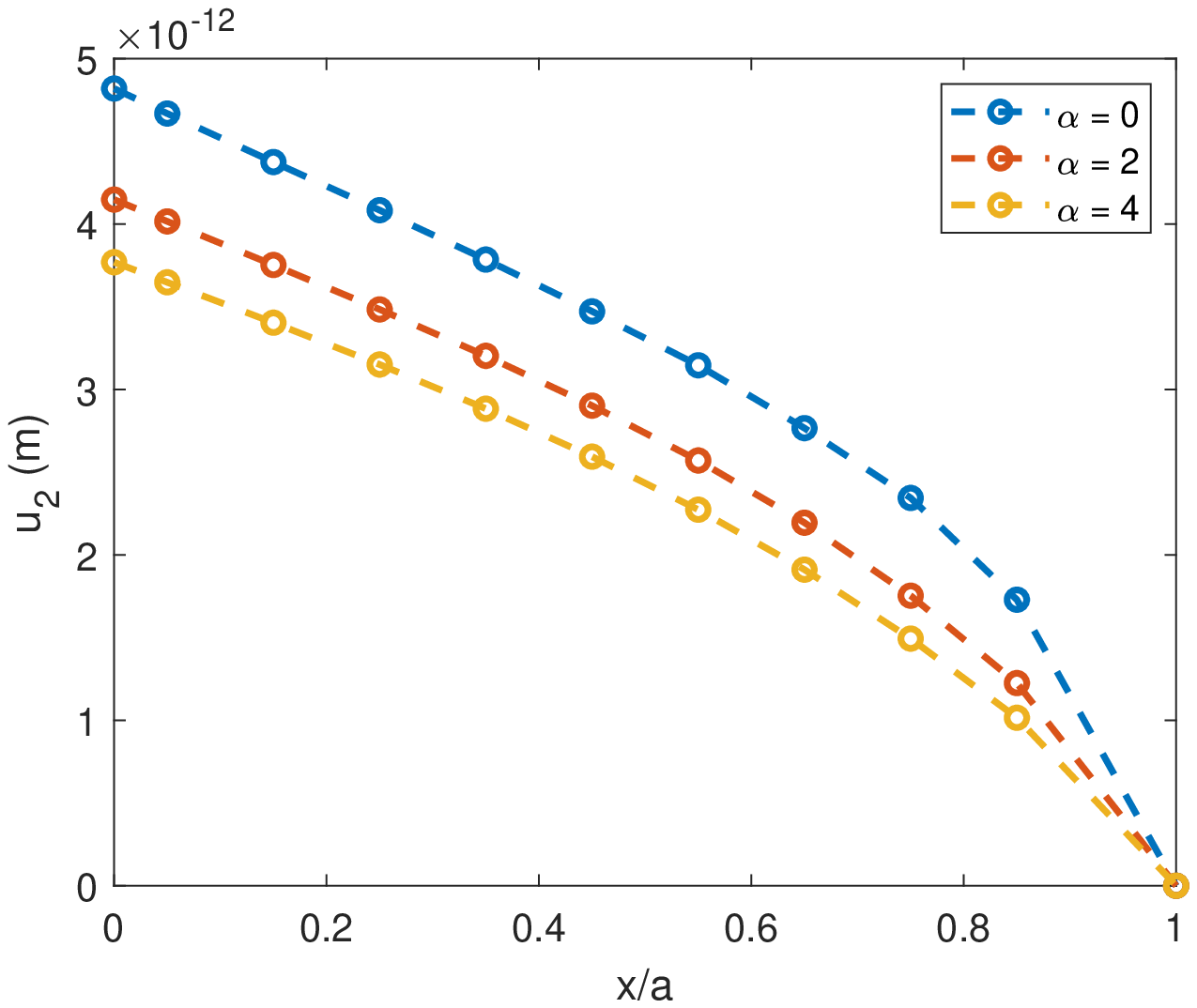}}  
    \subfigure[Electrically permeable]{ 
    \label{fig:Ex24_02_02} 
    \includegraphics[width=0.48\textwidth]{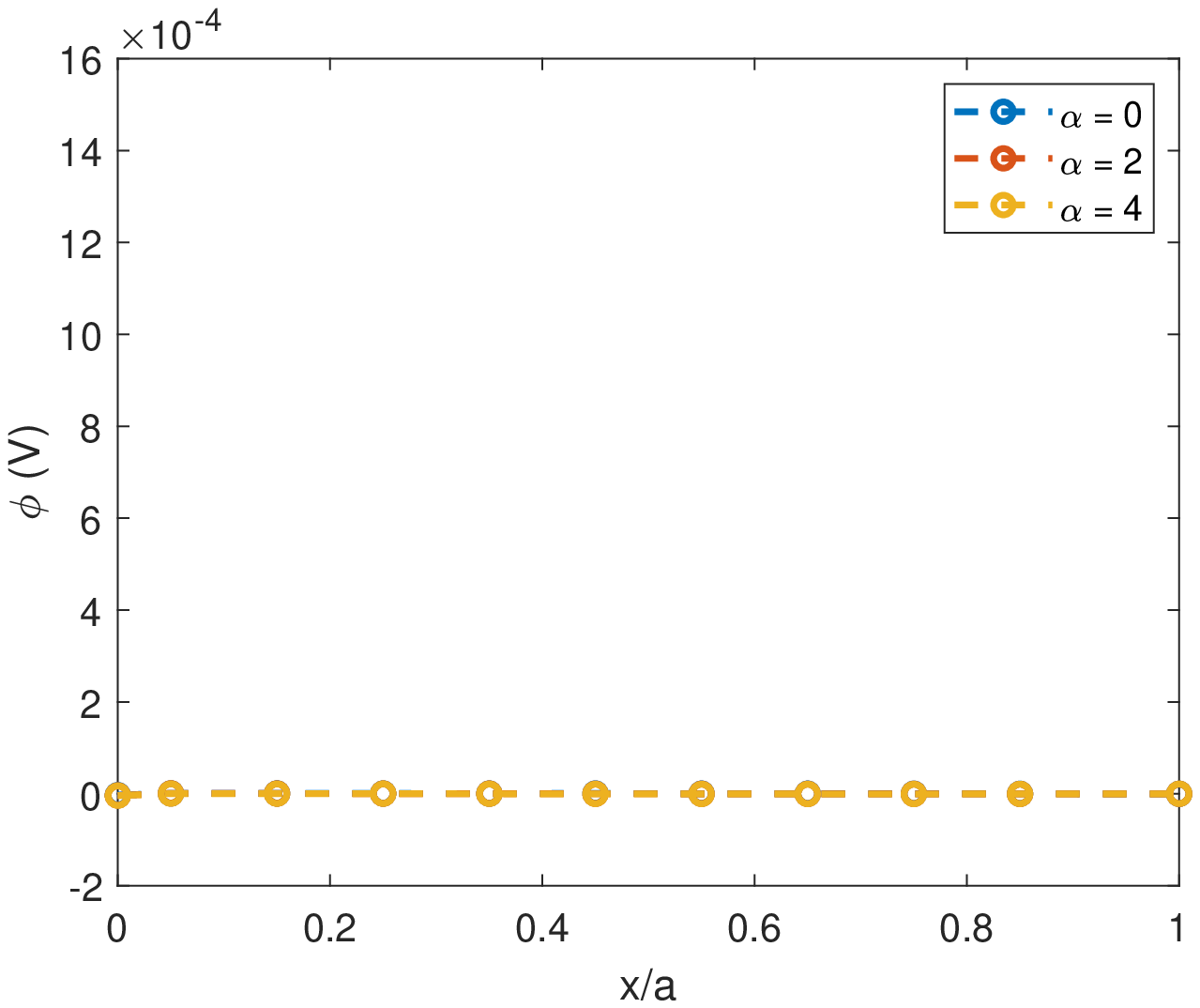}}   
  \caption{Distribution of displacement $u_2$ and electric potential $\phi$ on the upper crack-face. (a) (b) Electrically impermeable. (c) (d) Electrically permeable.} 
  \label{fig:Ex24_Solu} 
\end{figure}

The primal FPM is applied with 2052 points and quadrilateral partition. The point distribution is shown schematically in Fig.~\ref{fig:Ex24_Schem_02}. The computational parameters are the same as the previous example. Fig.~\ref{fig:Ex24_Solu} shows a comparison of the displacement $u_2$ and electric potential $\phi$ on the upper crack-face under the electrically permeable and impermeable crack-face boundary conditions. As can be seen, an electrically permeable crack results in a larger crack-opening displacement. And the results are also influenced by the strain and electric field gradient effect when $\alpha \not= 0$. As the system is symmetric with respect to the crack line, and the poling direction of the material is perpendicular to the symmetric axis, the electric potentials on both the upper and lower crack faces are zero. This is also verified by the computed solution in Fig.~\ref{fig:Ex24_02_02}. All the computed solutions for the three examples shown in this section are in good agreement with the study of \citet{Sladek2018, Sladek2018b} using FEM and MLPG method.

\section{Crack initiation and propagation}  \label{sec:crack_2}

\subsection{Simulations of crack propagation paths}

\begin{figure}[htbp] 
  \centering 
    \subfigure[]{ 
    \label{fig:Ex32_Schem_01} 
    \includegraphics[width=0.32\textwidth]{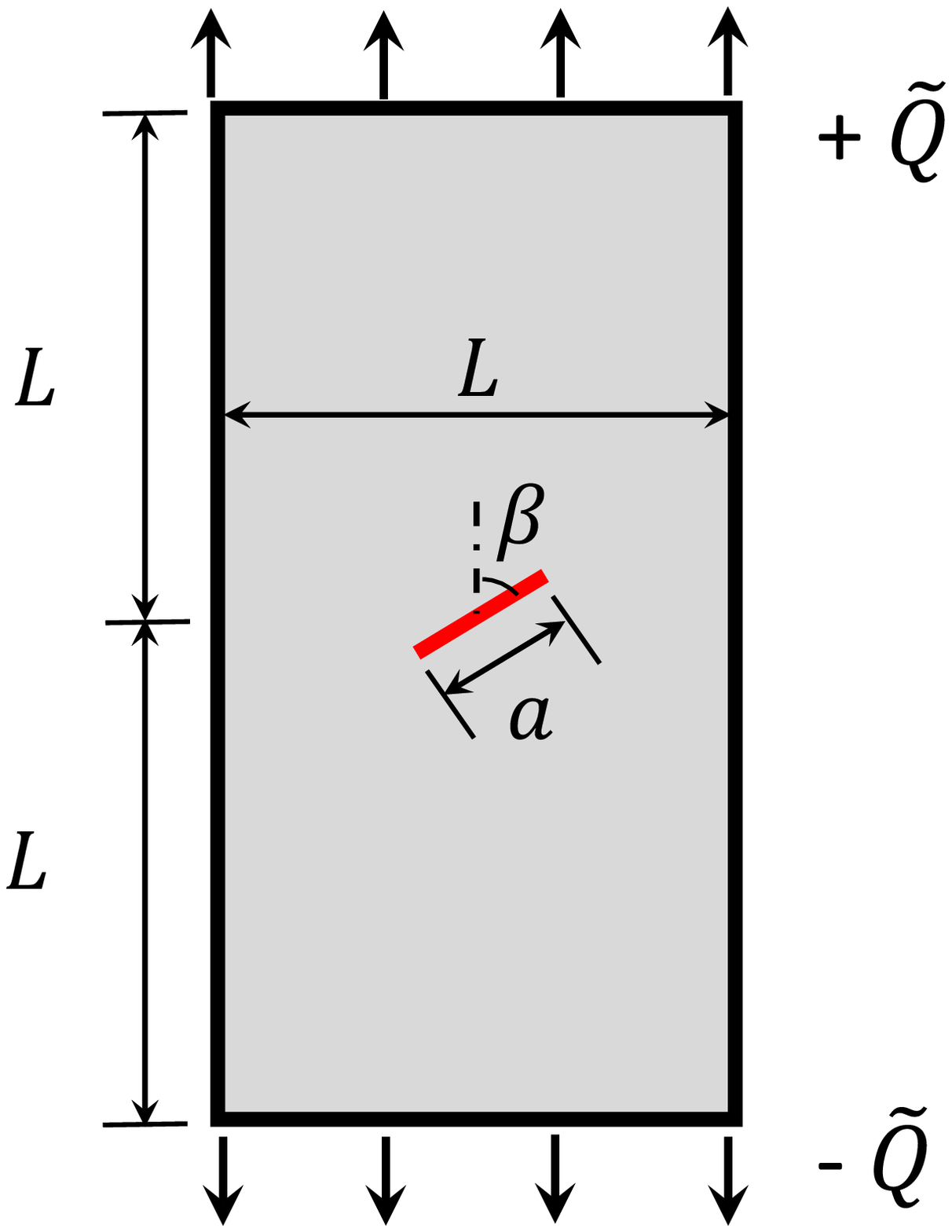}}  
    \subfigure[]{ 
    \label{fig:Ex32_Schem_02} 
    \includegraphics[width=0.32\textwidth]{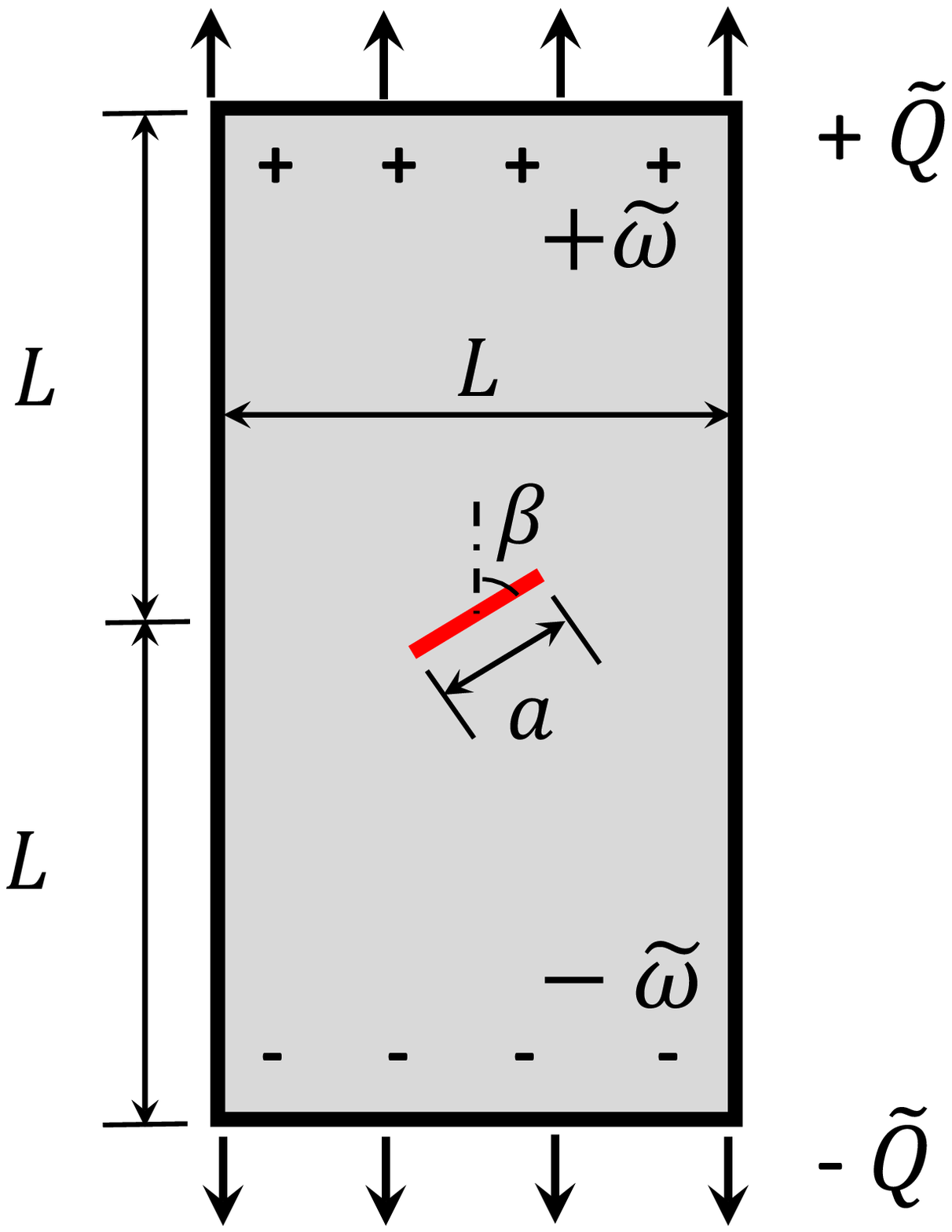}}  
        \subfigure[]{ 
    \label{fig:Ex32_Schem_03} 
    \includegraphics[width=0.32\textwidth]{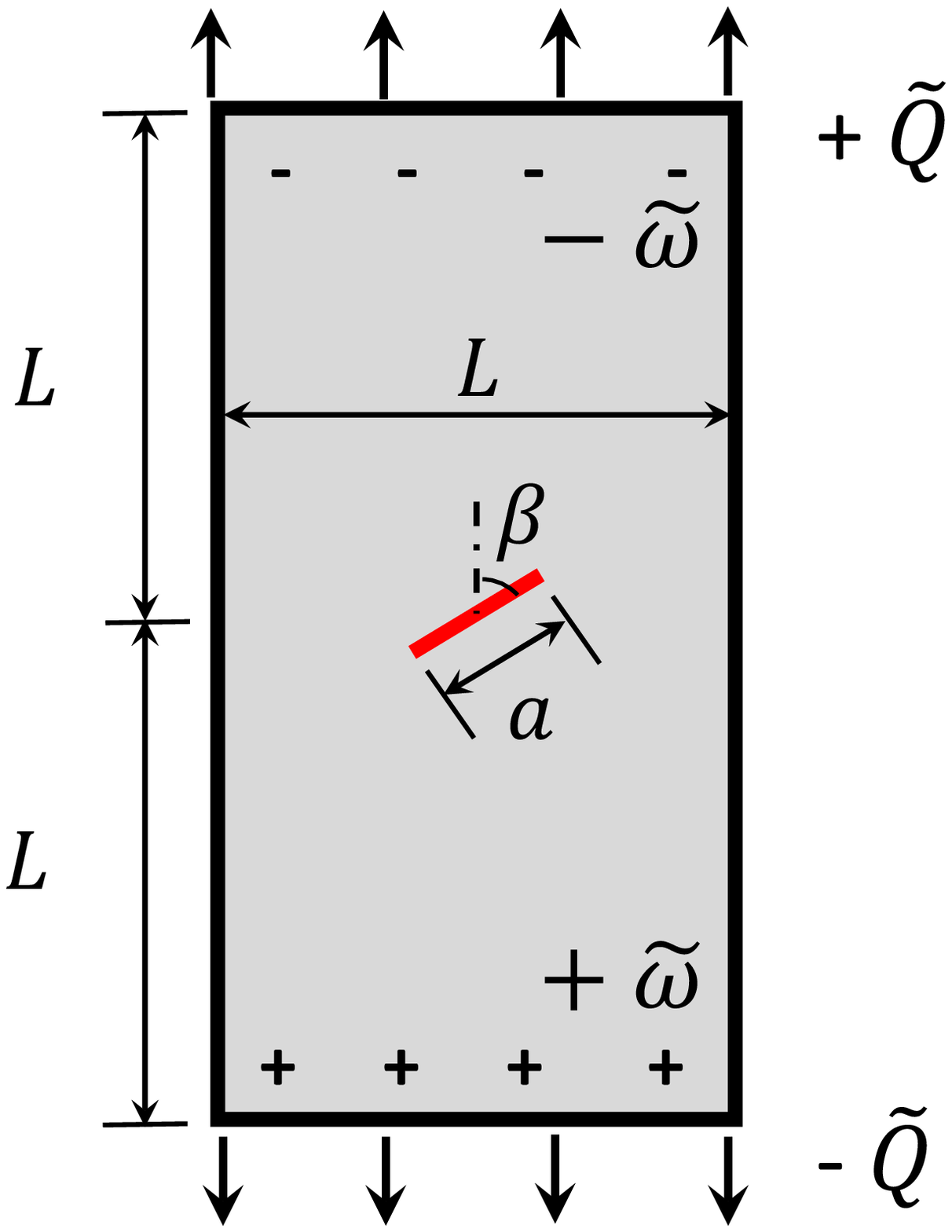}}  
  \caption{A plate with a pre-existing oblique crack at the center. (a) Pure mechanical loading. (b) Electric polarization parallel to $y$-axis. (c) Electric polarization antiparallel to $y$-axis.} 
  \label{fig:Ex32_Schem} 
\end{figure}

In this section, the FPM is employed in simulating crack propagation paths in dielectric solids under mechanical and electrical loadings. We consider a plate with a pre-existing oblique crack at the center. As shown in Fig.~\ref{fig:Ex32_Schem}, the plate is subjected to a tensile stress $\widetilde{Q}$ and/or an electrical surface charge $\widetilde{\omega}$. The example is initiated by Ex.~5, in which the flexoelectric effect leads to an asymmetric response at the tip of a hole, and therefore, may influence the crack propagation direction. We consider the same material as in Ex.~5, i.e., $E = 139~\mathrm{GPa}$, Poisson's ratio $\nu = 0.3$, $l = 3.33 \mathrm{\mu m}$, $\overline{\mu}_{11} = -6.3 \times 10^{-5}~\mathrm{C/m}$, $\overline{\mu}_{12} = 5.2 \times 10^{-6}~\mathrm{C/m}$, $\overline{\mu}_{44} = -3.4 \times 10^{-5}~\mathrm{C/m}$, $\Lambda_{11} = \Lambda_{33} = 4.9 \times 10^{-9}~\mathrm{F/m}$, $e_{31} = e_{33} = e_{15} = 0$. The plate width $L = 110~\mathrm{\mu m}$, the crack length $a = 20~\mathrm{\mu m}$ and the crack angle $\beta = 60^{\circ}$. The external loadings $\widetilde{Q} = 695~\mathrm{MPa}$, $\widetilde{\omega} = 0.5~\mathrm{C/m^2}$.

\begin{figure}[htbp] 
  \centering 
    \subfigure[]{ 
    \label{fig:Ex32_01_01} 
    \includegraphics[width=0.32\textwidth]{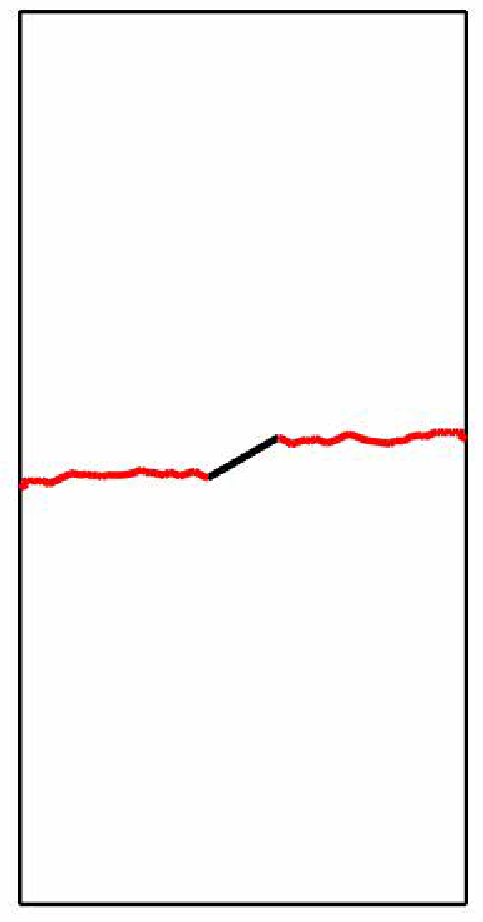}}  
    \subfigure[]{ 
    \label{fig:Ex32_01_02} 
    \includegraphics[width=0.32\textwidth]{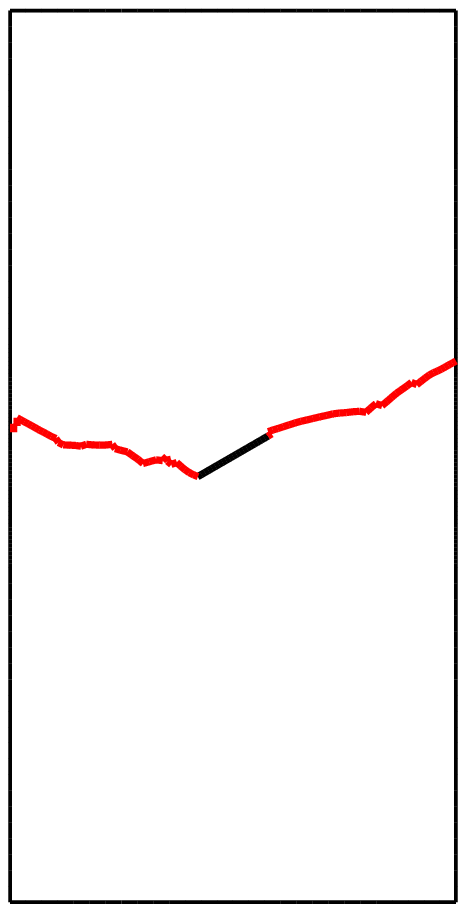}}  
    \subfigure[]{ 
    \label{fig:Ex32_01_03} 
    \includegraphics[width=0.32\textwidth]{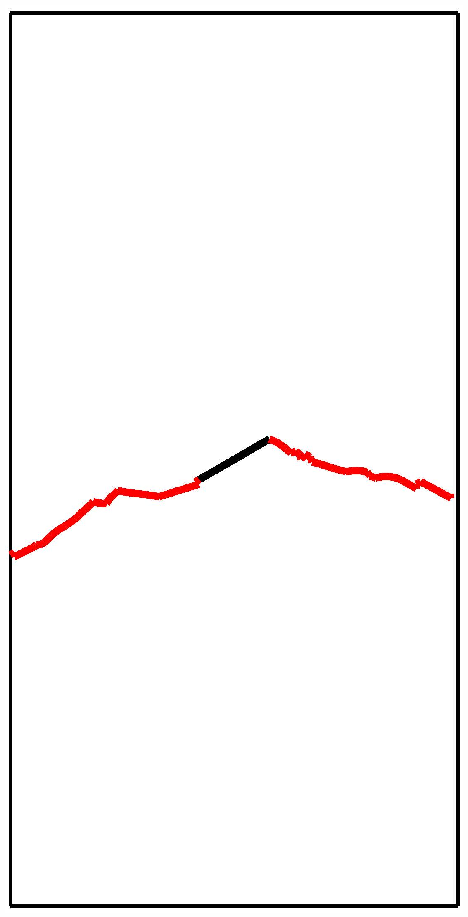}}  
    \subfigure[]{ 
    \label{fig:Ex32_01_04} 
    \includegraphics[width=0.32\textwidth]{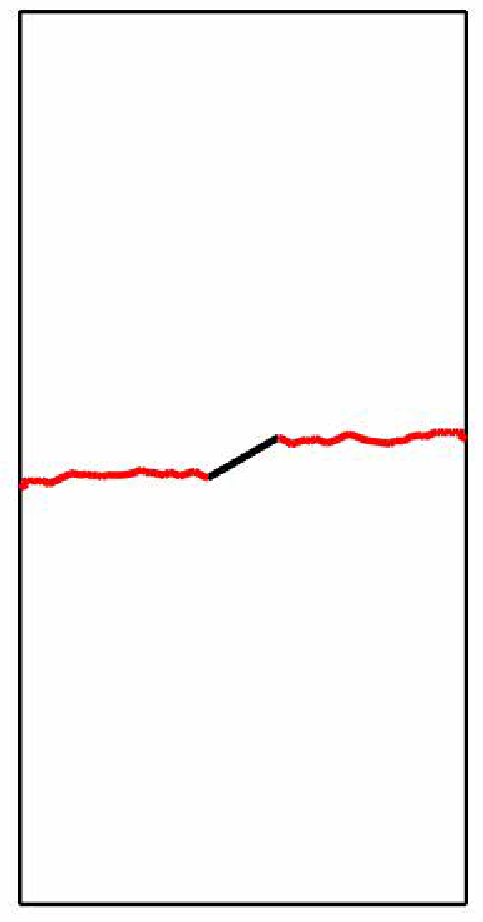}}  
    \subfigure[]{ 
    \label{fig:Ex32_01_05} 
    \includegraphics[width=0.32\textwidth]{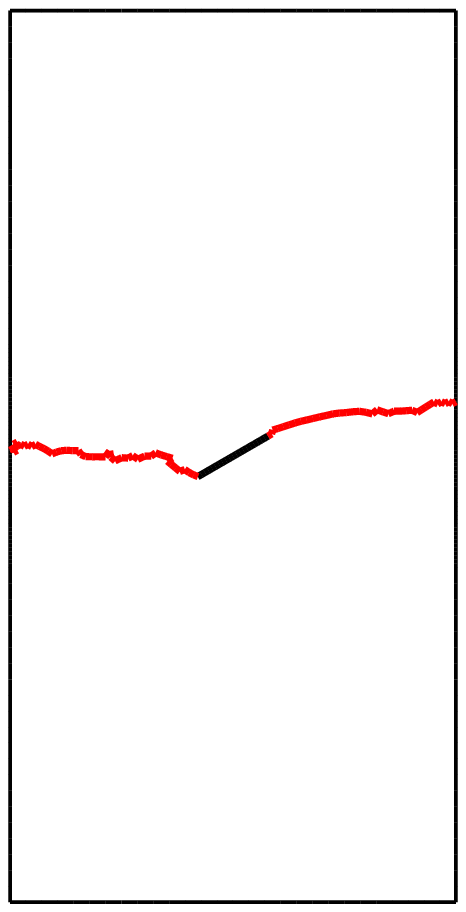}}  
    \subfigure[]{ 
    \label{fig:Ex32_01_06} 
    \includegraphics[width=0.32\textwidth]{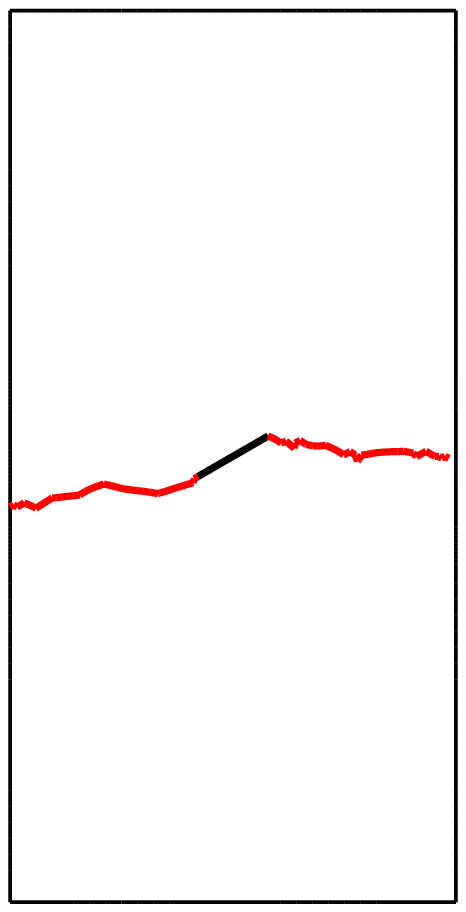}}  
  \caption{Crack propagation paths simulated by the FPM. (a) Pure mechanical loading ($l = 3.33~\mathrm{\mu m}$). (b) Electric polarization parallel to $y$-axis ($l = 3.33~\mathrm{\mu m}$). (c) Electric polarization antiparallel to $y$-axis ($l = 3.33~\mathrm{\mu m}$). (c) Pure mechanical loading ($l = 0$). (d) Electric polarization parallel to $y$-axis ($l = 0$). (e) Electric polarization antiparallel to $y$-axis ($l = 0$).} 
  \label{fig:Ex32_01} 
\end{figure}

A Maximum Hoop Stress criterion proposed by \citet{Erdogan1963} is applied to predict the crack propagation paths. In each load step, for all the internal boundaries connected with the current crack tip, the one with the maximum normal stress will be cracked. The computational parameters: $c_0 = \sqrt{20}$, $\eta_{21} = 2.0 E$, $\eta_{22} =50 E$, $\eta_{23} = 0 $. Fig.~\ref{fig:Ex32_01} presents the crack propagation paths simulated with the primal FPM under multiple electrical loadings with or without the strain gradient effect. When the plate is subjected to pure tensile loading, as shown in Fig.~\ref{fig:Ex32_01_04}, it is a simple linear elastic problem and shows the same crack propagation paths as its scaled-up model in a $110~\mathrm{mm} \times 220~\mathrm{mm}$ plate. The results are in good agreement with FPM simulation results based on linear elastic mechanics \cite{Yang2019} and experimental data \cite{Mageed1991}. When external electrical loadings are applied, as can be seen in Fig.~\ref{fig:Ex32_01_02} and \ref{fig:Ex32_01_03}, due to the flexoelectric phenomena, the  growth of crack propagation paths in both directions deflect significantly towards the positive electrode. In Fig.~\ref{fig:Ex32_01_04} and \ref{fig:Ex32_01_06}, the strain gradient effect disappears ($l = 0$). When comparing with the previous solutions, the absence of the strain gradient effect results in a decline of the deflection effect.

\begin{figure}[htbp] 
  \centering 
    \subfigure[]{ 
    \label{fig:Ex32_Part_01} 
    \includegraphics[width=0.32\textwidth]{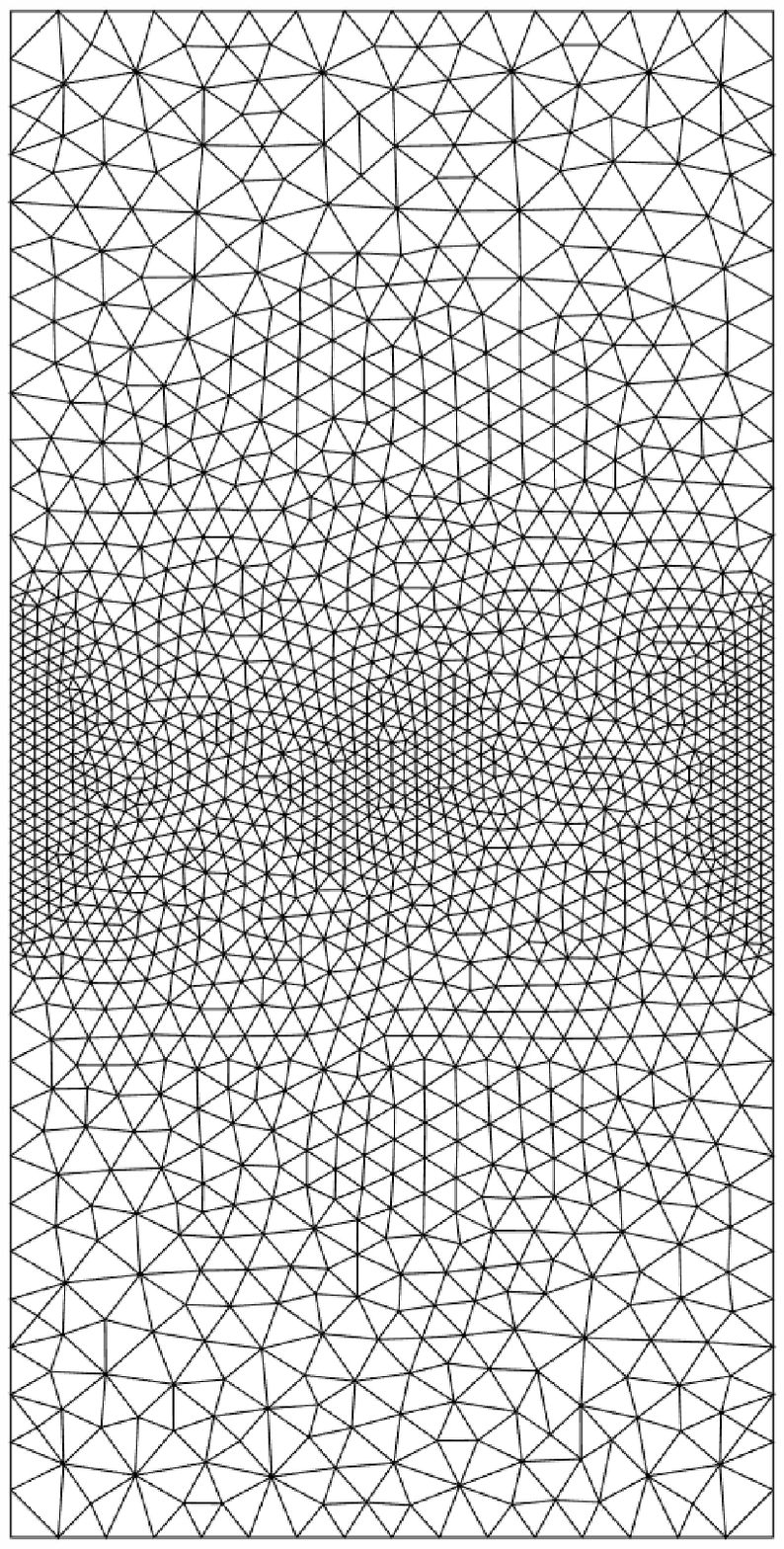}}  
    \subfigure[]{ 
    \label{fig:Ex32_Part_02} 
    \includegraphics[width=0.32\textwidth]{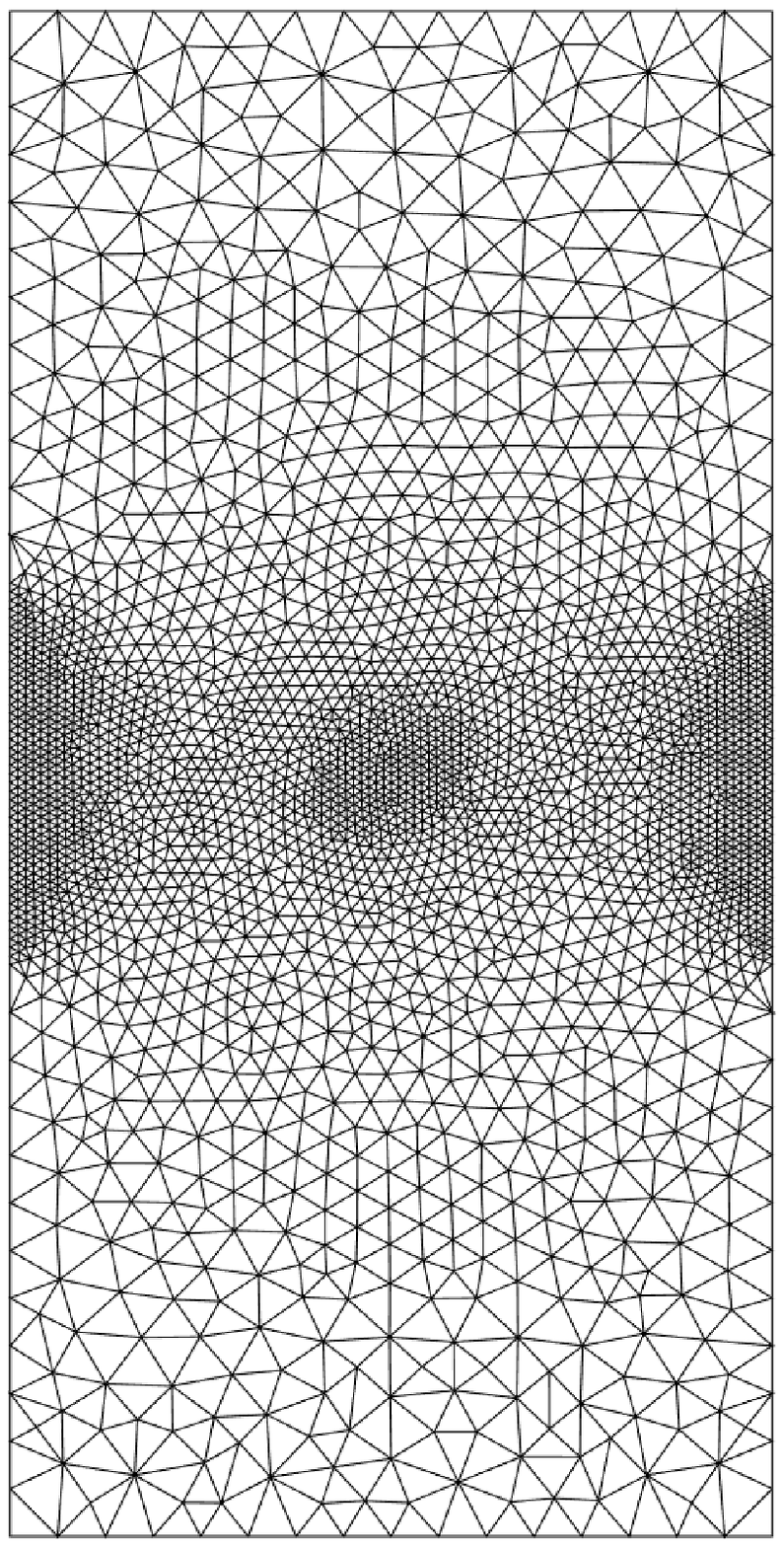}}  
        \subfigure[]{ 
    \label{fig:Ex32_Part_03} 
    \includegraphics[width=0.32\textwidth]{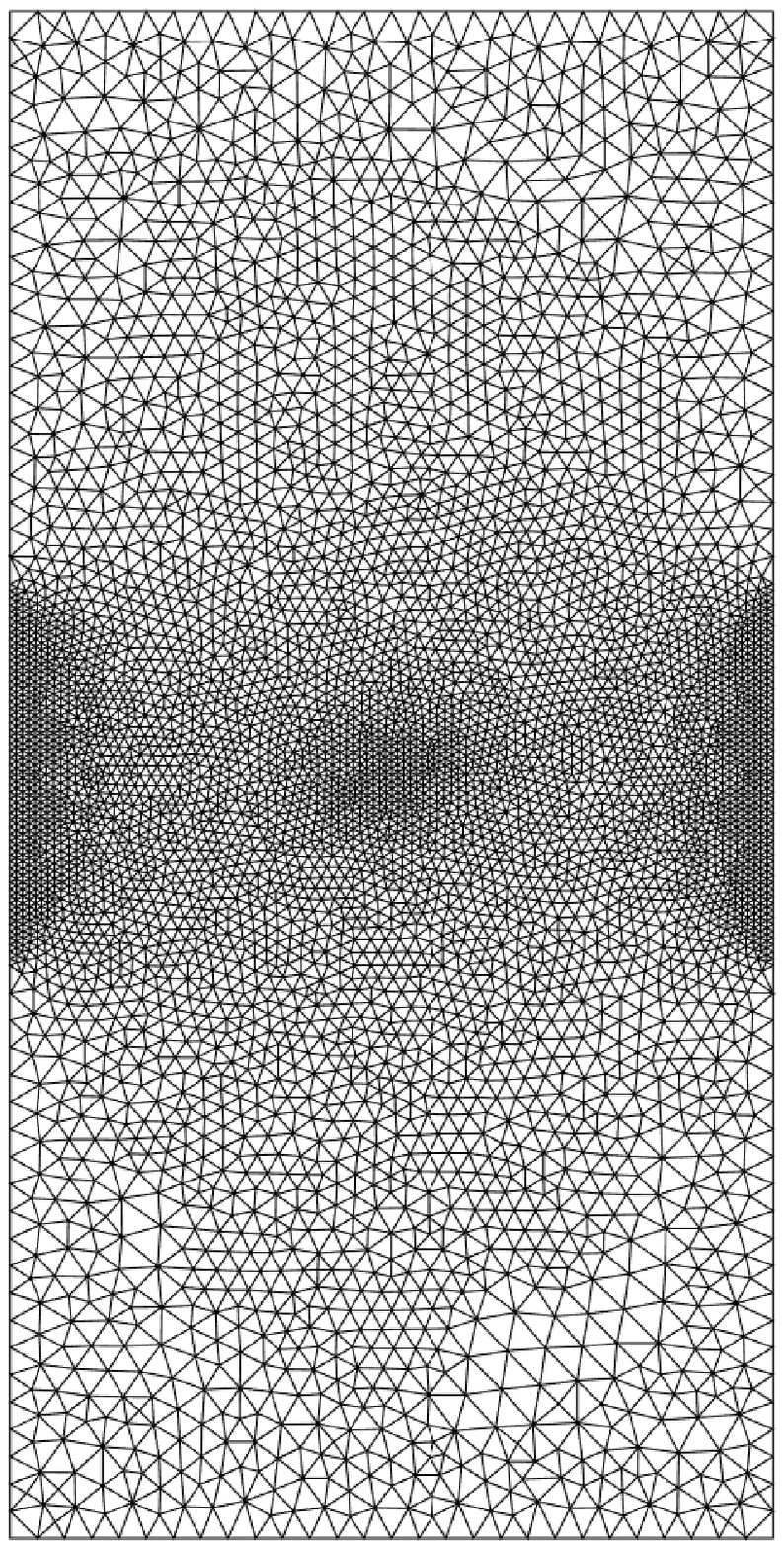}}  
  \caption{The domain partition used in the FPM with (a) 3178 points (b) 5346 points (c) 10396 points.} 
  \label{fig:Ex32_Part} 
\end{figure}

\begin{figure}[htbp] 
  \centering 
    \subfigure[]{ 
    \label{fig:Ex32_02_01} 
    \includegraphics[width=0.32\textwidth]{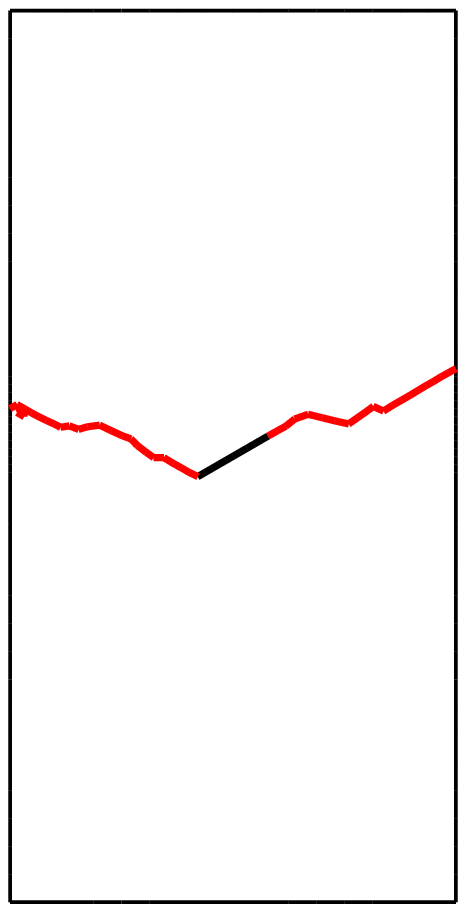}}  
    \subfigure[]{ 
    \label{fig:Ex32_02_02} 
    \includegraphics[width=0.32\textwidth]{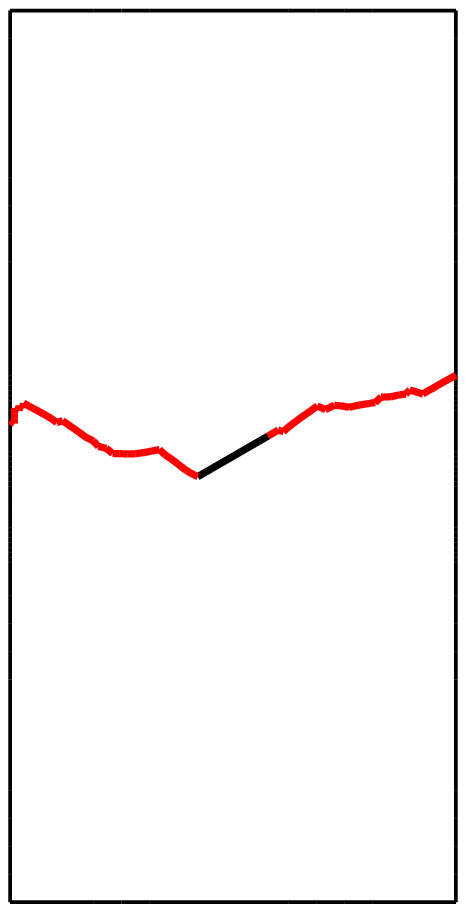}}  
        \subfigure[]{ 
    \label{fig:Ex32_02_03} 
    \includegraphics[width=0.32\textwidth]{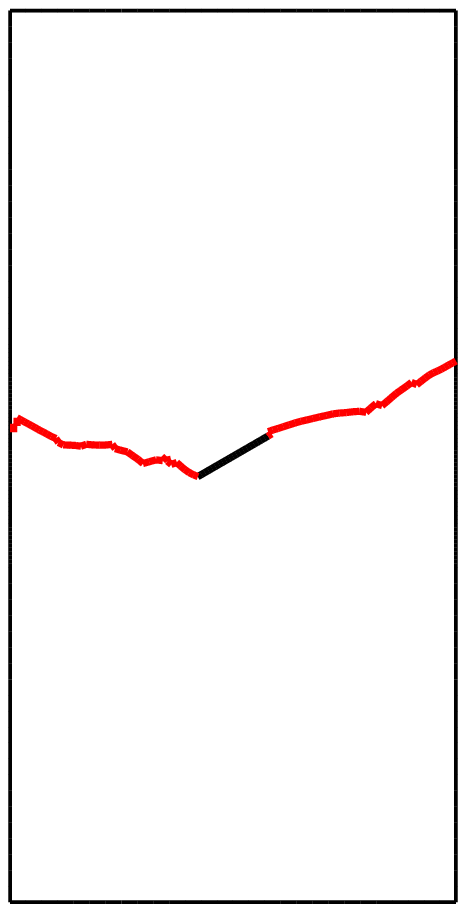}}  
  \caption{The crack propagation paths simulated by the FPM with (a) 3178 points (b) 5346 points (c) 10396 points. (Electric polarization parallel to $y$-axis, $l = 3.33~\mathrm{\mu m}$).} 
  \label{fig:Ex32_02} 
\end{figure}

Fig.~\ref{fig:Ex32_Part} shows three different domain partitions with 3178, 5346 and 10396 points used in the FPM analysis. The partitions are converted directly from ABAQUS meshing. The corresponding FPM Points are located at the centroid of each subdomain. The crack propagation paths simulated by FPM with these three partitions for Fig.~\ref{fig:Ex32_Schem_02} with $l = 3.33~\mathrm{\mu m}$ are presented in Fig.~\ref{fig:Ex32_02_01} -- \ref{fig:Ex32_02_03} respectively. As can be seen, when the number of points increases, the computed results approach a better estimate of the crack propagation path.

\subsection{Simulations of the crack initiation process}

Moreover, the electrical polarization can also result from the piezoelectric response of the material itself. Therefore, the coupling between the piezoelectric and flexoelectric effects may help, hinder or deflect the crack propagation in dielectric materials. This phenomenon is also observed in in experiments \cite{Cordero-Edwards2019}.

\begin{figure}[htbp] 
  \centering 
    \subfigure[]{ 
    \label{fig:Ex33_Schem_01} 
    \includegraphics[width=0.48\textwidth]{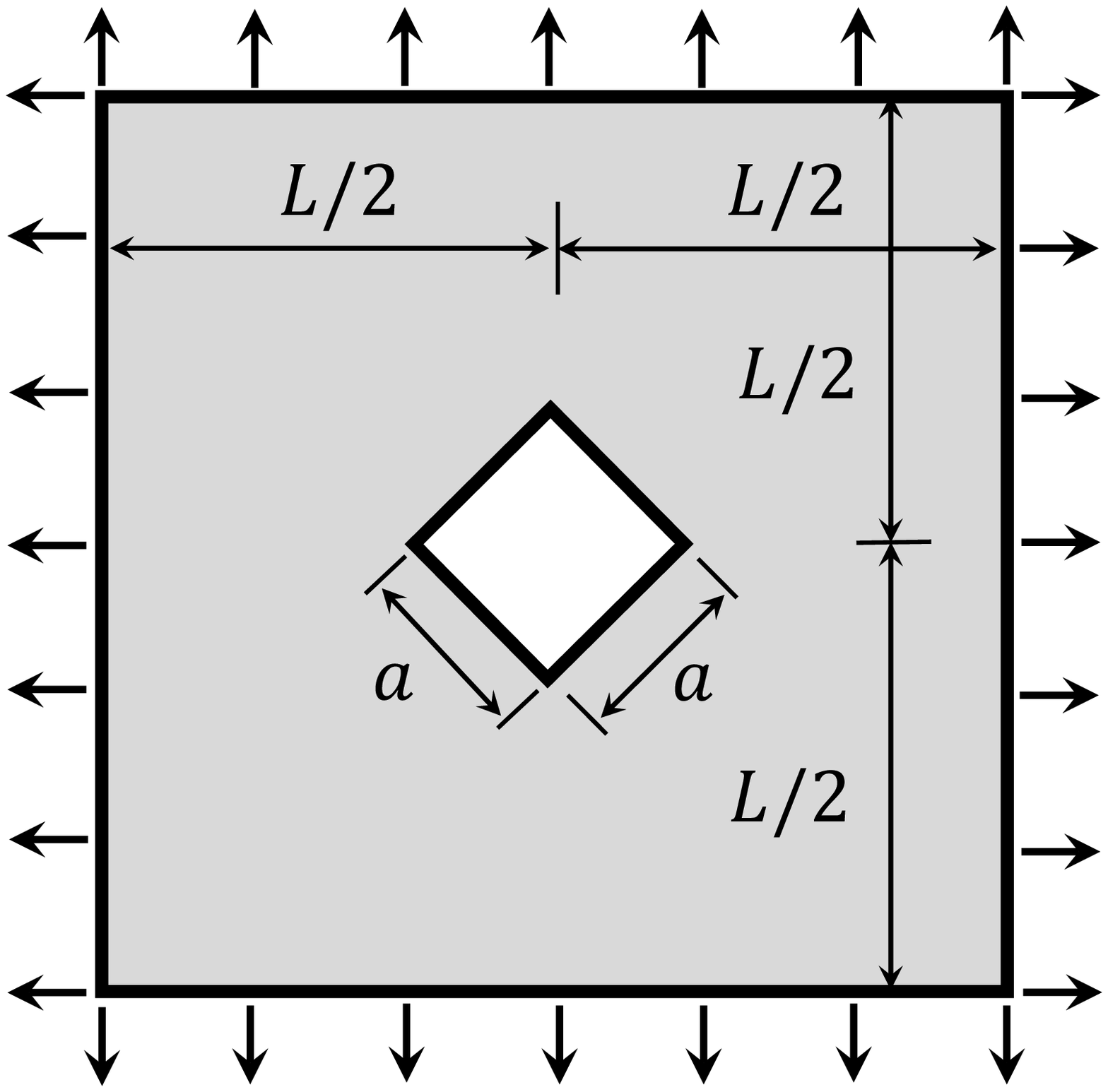}}  
    \subfigure[]{ 
    \label{fig:Ex33_Part_01} 
    \includegraphics[width=0.48\textwidth]{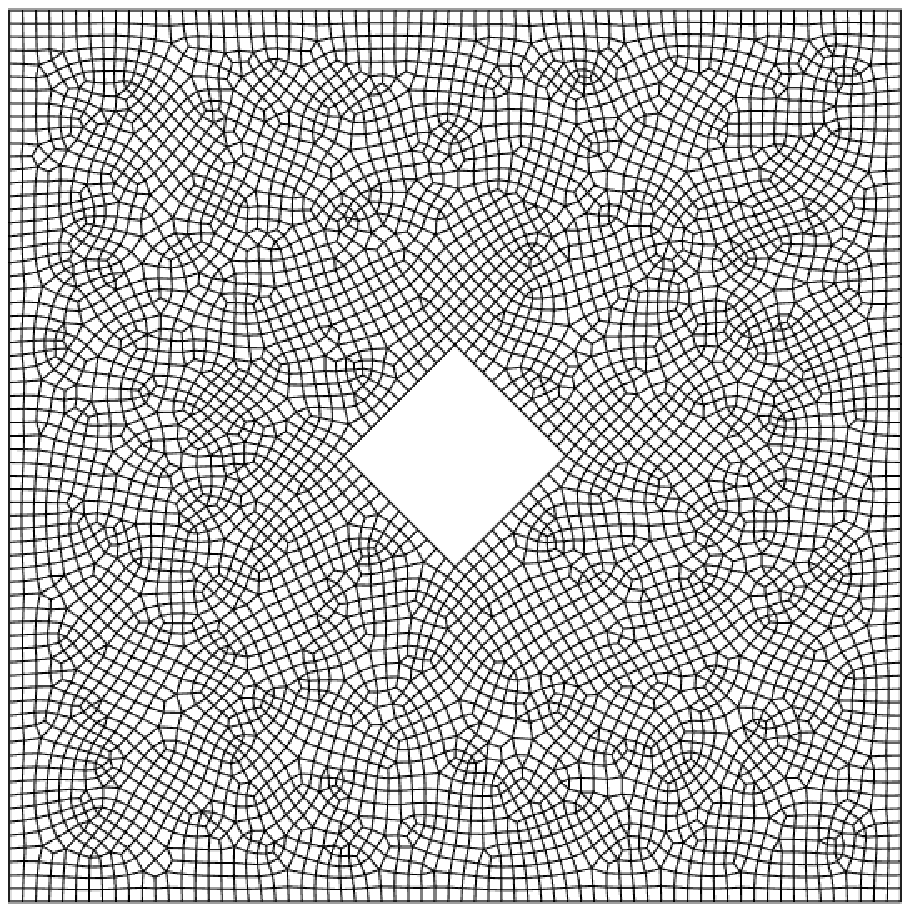}}  
  \caption{(a) A square plate with a central square hole subjected to biaxial loads. (a) The domain partition used in the FPM with 5278 points.} 
  \label{fig:Ex33_Schem} 
\end{figure}

In this section, we consider the process of crack initiation and development in a square plate with a central square hole with both piezoelectric and flexoelectric behaviors. As shown in Fig.~\ref{fig:Ex33_Schem_01}, the plate is under pure mechanical loadings. The biaxial tensile tractions in the two directions are identical. No external electrical loading is applied. The plate width $L = 40~\mathrm{\mu m}$, and the width of the square hole $a = \sqrt{2}/2~\mathrm{\mu m}$. The material properties: $E = 139~\mathrm{GPa}$, Poisson's ratio $\nu = 0.3$, $l = 0.1\mathrm{\mu m}$, $\overline{\mu}_{11} = \overline{\mu}_{12} = \overline{\mu}_{44} = 0.8 \times 10^{-5}~\mathrm{C/m}$, $\Lambda_{11} = \Lambda_{33} = 4.9 \times 10^{-9}~\mathrm{F/m}$, $e_{31} = 20~\mathrm{C/m^2}$, $e_{33} = e_{15} = 0$.

\begin{figure}[htbp] 
  \centering 
    \subfigure[]{ 
    \label{fig:Ex33_00_01} 
    \includegraphics[width=0.48\textwidth]{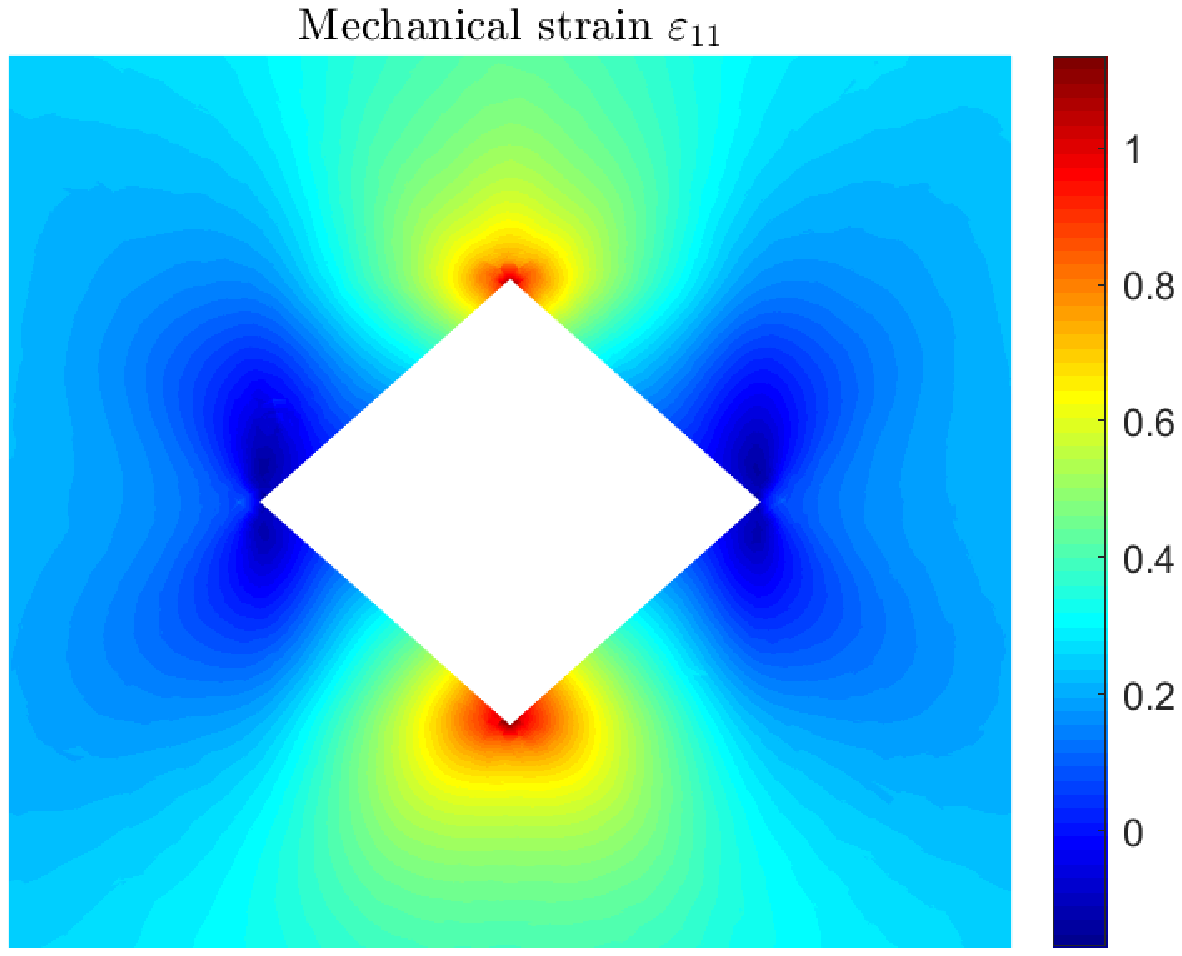}}  
    \subfigure[]{ 
    \label{fig:Ex33_00_02} 
    \includegraphics[width=0.48\textwidth]{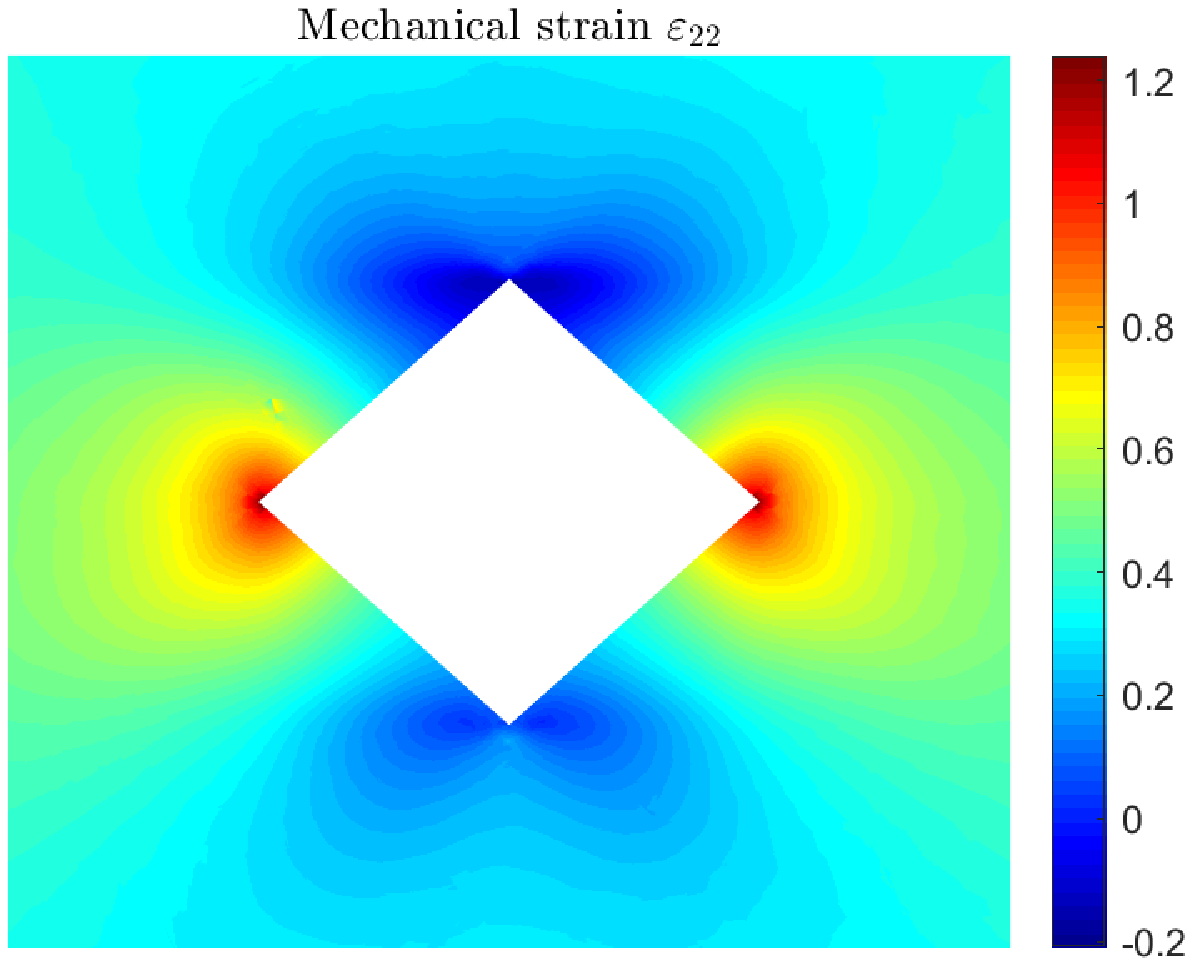}}  
  \caption{The computed solution of for a square plate with a central square hole. (a) Distribution of the mechanical strain $\varepsilon_{11}$. (b) Distribution of the mechanical strain $\varepsilon_{22}$.} 
  \label{fig:Ex33_00} 
\end{figure}

When the poling direction of the material is antiparallel to $y$-axis, the corresponding distribution of mechanical strain $\varepsilon_{11}$ and $\varepsilon_{22}$ are shown in Fig.~\ref{fig:Ex33_00}. Resulting from the coupling of piezoelectric and flexoelectric effects, the symmetry along $y$-axis breaks. Therefore, the crack propagation parallel to the poling direction of material can be helped by the coupling effect, whereas the crack propagation antiparallel to the poling direction will be hindered.

\begin{figure}[htbp] 
  \centering 
    \subfigure[]{ 
    \label{fig:Ex33_01_01} 
    \includegraphics[width=0.48\textwidth]{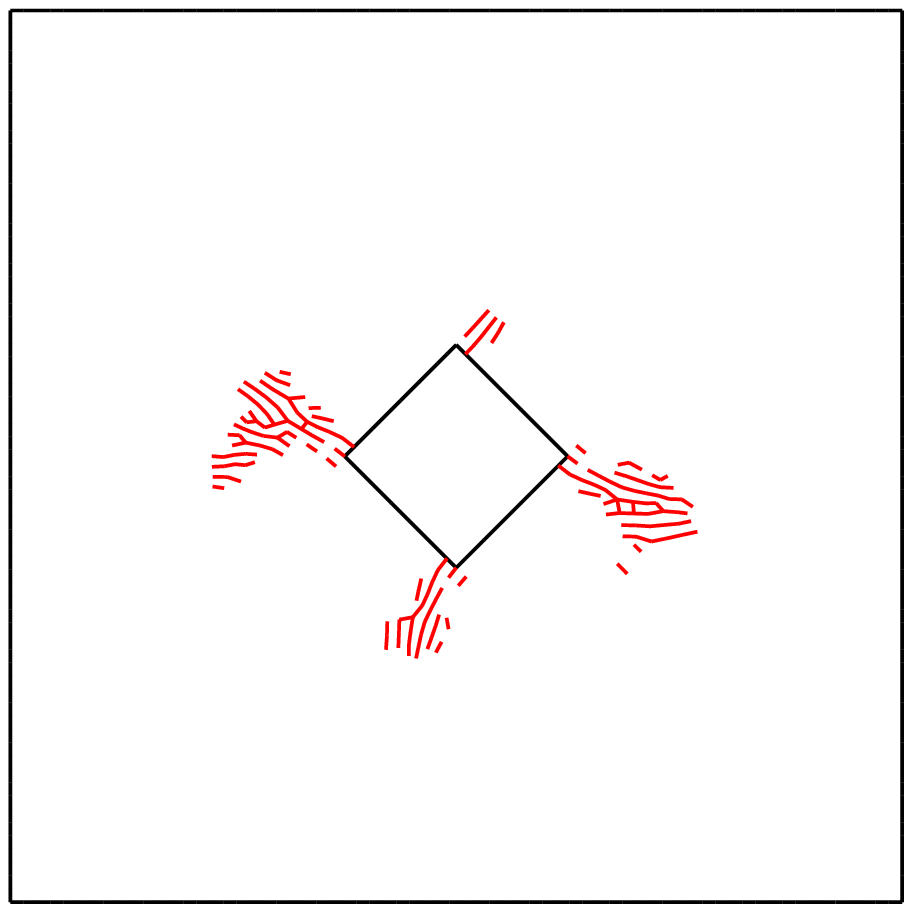}}  
    \subfigure[]{ 
    \label{fig:Ex33_01_02} 
    \includegraphics[width=0.48\textwidth]{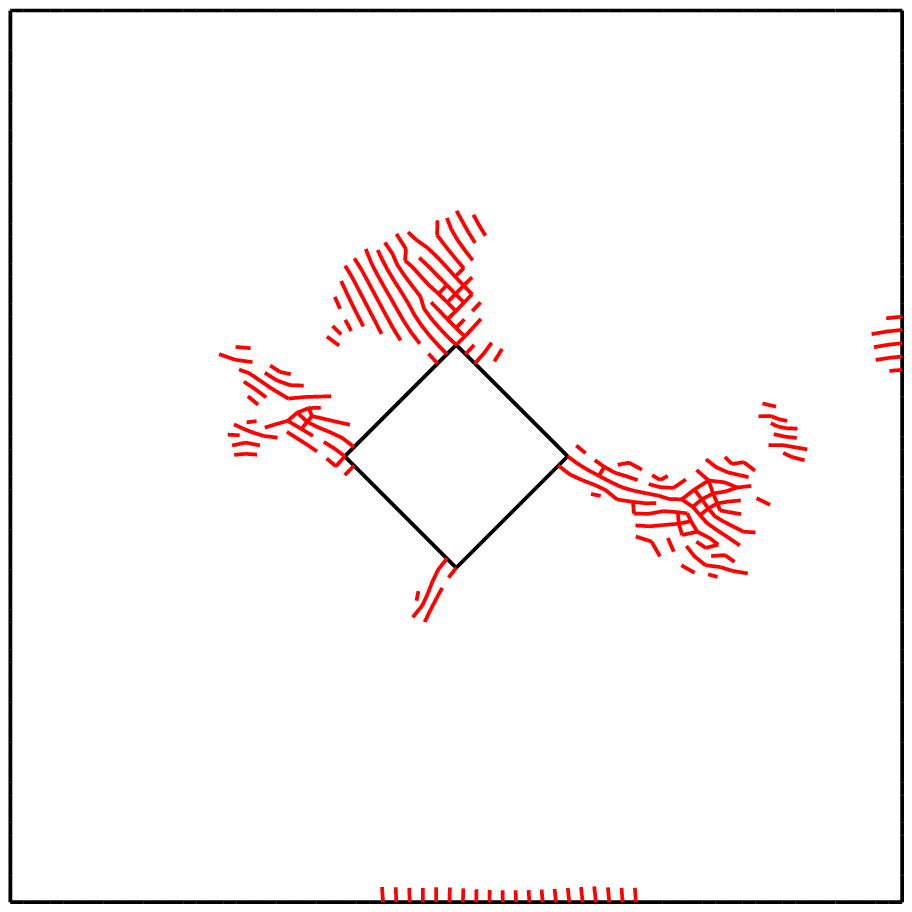}}  
        \subfigure[]{ 
    \label{fig:Ex33_01_03} 
    \includegraphics[width=0.48\textwidth]{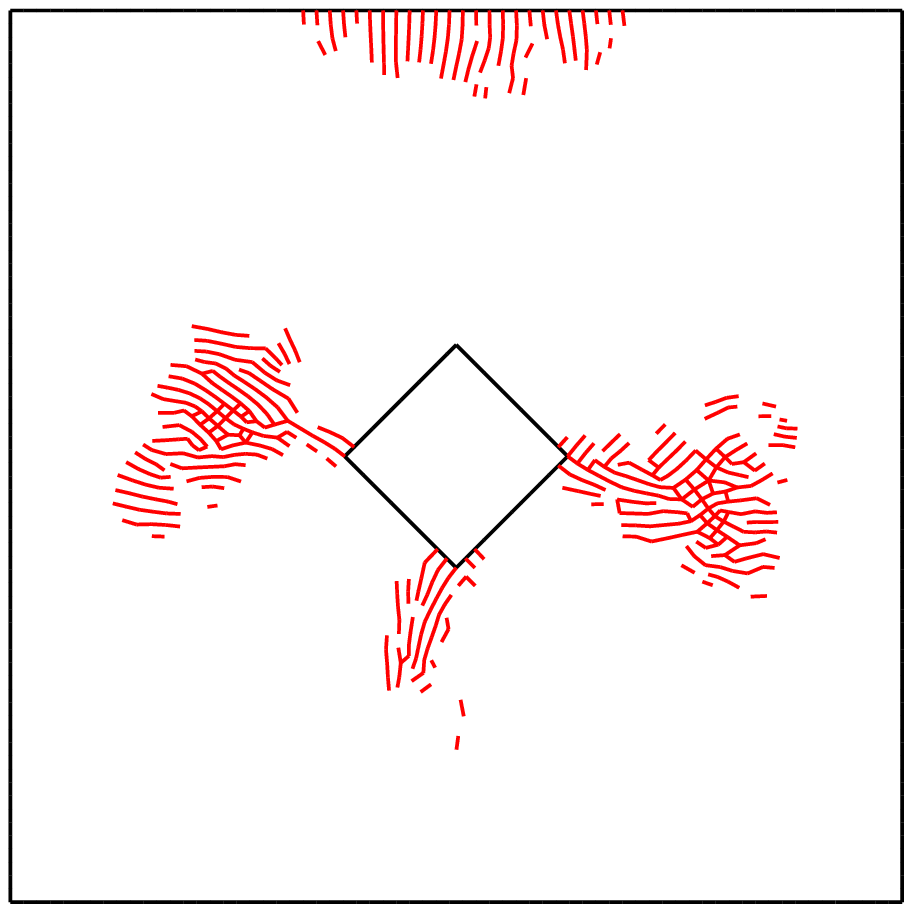}}  
  \caption{The crack development results simulated by FPM with Hoop-Stress-Based criterion for (a) non-piezoelectric material. (b) poling direction of material parallel to $y$-axis. (c) poling direction of material antiparallel to $y$-axis.} 
  \label{fig:Ex33_01} 
\end{figure}

Fig.~\ref{fig:Ex33_Part_01} presents the domain partition transformed from ABAQUS meshing and used in the FPM. The Points are distributed at the centroid of the subdomains. Computational parameters: $c_0 = \sqrt{20}$, $\eta_{21} = 2.0 E$, $\eta_{22} =50 E$, $\eta_{23} = 2.0 \Lambda_{33}$, $\eta_{24} = 0$. First, the Hoop-Stress-Based criterion is employed to simulate the crack initiation and development: in each analysis step, for all the internal boundaries, if the normal stress exceeds a prescribed critical value, the boundary is cracked. The simulated crack development results based on the Hoop-Stress-Based criterion with different poling directions are shown in Fig.~\ref{fig:Ex33_01}. As we have stated, when the poling direction of the material changes, the crack propagation in different directions are either helped or hindered.

\begin{figure}[htbp] 
  \centering 
    \subfigure[]{ 
    \label{fig:Ex33_02_01} 
    \includegraphics[width=0.48\textwidth]{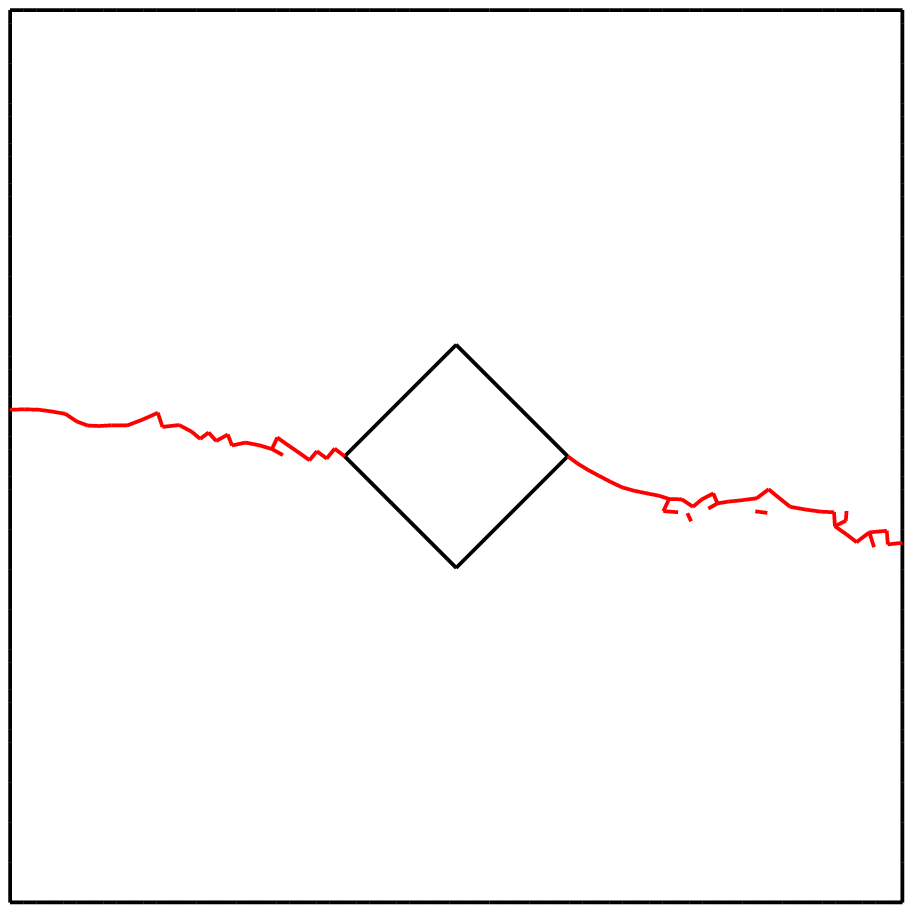}}  
    \subfigure[]{ 
    \label{fig:Ex33_02_02} 
    \includegraphics[width=0.48\textwidth]{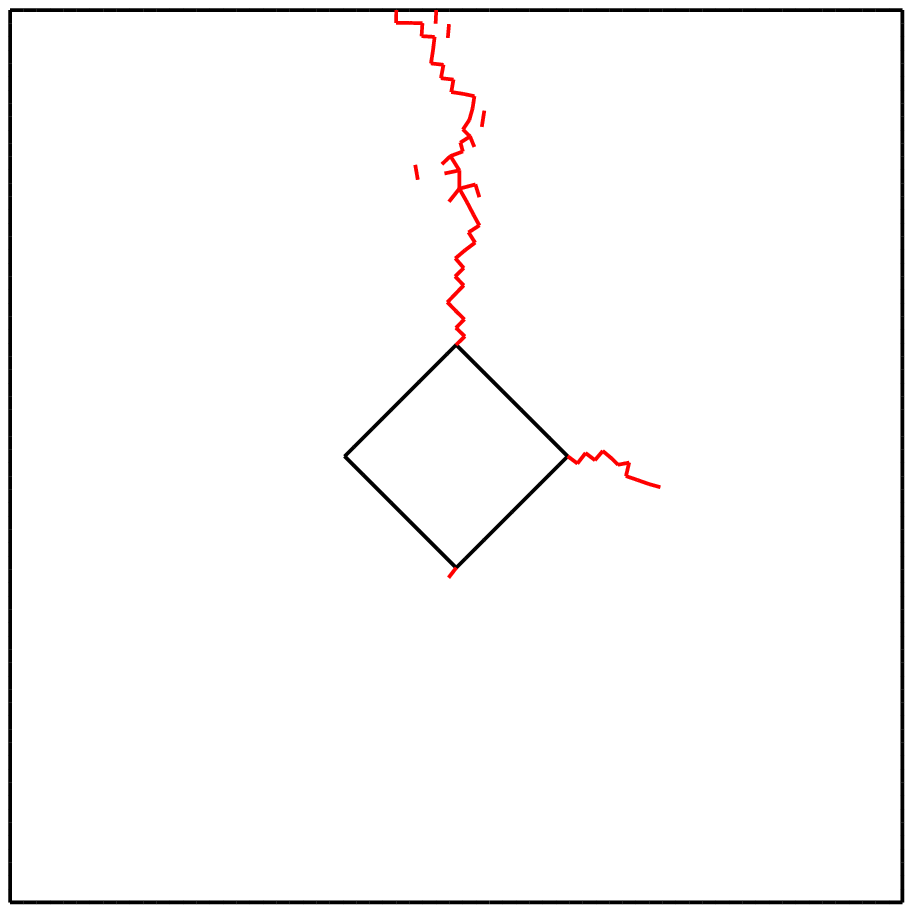}}  
        \subfigure[]{ 
    \label{fig:Ex33_02_03} 
    \includegraphics[width=0.48\textwidth]{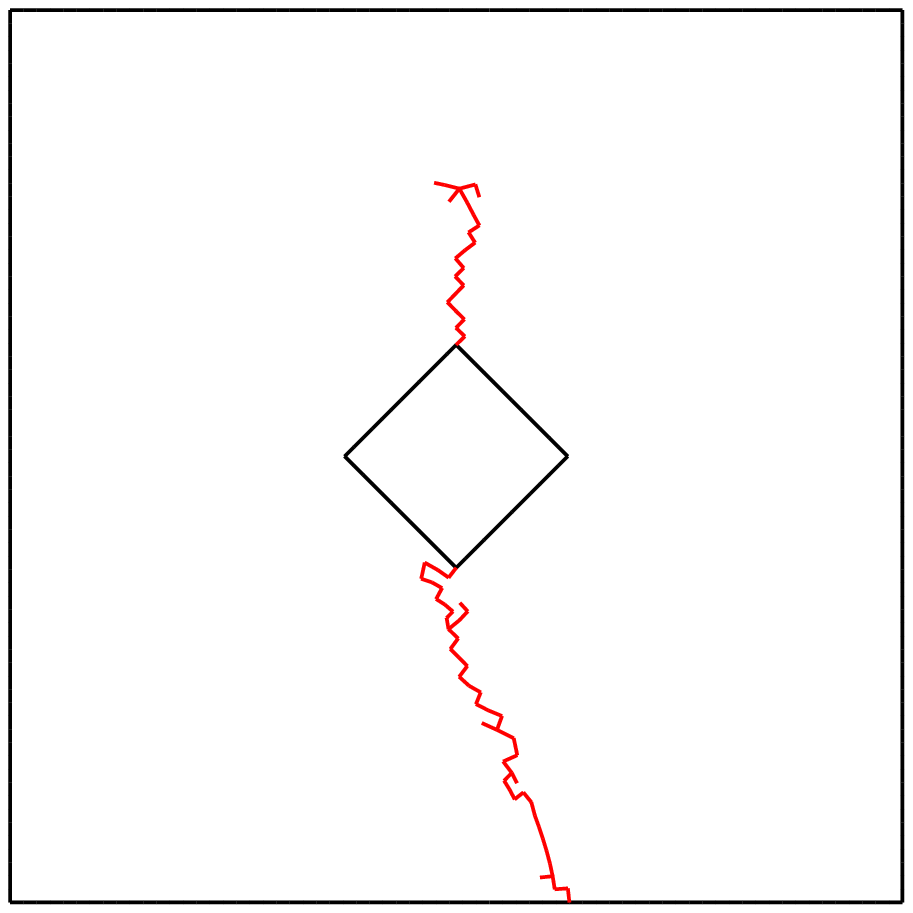}}  
  \caption{The crack development results simulated by FPM with BER-based criterion for (a) non-piezoelectric material. (b) poling direction of material parallel to $y$-axis. (c) poling direction of material antiparallel to $y$-axis.} 
  \label{fig:Ex33_02} 
\end{figure}

Furthermore, we simulate the same example under an inter-subdomain-boundary Bonding-Energy-Rate (BER) based criterion. The definition of BER and relating formulations can be found in Part. I of this series. In each analysis step, if the BER on an internal boundary exceeds a prescribed critical value, the boundary is cracked. Fig.~\ref{fig:Ex33_02} presents the simulated crack development results based on the BER-based criterion. While the same influence of the poling direction is observed, the simulation result based on BER-based criterion shows a more concise crack path comparing with the Hoop-Stress-Based criterion.

As can be seen, the FPM shows great results in simulation crack initiation and propagation problems with multiple crack criteria. Other continuum-physics-based criteria for crack initiation and development can also cooperate with the FPM. Based on discontinues trial and test functions, the process of simulating crack propagations in the FPM is much easier than the traditional element-based methods and other meshfree methods based on continuous trial and test functions.

\section{Discussion of the computational parameters} \label{sec:comp_para}

In this section, the recommended ranges for the computational parameters are given. The recommended values are based on a parametric study. Yet for brevity, we omit the details at here. A similar parametric study and discussion of the computational parameters on the FPM for heat transfer analysis can be seen in \cite{Guan2020a}.

The appropriate computational parameters may vary slightly for different problems, according to the continuity requirements. Table~\ref{table:Range} shows the recommended ranges for the computational parameters, in which $E$ is the Young’s modulus, and $\Lambda = \sqrt{\Lambda_{11} \Lambda_{33}}$ is the average permittivity of the dielectric in both directions. A recommended value is also given for all the computational parameters without any prior knowledge. Note that  very large values can be used for the penalty parameters enforcing the essential boundary conditions ($\eta_{11}$, $\eta_{12}$, $\eta_{13}$ and $\eta_{14}$). However, for the penalty parameters relating to the continuity requirement ($\eta_{21}$, $\eta_{22}$, $\eta_{23}$ and $\eta_{24}$), an excessively large value may be harmful for the accuracy and thus should be avoided. On the other hand, a too small penalty parameter can result in stability problem or discontinuous solutions. Therefore, an appropriate value of the penalty parameters is required. The constant parameter $c_0$ used in the local Radial Basis Function-based Differential Quadrature method within the recommended range has little influence on the final solutions. Yet a relatively small $c_0$ (e.g., $c_0 < 1$) may lead to inaccurate approximation of the derivatives and thus result in inaccurate solutions. Whereas a too large $c_0$ causes a stability problem. Therefore, in order to ensure the stability and consistency of the FPM and achieve a good accuracy, the computational parameters should be selected carefully in the recommended ranges.

\begin{table}[htbp]
\caption{Recommended ranges and values of the computational parameters (if applicable).}
\centering
{
\begin{tabular}{ | c | c | c |  }
\hline 
 Computational parameters & Recommended range & Recommended value \\
\hline
$c_0$ & $\sqrt{1} \sim \sqrt{25} $ & $\sqrt{20}$ \\
\hline
$\eta_{11} / E$ & $ > 10^3 $ & $10^{10}$ \\
\hline
$\eta_{12} / E$ & $ > 10^3 $ & $10^{10}$ \\
\hline
$\eta_{13} / \Lambda $ & $ > 10^3 $ & $10^{10}$ \\
\hline
$\eta_{14} / \Lambda $ & $ > 10^3 $ & $10^{10}$ \\
\hline
$\eta_{21} / E$ & $0 \sim 10$ & $2$ \\
\hline
$\eta_{22} / E$ & $0 \sim 10^2 $ & $50$ \\
\hline
$\eta_{23} / \Lambda $ & $ 0 \sim 10 $ & $0$ \\
\hline
$\eta_{34} / \Lambda $ & $ 0 \sim 10^2 $ & $0$ \\
\hline
\end{tabular}}
\label{table:Range}
\end{table}

\section{Conclusion}

In this paper, as the second part of our current work, extensive numerical results and validation of the FPM in analyzing flexoelectric effects in dielectric solids are given. Numerous examples are presented for both the primal and mixed FPM with full or reduced theories. The present FPMs have passed the linear and quadratic displacement patch tests successfully, and have shown excellent solutions in analyzing flexoelectric and converse flexoelectric effects in continuous domains. Scaling studies have shown that the difference between the full and reduced theories is negligible at micro or larger scales. Whereas for nano-scale structures, the influence of the electroelastic stress, as well as the electric field gradient effects become significant and thus the full theory is recommended to be employed. The FPMs are also applied in analyzing flexoelectric effects on stationary cracks where a concentration of strain occurs. Finally, we present the simulation of crack initiation propagation problems involving flexoelectric effects using the FPM, which has never been studied in any previous literature. Both the simple Hoop-stress-criterion and a BER-based crack criterion are utilized. Numerical solutions have shown that the FPM can successfully predict the crack propagation paths and the results are not sensitive to mesh refinement. The flexoelectric effect, coupling with the piezoelectricity, may help, hinder, or deflect the crack propagation in dielectric materials. A recommended range of the computational parameters is given at last. We can conclude that, both the presently developed primal and mixed FPMs have shown excellent performance in flexoelectric analysis. Finally, the present FPMs are especially suitable for crack initiation and propagation simulations.

\section*{References}
\bibliography{Part_2}

\section*{Acknowledgment}

Yue Guan thankfully acknowledges the financial support for her work, provided through the funding for Professor Atluri’s Presidential Chair at TTU.

\appendix

\section{Flexoelectric material with cubic symmetry} \label{app:sym}

In this appendix, the constitutive matrices are given for a material with isotropic elasticity, cubic symmetry for the flexoelectric tensor, and tetragonal symmetry for the piezoelectric tensor.
\begin{align}
\begin{split}
 & \boldsymbol{\Lambda} = \left[ \begin{matrix} \Lambda_{11} & 0 \\ 0 & \Lambda_{33} \end{matrix} \right], \quad \mathbf{C}_{\sigma \varepsilon} = \left[ \begin{matrix} \lambda & \lambda + 2 G & 0 \\ \lambda + 2 G & \lambda & 0 \\ 0 & 0 & G \end{matrix} \right], \quad \boldsymbol{\Phi} = \left[ \begin{matrix} \Phi_{11} & 0 & 0 & 0 \\ 0 & \Phi_{33} & 0 & 0 \\ 0 & 0 & \Phi_{11} & 0 \\ 0 & 0 & 0 & \Phi_{33} \end{matrix} \right], \\ 
& \mathbf{C}_{\mu \kappa} = l^2 \left[ \begin{matrix} \lambda + 2 G & 0 & 0 & 0 & 0 & \lambda/2 \\ 0 & \lambda + 2 G & 0 & 0 & \lambda/2 & 0 \\ 0 & 0 &  G & 0 & 0 & G/2 \\ 0 & 0 & 0 & G & G/2 & 0 \\ 0 & \lambda/2  & 0 & G/2  & \left( \lambda + 3 G \right) /4 & 0 \\ \lambda/2 & 0 & G/2 & 0 & 0 & \left( \lambda + 3 G \right) /4  \end{matrix} \right], \mathbf{b} = \mathbf{0},\\ 
& \mathbf{e} = \left[ \begin{matrix} 0 & 0 & e_{15} \\ e_{31} & e_{33} & 0 \end{matrix} \right]^\mathrm{T}, \quad \mathbf{a} = \left[ \begin{matrix}  \overline{\mu}_{11} & 0 & \overline{\mu}_{44} & 0 & 0 & \left( \overline{\mu}_{12} + \overline{\mu}_{44} \right) /2 \\ 0 & \overline{\mu}_{11} & 0 & \overline{\mu}_{44} & \left( \overline{\mu}_{12}  + \overline{\mu}_{44}  \right) / 2 & 0  \end{matrix} \right]^\mathrm{T}.
\end{split}
\end{align}
where $\left( \lambda , G \right)$ are the Lamé parameters, $l$ is the internal material length, $\boldsymbol{\Lambda}$, $\mathbf{e}$, and $\mathbf{a}$ are related to the permittivity of the dielectric, the piezoelectric tensor and the flexoelectric tensor. Matrices $\mathbf{b}$ and $\boldsymbol{\Phi}$ donates the quadrupole-strain coefficients and higher-order electric parameters respectively.

\section{Material Properties for PZT-5H} \label{app:PZT-5H}

\begin{align}
\begin{split}
 & \boldsymbol{\Lambda} = \left[ \begin{matrix} \Lambda_{11} & 0 \\ 0 & \Lambda_{33} \end{matrix} \right], \quad \boldsymbol{\Phi} = q^2 \left[ \begin{matrix} \boldsymbol{\Lambda} & \mathbf{0} \\ \mathbf{0} & \boldsymbol{\Lambda}   \end{matrix} \right], \quad \mathbf{e} = \left[ \begin{matrix} 0 & 0 & e_{15} \\ e_{31} & e_{33} & 0 \end{matrix} \right]^\mathrm{T}, \quad \mathbf{b} = b^2 \left[ \begin{matrix} \mathbf{e} & \mathbf{e} \end{matrix} \right], \\ 
& \mathbf{C}_{\sigma \varepsilon} = \left[ \begin{matrix} c_{11} & c_{13} & 0 \\ c_{13} & c_{33} & 0 \\ 0 & 0 & c_{44} \end{matrix} \right], \quad  \mathbf{a} = a^2 \left[ \begin{matrix} 0 & e_{31} \\ 0 & e_{33} \\ 0 & e_{33} \\ 0 & e_{31} \\ e_{15} & 0 \\ e_{15} & 0 \end{matrix} \right], \quad \mathbf{C}_{\mu \kappa} = l^2 \left[ \begin{matrix} c_{11} & 0 & c_{13} & 0 & 0 & 0 \\ 0 & c_{33} & 0 & c_{13} & 0 & 0 \\ c_{13} & 0 & c_{33} & 0 & 0 & 0 \\ 0 & c_{13} & 0 & c_{11} & 0 & 0 \\ 0 & 0 & 0 & 0 & c_{44} & 0 \\ 0 & 0 & 0 & 0 & 0 & c_{44} \end{matrix} \right].
\end{split}
\end{align}
\begin{align} \nonumber
\text{where :} \; & c_{11} = 12.6 \times 10^{10}~\mathrm{Pa}, \;  c_{13} = 5.3 \times 10^{10}~\mathrm{Pa}, \; c_{33} = 11.7 \times 10^{10}~\mathrm{Pa}, \; c_{44} = 3.53 \times 10^{10}~\mathrm{Pa}; \\ \nonumber
& e_{31} = -6.5~\mathrm{C/m^2}, \;  e_{33} = 23.3~\mathrm{C/m^2}, \; e_{15} = 17.0~\mathrm{C/m^2}; \\ \nonumber
& \Lambda_{11} = 15.1 \times 10^{-9}~\mathrm{F/m}, \; \Lambda_{33} = 13.0 \times 10^{-9}~\mathrm{F/m}.
\end{align}
The size factors $l$, $q$, $a$ and $b$ are:
\begin{align} \nonumber
\begin{split}
l^2 = \alpha l_0^2, \; \; q^2 = \alpha q_0^2, \; \; a^2 = \alpha a_0^2, \; \;, b^2 = \alpha b_0^2, 
\end{split}
\end{align}
where $l_0 = 5 \times 10^{-9}~\mathrm{m}$,  $q_0 = 3 \times 10^{-10}~\mathrm{m}$,  $a_0^2 = 5 \times 10^{-12}~\mathrm{m}$ and  $b_0^2 = 5 \times 10^{-12}~\mathrm{m}$.

\end{document}